\documentclass[11pt]{amsart}

\usepackage{amsmath,amssymb}
\usepackage{mathtools}
\usepackage{mathrsfs}
\usepackage{eucal}
\usepackage{tikz-cd}
\usetikzlibrary{arrows.meta,decorations.markings,calc}
\usepackage{microtype}
\usepackage[english]{babel}
\usepackage[T1]{fontenc}
\usepackage{lmodern}
\usepackage[utf8]{inputenc}

\usepackage[
  top=3cm,
  bottom=2cm,
  left=3cm,
  right=3cm,
  marginparwidth=1.75cm
]{geometry}
\usepackage{graphicx}
\usepackage{xcolor}
\graphicspath{{pictures/}}

\usepackage{enumitem}
\setlist[1]{itemsep=3pt}

\definecolor{linkblue}{RGB}{25,65,110}
\colorlet{cycleblue}{linkblue}
\definecolor{cycleteal}{RGB}{30,114,106}
\definecolor{cyclerust}{RGB}{166,78,54}
\definecolor{cyclegray}{RGB}{103,111,119}
\usepackage[colorlinks=true,allcolors=linkblue]{hyperref}
\usepackage[nameinlink]{cleveref}

\crefname{theorem}{Theorem}{Theorems}
\Crefname{theorem}{Theorem}{Theorems}
\crefname{lemma}{Lemma}{Lemmata}
\Crefname{lemma}{Lemma}{Lemmata}
\crefname{proposition}{Proposition}{Propositions}
\Crefname{proposition}{Proposition}{Propositions}
\crefname{corollary}{Corollary}{Corollaries}
\Crefname{corollary}{Corollary}{Corollaries}
\crefname{definition}{Definition}{Definitions}
\Crefname{definition}{Definition}{Definitions}
\crefname{remark}{Remark}{Remarks}
\Crefname{remark}{Remark}{Remarks}
\crefname{example}{Example}{Examples}
\Crefname{example}{Example}{Examples}
\crefname{section}{Section}{Sections}
\Crefname{section}{Section}{Sections}
\crefname{figure}{Figure}{Figures}
\Crefname{figure}{Figure}{Figures}
\crefname{equation}{}{}
\Crefname{equation}{}{}

\theoremstyle{plain}
\newtheorem{theorem}{Theorem}[section]

\newtheorem{proposition}[theorem]{Proposition}
\newtheorem{corollary}[theorem]{Corollary}
\newtheorem{lemma}[theorem]{Lemma}

\newtheorem*{theorem*}{Theorem}
\newtheorem{introtheorem}{Theorem}

\crefname{introtheorem}{Theorem}{Theorems}
\Crefname{introtheorem}{Theorem}{Theorems}

\theoremstyle{definition}
\newtheorem{definition}[theorem]{Definition}

\theoremstyle{remark}
\newtheorem{remark}[theorem]{Remark}

\newtheorem*{introexample}{Example}
\newtheorem{example}[theorem]{Example}

\newcommand{\R}{\mathbb{R}}
\newcommand{\C}{\mathbb{C}}
\newcommand{\N}{\mathbb{N}}

\newcommand{\llbracket}{\mathopen{[\![}}
\newcommand{\rrbracket}{\mathclose{]\!]}}

\newcommand{\CP}{\C\mathrm{P}}
\newcommand{\PP}{\mathbb{P}}

\newcommand{\dd}{\mathrm{d}}

\newcommand{\dM}{\mathsf d_M}
\newcommand{\W}{\mathcal{W}}

\newcommand{\cC}{\mathcal{C}}
\newcommand{\cP}{\mathcal{P}}
\newcommand{\mug}{\mu}

\newcommand{\owp}{\Omega}

\DeclareMathOperator{\supp}{supp}
\DeclareMathOperator{\vol}{vol}

\title{Wasserstein geometry of analytic cycle spaces}
\author{Fabio Cavalletti}
\author{Lorenzo Cecchi}
\author{Antonio Lerario}
\begin{document}

\begin{abstract}
We use optimal transport to define an intrinsic distance on connected
components of the Barlet space of analytic cycles on a compact K\"ahler
manifold, by associating to each cycle its normalized volume measure.
We prove that any two cycles in the same connected component can be joined by a deformation of finite Wasserstein length, and that the resulting metric is complete and geodesic. On the locus of smooth reduced cycles at regular
parameter points, Wasserstein lengths and energies coincide with those
of the normalized Weil--Petersson metric. The construction thus gives
a transport interpretation of Weil--Petersson geometry and extends
beyond it to components containing no smooth reduced cycles.
\end{abstract}
\maketitle
\setcounter{tocdepth}{1}
\tableofcontents

\section{Introduction}\label{sec:introduction}

What is the optimal way to deform a complex projective hypersurface into another one? In \cite{ACL}, Antonini, Cavalletti and Lerario approached this question using optimal transport. The starting point was to associate to a hypersurface its normalized volume measure, obtaining a topological embedding of the space of degree-$d$ hypersurfaces in $\CP^n$ into the Wasserstein space of probability measures on $\CP^n$. Singular hypersurfaces were included by integrating over their regular parts and keeping track of multiplicities. One could then measure a deformation by the Wasserstein length of the corresponding curve of measures, and define a distance by minimizing this length among deformations through hypersurfaces of the same degree. A main result of that work was that this distance is finite, complete and geodesic. Moreover, on the open set of smooth hypersurfaces, the construction gives a natural K\"ahler structure.

This question is part of a more general problem: how should we measure distances in a parameter space of geometric objects? For hypersurfaces, the parameter space is itself a complex projective space, namely the projectivization of the space of homogeneous polynomials of degree $d$. A first possibility is therefore to use the Fubini--Study distance associated with a Hermitian inner product on the coefficient space. This geometry is homogeneous: it does not distinguish points representing smooth hypersurfaces from points representing singular ones, and it does not directly measure the motion of their zero sets. One way to incorporate this distinction is to make motion more expensive near the discriminant. This is the principle behind condition geometry in polynomial system solving, developed from the work of Shub and Smale; see \cite{shub-condition,shub-smale-geometric}. For example, when the discriminant $\Sigma$ of singular hypersurfaces is nonempty, one can rescale the Fubini--Study Riemannian metric by the inverse square of the distance to $\Sigma$. The resulting metric has locally Lipschitz coefficients, and approaching $\Sigma$ requires infinite length. (The usual condition metric is defined on a space of problem--solution pairs, with the condition number as a weight.) For the geometric deformation problem, however, we also want to include singular objects and allow paths to reach them. The Wasserstein construction from \cite{ACL} does this while measuring the motion of the hypersurfaces themselves.

In this paper we extend the construction of \cite{ACL} from projective
hypersurfaces to analytic cycles on compact K\"ahler manifolds,
developing a \emph{Wasserstein geometry of analytic cycle spaces}.
Barlet's cycle space provides the parameter space for these objects
(\Cref{sec:cycle-space}). We represent each cycle by its normalized
volume measure and measure a deformation by the Wasserstein length of
the resulting curve of measures. Minimizing these lengths gives a
finite intrinsic distance on each connected component of the Barlet space. On the locus of
smooth reduced cycles at regular parameter points, the corresponding
lengths and energies agree with those of the classical Weil--Petersson
metric. The distance is also defined on components containing no
smooth reduced cycles, where there is no such smooth geometry from
which to start. We explain now the main ideas of the construction.

Fix a compact connected K\"ahler manifold $(M,\omega)$ of complex dimension $n$ and an integer $0\leq q\leq n$. We denote the associated Riemannian metric by $g$ and its distance by $\dM$. An \emph{effective analytic $q$-cycle} is a finite formal sum
\[
 \gamma=\sum_{j=1}^N m_j[Z_j],
\]
where the $Z_j\subset M$ are distinct irreducible analytic subsets of complex dimension $q$, the positive integers $m_j$ are their multiplicities, and $[Z_j]$ denotes the cycle with support $Z_j$ and multiplicity one. Thus, for example, a pair of distinct projective lines gives the cycle $[L_1]+[L_2]$, whereas a double line gives $2[L]$. A cycle is called \emph{reduced} if all its multiplicities are one. Integration over $\gamma$ means integration over the smooth parts $(Z_j)_{\mathrm{reg}}$, weighted by $m_j$.

We denote by $\mathcal C_q(M)$ the Barlet space of these cycles. It carries a natural reduced complex-space structure, locally modeled on complex analytic sets, and its connected components are compact. If $M$ is projective, then, for a fixed projective embedding, the cycles of a given degree form the analytification of the corresponding reduced Chow variety (\Cref{remark:projective}). Thus our construction also equips connected components of Chow varieties with an intrinsic transport geometry. The homology class of a cycle is constant on each component (although distinct components may have the same homology class). In particular, on a connected component $B\subset\mathcal C_q(M)$ other than that of the zero cycle, the volume $V_B:=\int_\gamma\omega^q/q!$ is positive and independent of $\gamma\in B$. This follows from the closedness of $\omega$; the relevant facts about cycles, their topology and integration are collected in \Cref{sec:cycle-space}.

Let $\mathcal P(M)$ denote the space of Borel probability measures on $M$. The central construction is the \emph{cycle-to-measure} map
\[
 \mu:B\longrightarrow\mathcal P(M),\qquad
 \gamma=\sum_jm_j[Z_j]\longmapsto
 \frac1{V_B}\sum_jm_j\operatorname{vol}_{(Z_j)_{\mathrm{reg}}},
\]
where $\operatorname{vol}_{(Z_j)_{\mathrm{reg}}}$ is the Riemannian volume measure induced by $g$. The mass assigned to a region of $M$ is the fraction of the total cycle volume contained in that region, with multiplicities included. Thus integration against $\mu(\gamma)$ is normalized integration along the cycle (\Cref{def:mu-gamma}). We equip $\mathcal P(M)$ with the $2$-Wasserstein distance, denoted by $\W$, which is the square root of the least quadratic transport cost, with cost $\dM(x,y)^2$ for moving a unit of mass from $x$ to $y$ (\Cref{def:wasserstein-distance}). Since $M$ is compact, this distance induces the narrow topology, namely convergence against continuous test functions.

The first observation is that the map $\mu$ is a topological embedding, see \Cref{prop:mu-gamma}. This will allow us to identify the component $B$ with a subset of the Wasserstein space. Indeed, integration of continuous $(q,q)$-forms varies continuously in the cycle-space topology, so the preceding formula gives continuity of $\mu$. Together with the requirement that cycles avoiding a fixed closed subset form an open set, these integration maps characterize the topology; see \Cref{sec:cycle-space}. The measure determines the support and, since $V_B$ is fixed, the multiplicities of the cycle. Thus $\mu$ is injective, and compactness of $B$ gives the embedding.

To relate the measure construction to classical geometry, we next recall the K\"ahler structures on the cycle space. Results of Fujiki and Varouchas show that it is K\"ahler in the sense of complex spaces; see \cite[Chapter~12]{cycles2} and \cite{axelsson-schumacher}. This notion is formulated using local strictly plurisubharmonic potentials, so it makes sense even at singular parameter points. We will use its smooth counterpart on a suitable open locus; \Cref{remark:current} explains the distinction between the K\"ahler structure on the complex space and the smooth form considered here.

Fix an irreducible component $C\subset B$, with its reduced structure. In this paper, a point $\gamma\in C$ is called \emph{good} if it is a smooth point of the parameter space $C$ and represents a reduced cycle $\gamma=[Z]$ whose support $Z$ is a smooth compact complex submanifold of $M$, possibly disconnected (\Cref{def:two-regularities}). One condition concerns the parameter space, the other the cycle it represents. The good locus, denoted by $\widetilde C$, is the complement of a closed analytic subset of $C$. If nonempty, it is a dense connected complex manifold, and its complement has real codimension at least two (\Cref{prop:bad-locus-analytic}). Note that a connected component $B$ may also contain irreducible components with no good points.

To describe the Weil--Petersson geometry on $\widetilde C$, we first need to understand its tangent vectors. At a good point $\gamma=[Z]$, let $N_{Z/M}$ denote the holomorphic normal bundle (ambient tangent directions modulo those tangent to $Z$) and $H^0(Z,N_{Z/M})$ its space of global holomorphic sections. Classical deformation theory associates with each tangent vector $v\in T_\gamma C$ a section $\nu_v$, and this correspondence is injective (\Cref{Ss:tangentgood}). To see this geometrically, represent the deformation by a family of embeddings of $Z$ into $M$. Changing the parametrization of $Z$ changes the velocity by a tangential field. Its normal component is therefore independent of that choice and, using the metric $g$ to represent normal classes by orthogonal normal vectors, gives $\nu_v$, see \Cref{lemma:normalnu}. (A holomorphic normal section represents an infinitesimal deformation, but it need not extend to an actual family. The injective correspondence above is therefore not always surjective; see \Cref{rem:douady} for the deformation-theoretic details.)

We can now average the ambient form $\omega$ over $Z$. For $v,w\in T_\gamma C$, set
\[
 \Omega_\gamma(v,w)
 :=\frac1{V_B}\int_Z\omega(\nu_v,\nu_w)\,\frac{\omega^q}{q!}.
\]
This is the normalized Weil--Petersson K\"ahler form; see \cite{biswas-schumacher,axelsson-schumacher} and \Cref{def:WP}. If $J$ denotes the complex structure of $\widetilde C$, the associated Riemannian metric is $\Theta(v,w)=\Omega(v,Jw)$: the squared norm of a deformation is its mean squared normal velocity. Equivalently, $\Omega$ is the normalized fiber integral of the pullback of $\omega^{q+1}/(q+1)!$ to the \hyperref[sec:tautological]{tautological family} (whose fiber over $[Z]$ is $Z$); see \Cref{sec:WPgood}. Its Riemannian distance on $\widetilde C$, denoted by $\mathsf d_{\mathrm{WP}}$, need not be complete.

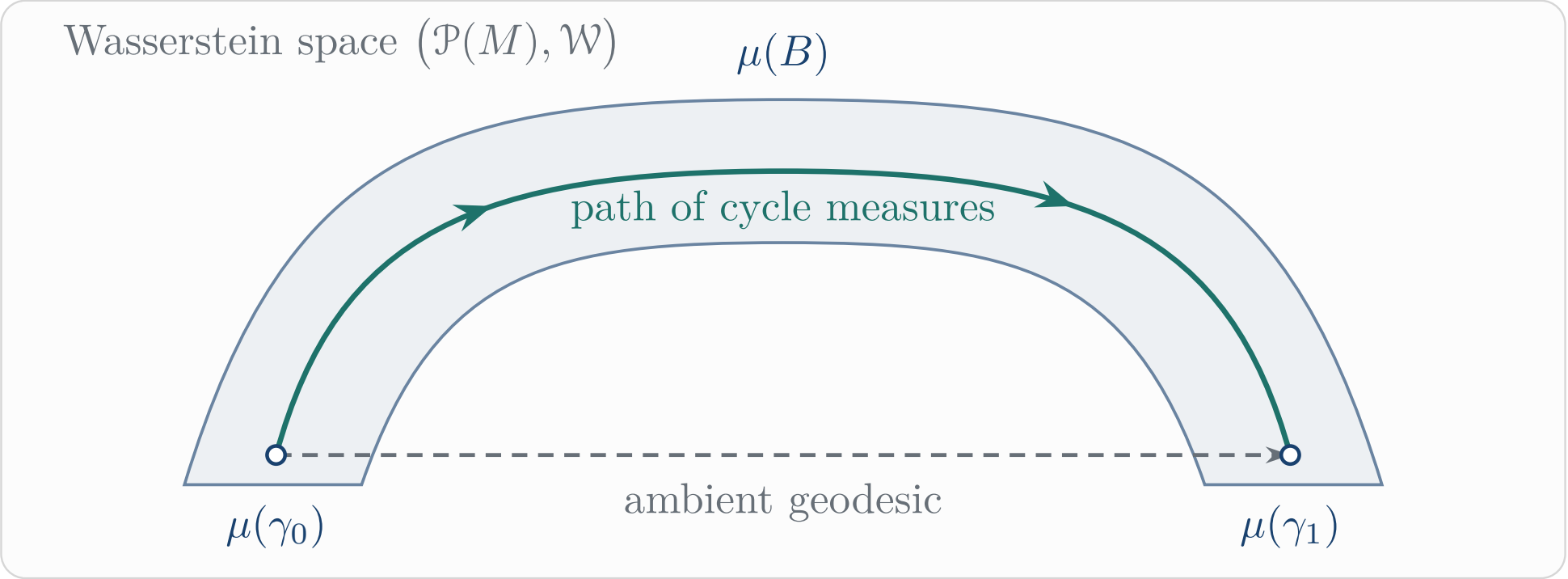
\begin{figure}[!t]
\centering
\begin{tikzpicture}[x=1cm,y=.9cm,font=\small,
  every node/.style={align=center},
  >={Stealth[length=2mm,width=1.35mm]}]
  \filldraw[fill=black!1,draw=black!16,line width=.45pt,rounded corners=2mm]
    (-6.4,-1.0) rectangle (6.4,4.25);
  \node[anchor=west,text=cyclegray] at (-6.03,3.86)
    {Wasserstein space $\bigl(\mathcal P(M),\W\bigr)$};

  \path[fill=cycleblue!8,draw=cycleblue!65,line width=.7pt]
    (-4.9,-.15)
    .. controls (-4.05,2.95) and (-2.75,3.35) .. (0,3.35)
    .. controls (2.75,3.35) and (4.05,2.95) .. (4.9,-.15)
    -- (3.45,-.15)
    .. controls (2.83,1.85) and (1.85,2.05) .. (0,2.05)
    .. controls (-1.85,2.05) and (-2.83,1.85) .. (-3.45,-.15)
    -- cycle;
  \node[text=cycleblue] at (0,3.76) {$\mu(B)$};

  \coordinate (leftmeasure) at (-4.15,.12);
  \coordinate (rightmeasure) at (4.15,.12);
  \draw[cyclegray,line width=.9pt,dash pattern=on 3.6pt off 2.5pt,->]
    (leftmeasure) -- (rightmeasure);
  \node[text=cyclegray,fill=black!1,inner sep=2pt] at (0,-.31)
    {ambient geodesic};

  \draw[cycleteal,line width=1.25pt,
    postaction={decorate},
    decoration={markings,
      mark=at position .27 with {\arrow{Stealth}},
      mark=at position .73 with {\arrow{Stealth}}}]
    (leftmeasure)
    .. controls (-3.65,2.15) and (-2.5,2.7) .. (0,2.7)
    .. controls (2.5,2.7) and (3.65,2.15) .. (rightmeasure);
  \node[text=cycleteal] at (0,2.35) {path of cycle measures};

  \filldraw[fill=white,draw=cycleblue,line width=.85pt]
    (leftmeasure) circle[radius=2.1pt];
  \filldraw[fill=white,draw=cycleblue,line width=.85pt]
    (rightmeasure) circle[radius=2.1pt];
  \node[text=cycleblue] at (-4.15,-.53) {$\mu(\gamma_0)$};
  \node[text=cycleblue] at (4.15,-.53) {$\mu(\gamma_1)$};
\end{tikzpicture}
\caption[The path constraint defining the inner Wasserstein distance]{%
The inner distance constrains the whole deformation (schematic view).
The dashed path illustrates an ambient Wasserstein geodesic that may
leave $\mu(B)$. Every point of the solid path represents a cycle in $B$.
The inner Wasserstein distance is the infimum of the Wasserstein lengths
of such constrained paths; the ambient distance is a lower bound.}
\label{fig:inner-wasserstein}
\end{figure}

We now return to the problem of measuring deformations on the whole component $B$, seen via the embedding $\mu:B\hookrightarrow \mathcal P(M)$. A path of cycles gives a curve of probability measures through $\mu$, and we use its Wasserstein length, as explained below, as the cost of the deformation. For a curve of measures, we denote this length by $L_\W$: it is the supremum, over finite partitions of the parameter interval, of sums of successive Wasserstein distances (\Cref{sec:length-spaces}). Our aim is to minimize this cost among paths that remain in $B$ throughout the deformation.

The constraint on the path is essential. Simply restricting $\W$ to the cycle measures gives the \emph{outer Wasserstein distance}, which measures the ambient transport cost between the endpoints and supplies a lower bound for our problem. The shortest ambient paths, however, may leave the image of the cycle space. We impose the cycle constraint by taking the infimum of lengths over paths of cycles themselves. This defines the \emph{inner Wasserstein distance}, the distance studied in this paper:
\[
 \W_B^{\mathrm{in}}(\gamma_0,\gamma_1)
 :=\inf_\gamma L_\W(\mu\circ\gamma),
\]
where $\gamma:[0,1]\to B$ ranges over continuous paths joining $\gamma_0$ to $\gamma_1$ (\Cref{def:inner-wasserstein-distance}). Only rectifiable paths, those with finite length, contribute a finite value to this infimum. The difficulty is to show that such a path exists between every pair of cycles. The path constraint is illustrated in \Cref{fig:inner-wasserstein}.

\begin{samepage}
\begin{introtheorem}\label{intro:theorem-A}
The inner Wasserstein distance is finite on $B$ and makes $B$ a complete geodesic metric space.
\end{introtheorem}

This is proved in \Cref{thm:finiteness,thm:inner-wasserstein-complete-geodesic}, in the more general form where paths are constrained to any connected closed analytic subset of $B$. The cycles may be singular or nonreduced throughout the deformation. Thus deformation cost is defined through the motion of measures even where no good-locus geometry is available.
\end{samepage}

We now explain the connection with Weil--Petersson geometry recalled above. An absolutely continuous curve $t\mapsto\mu_t$ in $(\mathcal P(M),\W)$ has a \emph{metric speed}, denoted by $|\dot\mu_t|_\W$. This is the metric analogue of the norm of a velocity vector: it is defined almost everywhere by taking the limit of the distance between nearby measures divided by the elapsed time (\Cref{sec:length-spaces}, \eqref{eq:metricspeed}). Its time integral gives the length of the curve. Optimal transport describes this motion by a time-dependent vector field $v_t$ satisfying the \emph{continuity equation}
\[
 \partial_t\mu_t+\operatorname{div}(v_t\mu_t)=0
\]
in the sense of distributions. Testing this equation against a smooth function expresses the rate of change of its integral as the integral of its directional derivative along $v_t$. This is conservation of mass expressed in terms of the velocity $v_t$ of the moving measure. Under the integrability assumptions in \Cref{lem:continuity-speed-criterion}, the continuity equation implies absolute continuity and bounds the metric speed by the $L^2(\mu_t)$ norm of $v_t$. Conversely, an absolutely continuous curve admits a velocity field of minimal $L^2$ norm, and this norm equals its metric speed. We recall this correspondence in \Cref{sec:wasserstein-tangents}, following \cite[Chapter~8]{ags}.

For a fixed probability measure $\mu$, let $L^2(\mu;TM)$ be the real Hilbert space of ambient vector fields that are square-integrable with respect to $\mu$ and $g$, identifying fields that agree $\mu$-almost everywhere. The \emph{Wasserstein tangent space} at $\mu$ is
\[
 T_\mu\mathcal P(M)
 :=\overline{\{\nabla f:f\in C^\infty(M,\R)\}}^{\,L^2(\mu;TM)},
\]
where $\nabla f$ is the Riemannian gradient with respect to $g$. This definition singles out the velocities that change the measure to first order. Indeed, a field $w$ is orthogonal to every gradient precisely when $\operatorname{div}(w\mu)=0$ in the distributional sense. Adding such a field to a velocity leaves the continuity equation unchanged; orthogonal projection onto $T_\mu\mathcal P(M)$ removes this freedom and gives the velocity of minimal norm. Thus the distinguished velocity of an absolutely continuous curve belongs to $T_{\mu_t}\mathcal P(M)$ at almost every time. It is this linear tangent space that we use throughout.

For a good cycle $\gamma=[Z]$ and its measure $\mu(\gamma)$, this tangent space has a concrete description: it is the orthogonal sum of the closure of tangential gradients on $Z$ and all square-integrable normal fields along $Z$ (see \eqref{eq:tangentsplitting}). The first summand describes redistribution of mass along the fixed support; the second describes its normal motion. Every smooth normal field is, in fact, the restriction to $Z$ of the gradient of a smooth ambient function vanishing on $Z$. The normal deformation fields of classical deformation theory therefore give the inclusions
\[
 T_\gamma C\hookrightarrow H^0(Z,N_{Z/M})
 \hookrightarrow T_{\mu(\gamma)}\mathcal P(M),
 \qquad v\longmapsto\nu_v.
\]
The second inclusion is understood over the real numbers, using the metric identification with the orthogonal normal bundle. In this way, deformations of cycles determine Wasserstein tangent vectors; their compatibility with the complex and symplectic structures is discussed in \Cref{rem:infinitesimal-wasserstein-kahler}.

The key point is that normal deformation fields transport the cycle measures themselves. Along a smooth curve of good cycles, the normal field solves the continuity equation. Since it belongs to the Wasserstein tangent space, it is the velocity of minimal norm. Its $L^2$ norm is exactly the Weil--Petersson norm, which gives our second main result.

\begin{samepage}
\begin{introtheorem}\label{intro:theorem-B}
On the good locus, motion of cycles has the same speed in Wasserstein and Weil--Petersson geometry. More precisely, let $t\mapsto\gamma_t$ be a smooth curve in $\widetilde C$ and put $\mu_t:=\mu(\gamma_t)$. Then
\[
 |\dot\mu_t|_\W^2
 =\|\nu_{\dot\gamma_t}\|_{L^2(\mu_t)}^2
 =\Theta_{\gamma_t}(\dot\gamma_t,\dot\gamma_t).
\]
\end{introtheorem}

This is proved in \Cref{thm:isometry}. The identity holds at every time of a smooth curve, with the one-sided interpretation at endpoints. Its three terms express the same quantity as the squared metric speed of the measure, the mean squared normal velocity of the cycle, and the squared Riemannian speed in the parameter space. Integrating in time identifies Wasserstein and Weil--Petersson energies; integrating their square roots identifies lengths.
\end{samepage}

\Cref{intro:theorem-B} identifies classical deformation energy with the transport energy of the associated volume measures. It also explains the fiber-integral formula for the Wasserstein Hermitian structure obtained in \cite{ACL}: the pullback and pushforward operations appearing there are precisely those defining the Weil--Petersson form (\Cref{ex:hypersurfaces-cpn}).

To compare distances, we must also specify which deformations are allowed. If all paths are required to remain in $\widetilde C$, minimizing their Wasserstein lengths recovers exactly $\mathsf d_{\mathrm{WP}}$. Allowing them to use the whole irreducible component $C$ gives a possibly smaller distance. Nevertheless, sufficiently close good cycles have the same distance in the two geometries: passing through the singular or nonreduced cycles of $C$ gives no shorter connection between them (\Cref{prop:local-wp-inner-wasserstein}). Thus the intrinsic transport geometry retains the classical Riemannian geometry locally, while allowing additional deformations on the full component.

The case of zero-dimensional cycles gives a particularly transparent example. Fix a positive integer $d$. The component of degree-$d$ zero-cycles is the symmetric product $B=\operatorname{Sym}^d(M)$, whose points are unordered configurations of $d$ points, with repetitions allowed (\Cref{sec:zerocycles}). A cycle $\gamma=[x_1]+\cdots+[x_d]$ is represented by the probability measure that assigns mass $1/d$ to each point, counted with multiplicity. For another cycle $\eta=[y_1]+\cdots+[y_d]$, transport is described by an optimal matching. With $\mathfrak S_d$ denoting the permutation group on $d$ elements, we have
\[
 \W\bigl(\mu(\gamma),\mu(\eta)\bigr)^2
 =\frac1d\min_{\sigma\in\mathfrak S_d}
       \sum_{i=1}^d\dM(x_i,y_{\sigma(i)})^2.
\]
\begin{samepage}
\begin{introtheorem}\label{intro:theorem-C}
For degree-$d$ zero-cycles, the inner and outer Wasserstein distances coincide. If the good locus is nonempty, its Weil--Petersson completion is canonically isometric to the full symmetric product equipped with this distance.
\end{introtheorem}

These statements are proved in \Cref{thm:zero-cycles-inner-equals-outer,thm:zero-cycles-completion}. An optimal matching realizes the first assertion by moving paired points along minimizing ambient geodesics. When the endpoint configurations have distinct points, these motions have no collisions at intermediate times, which also yields the completion statement.
\end{samepage}

The matching description also connects this example with Almgren's geometry of unordered tuples: using the embedding of De Lellis and Spadaro \cite{delellis-spadaro}, we obtain a bi-Lipschitz embedding of the symmetric product into a Euclidean space (\Cref{thm:zero-cycles-bilipschitz-embedding}). This embedding is definable in the globally subanalytic structure recalled in \Cref{subsect:subanalytic}.

\begin{figure}[!t]
\centering
\begin{tikzpicture}[x=1cm,y=.9cm,font=\small,
  every node/.style={align=center},
  >={Stealth[length=2mm,width=1.35mm]}]
  \path[use as bounding box] (-6.8,-3.0) rectangle (6.8,2.65);
  \draw[black!16,line width=.45pt] (0,-2.6) -- (0,2.15);

  \begin{scope}[xshift=-3.6cm]
    \node[anchor=north west] at (-3.0,2.53) {(a)};

    \begin{scope}[shift={(-1.6,1.4)},rotate=-10]
      \filldraw[fill=cycleblue!3,draw=cyclegray!45,line width=.45pt]
        (-.85,-.45) -- (.5,-.45) -- (.85,.45) -- (-.5,.45) -- cycle;
      \draw[cycleblue,line width=1.05pt] (0,0) ellipse[x radius=.57,y radius=.27];
    \end{scope}
    \node[text=cycleblue] at (-1.6,.6) {smooth conic};

    \draw[cycleteal,line width=1.05pt] (-2.25,-1.34) -- (-.95,-.66);
    \draw[white,line width=4pt] (-2.15,-.55) -- (-1.05,-1.45);
    \draw[cycleteal,line width=1.05pt] (-2.15,-.55) -- (-1.05,-1.45);
    \node[text=cycleteal] at (-1.6,-1.9) {skew lines};

    \draw[cyclerust,line width=1.05pt] (.95,-.23) -- (2.25,.53);
    \draw[cyclerust,line width=1.05pt] (.95,.53) -- (2.25,-.23);
    \fill[cyclerust] (1.6,.15) circle[radius=1.75pt];
    \node[text=cyclerust] at (1.6,-.72) {intersecting\\lines};

    \draw[cycleblue,line width=.85pt,->]
      (-.65,1.2) .. controls (-.03,1.06) and (.2,.62) .. (.7,.34);
    \draw[cycleteal,line width=.85pt,->]
      (-.65,-.89) .. controls (-.03,-.73) and (.2,-.32) .. (.7,-.04);
    \node[text=cyclegray,font=\footnotesize] at (.35,1.48) {degenerate};
    \node[text=cyclegray,font=\footnotesize,text width=5.5cm] at (0,-2.65)
      {one cycle in the full space,\\separate limits in the two completions};
  \end{scope}

  \begin{scope}[xshift=3.6cm]
    \node[anchor=north west] at (-3.0,2.53) {(b)};

    \draw[cycleblue,line width=1.05pt] (-2.55,1.2) -- (2.55,1.2);
    \node[text=cycleblue,anchor=south] at (1.95,1.28) {fixed $S$};
    \draw[cyclegray!70,line width=.9pt,dash pattern=on 3pt off 2.2pt]
      (-.9,.06) -- (-.9,2.02);
    \draw[cycleteal,line width=1.15pt] (.9,.06) -- (.9,2.02);
    \node[text=cyclegray] at (-.9,2.24) {$F_s$};
    \node[text=cycleteal] at (.9,2.24) {$F_t$};
    \draw[cyclegray,line width=.7pt,->] (-.62,1.8) -- (.62,1.8);
    \draw[cyclegray,line width=.65pt,fill=white] (-.9,1.2) circle[radius=1.7pt];
    \fill[cyclerust] (.9,1.2) circle[radius=2.0pt];

    \filldraw[fill=cycleblue!3,draw=cycleblue!55,line width=.65pt]
      (0,-1.5) circle[radius=1.0];
    \draw[cycleblue!35,line width=.45pt,dashed]
      (1,-1.5) arc[start angle=0,end angle=180,x radius=1,y radius=.24];
    \draw[cycleblue!40,line width=.5pt]
      (-1,-1.5) arc[start angle=180,end angle=360,x radius=1,y radius=.24];
    \draw[cyclegray!60,line width=.65pt,densely dotted,->]
      (-.9,.02) .. controls (-.9,-.48) and (-.65,-.72) .. (-.65,-1.68);
    \draw[cyclegray!60,line width=.65pt,densely dotted,->]
      (.9,.02) .. controls (.9,-.48) and (.65,-.72) .. (.65,-1.68);
    \node[text=cyclegray] at (-1.48,-.22) {$\pi$};
    \draw[cycleteal,line width=.95pt,->]
      (-.65,-1.68) .. controls (-.3,-1.77) and (.3,-1.77) .. (.65,-1.68);
    \fill[cyclegray] (-.65,-1.68) circle[radius=1.65pt];
    \fill[cycleteal] (.65,-1.68) circle[radius=1.65pt];
    \node[text=cyclegray] at (-.76,-1.98) {$s$};
    \node[text=cycleteal] at (.76,-1.98) {$t$};
    \node[text=cycleblue] at (0,-2.78) {parameter sphere $\CP^1$};
  \end{scope}
\end{tikzpicture}
\caption[Two roles of singular cycles]{%
Two roles of singular cycles, shown as incidence schematics.
(a) Both families reach the same pair of intersecting lines in the full
cycle space, while their Weil--Petersson completions retain separate limits
(\Cref{ex:conics-p3-completion}).
(b) $S$ is the fixed section, $\pi$ the projection to the base, and
$F_s,F_t$ the fibers over $s,t\in\CP^1$. The dashed and solid fibers show
two positions of the moving fiber in $[S]+[F_t]$. Every cycle is singular,
while the parameter sphere has a smooth metric, round for a rotationally
invariant ambient form (\Cref{ex:hirzebruch-inner-geometry}).}
\label{fig:singular-cycles}
\end{figure}
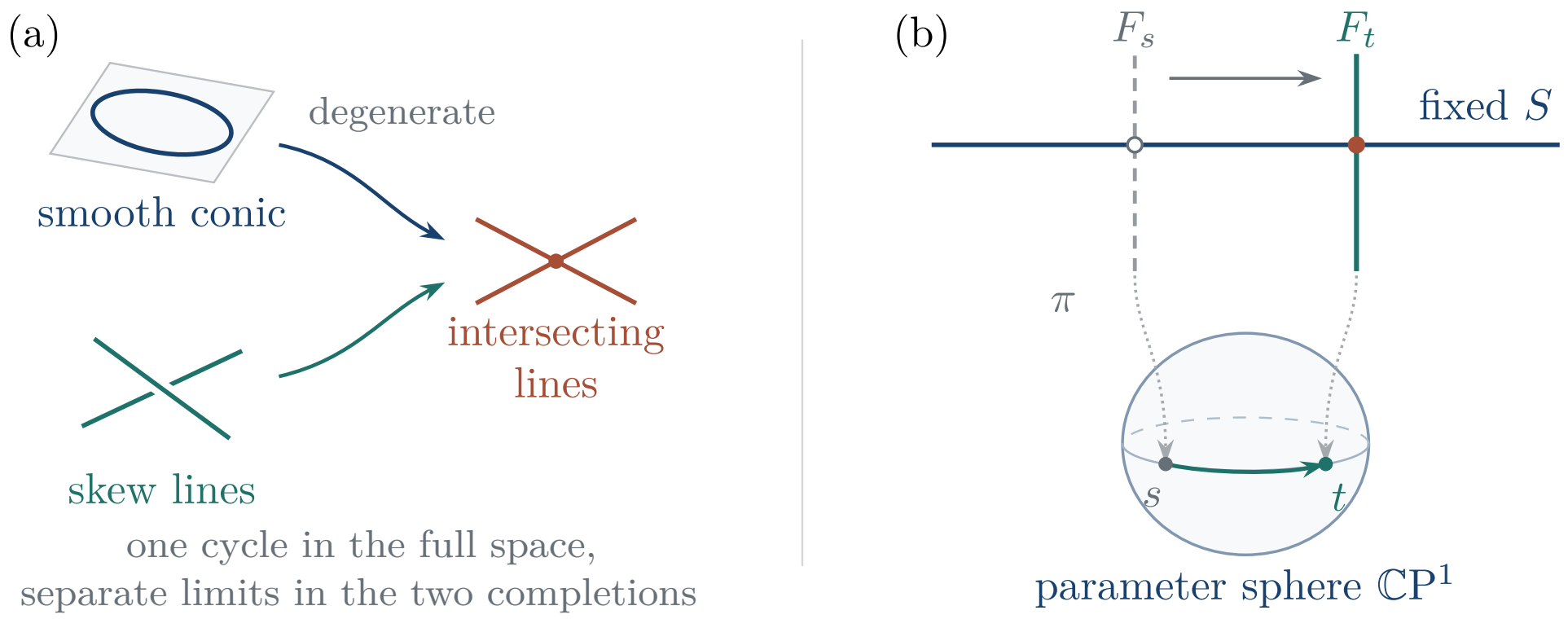

The zero-dimensional case suggests comparing the full cycle space with the metric completion of its good locus. For an irreducible component $C$ with nonempty good locus, we prove that the Weil--Petersson completion admits a canonical surjection onto $C$ equipped with its own inner Wasserstein distance, and that this map does not increase distances (\Cref{thm:wp-completion-comparison}). 

\begin{introexample}
The space of degree-two curves in $\CP^3$ contains both smooth planar conics and pairs of skew lines, that is, disjoint projective lines. These lie in two different irreducible components of the same connected cycle space. A conic can degenerate, within its plane, to a pair of intersecting lines. The same singular cycle can also be reached by moving two skew lines until they meet. Following one degeneration and reversing the other gives a finite-length path from a smooth conic to a pair of skew lines. Along this path the support changes from a connected smooth curve to a disconnected smooth curve, passing through their common singular limit. The two good loci carry their own Weil--Petersson length distances: no path in their union joins the two families. Taking the metric completion of each gives one limiting point from the conics and one from the skew lines. In the disjoint union of these two completions the points remain distinct, although they represent the same cycle in the full cycle space. The inner Wasserstein geometry on the whole connected component joins the families through that common cycle; see \Cref{ex:conics-p3-completion} for the explicit degenerations and completion maps. See also \Cref{fig:singular-cycles}(a). 
\end{introexample}

The construction also applies when there are no good cycles at all. In that situation there is no good-locus Weil--Petersson metric from which to start, but motion of the cycle measures still determines an intrinsic geometry.

\begin{introexample}
Start with the holomorphic tangent bundle of the projective line $\CP^1$. Each fiber is a complex line; adding one point at infinity to every fiber produces a compact complex surface fibered by projective lines. The added points form a smooth curve meeting every fiber once, called the section at infinity. Keep this section fixed and take its union with one moving fiber. Every union has a transverse intersection, so all the resulting cycles are singular. They form an entire irreducible connected component of the cycle space, parametrized by the base $\CP^1$. Only the fiber moves, while the section contributes a fixed part of the volume measure. Nevertheless, the induced inner Wasserstein distance is the distance of a smooth Riemannian metric on the parameter sphere; for a rotationally invariant ambient K\"ahler form, this sphere is round up to scale. Thus an entire component of singular cycles can have a smooth metric geometry, even though it contains no good cycles whose Weil--Petersson geometry could be completed. The surface in this construction is the second Hirzebruch surface; the construction and metric computation are given in \Cref{ex:hirzebruch-inner-geometry}. The fixed section and moving fiber are illustrated in \Cref{fig:singular-cycles}(b).
\end{introexample}

\medskip
\noindent\textbf{Plan of the paper.}
\Cref{sec:preliminaries} collects the background on analytic cycles,
deformation theory, Weil--Petersson geometry, and optimal transport,
together with the analytic tools used in the proofs.
In \Cref{sec:wasserstein-cycle-geometry} we construct the inner
Wasserstein distance, prove its finiteness, completeness and geodesicity,
and relate it to Weil--Petersson geometry and the completion of the
good locus. \Cref{subsec:examples} presents examples, including a
component consisting entirely of singular cycles.
Finally, \Cref{sec:zerocycles} treats zero-dimensional cycles,
their optimal-matching description, and their realization in Euclidean
space.

\medskip
\noindent\textbf{History of this work.}
The present paper is part of a research programme developed over several
years, including Lorenzo Cecchi's doctoral research. Its starting point
is the optimal-transport approach to projective hypersurfaces introduced
by Antonini, Cavalletti, and Lerario \cite{ACL}. 
The results from this paper, containing the extension of this programme to algebraic and analytic cycle spaces,
including its relationship with Weil--Petersson geometry, were announced
by Antonio Lerario in his plenary FoCM lecture in Vienna in July 2026 \cite{lerario-focm2026}.

\medskip
\noindent\textbf{Use of artificial intelligence.}
The main ideas, definitions, and constructions of this paper originate in a research programme that predates the widespread availability of generative AI tools. We regard the main contribution of this work as the organization of these ideas into a coherent geometric framework, connecting optimal transport with classical deformation theory. Within this framework, the Weil--Petersson geometry of smooth families and the transport of singular cycles emerge as aspects of the same construction.

ChatGPT 5.6 helped us complete our proof of \cref{thm:finiteness} by suggesting a result from geometric measure theory (\cref{lem:gmt-removable-frontier}), and derived the lower bound in \cref{thm:local-bilipschitz}. Afterwards, ChatGPT 5.6 was also used to simplify some arguments and to revise the language and exposition of the manuscript. The authors evaluated and revised its contributions, made the final mathematical and editorial decisions, and take full responsibility for the correctness of the results, the proofs, and the final text.

\section{Preliminaries}\label{sec:preliminaries}
\subsection{The cycle space}\label{sec:cycle-space}

Let \((M,\omega)\) be a compact connected K\"ahler manifold of complex dimension \(n\), and fix an integer
\(q\in\{0,\dots,n\}\). We denote by \(\mathcal{C}_q(M)\) the Barlet space of
\(q\)-dimensional analytic cycles in \(M\) \cite[Chapter~5]{cycles2}. A point of \(\mathcal{C}_q(M)\) is a finite formal sum
\[
\gamma=\sum_{j=1}^N m_j[Z_j],
\]
where the \(Z_j\subset M\) are distinct irreducible complex analytic subsets of pure
dimension \(q\) and \(m_j\in\mathbb{N}_{\ge 1}\).
We write
\[
|\gamma| \coloneqq \bigcup_{j=1}^N Z_j
\]
for the support of \(\gamma\). We will say that $\gamma$ is \emph{smooth} if its support is a smooth complex submanifold, and that it is \emph{reduced} if,  in the above decomposition, we have  $m_j=1$ for all $j=1, \ldots, N$. If
\[
\gamma=\sum_{j=1}^N m_j[Z_j]
\]
is such a cycle and \(x\in M\), its multiplicity at \(x\) is defined by
\begin{equation}\label{eq:multiplicity}
\operatorname{mult}_x(\gamma)
 \coloneqq 
\sum_{j:\,x\in Z_j} m_j\,\operatorname{mult}_x(Z_j),
\end{equation}
where \(\operatorname{mult}_x(Z_j)\) denotes the multiplicity of the analytic set \(Z_j\) at
\(x\) (see \cite[Definition~3.1.32]{cycles1}). In particular, \(\operatorname{mult}_x(\gamma)=0\) for \(x\notin |\gamma|\), and if
\(x\) is a smooth point of a unique irreducible component \(Z_j\), then
\(\operatorname{mult}_x(\gamma)=m_j\). Thus the function
\(x\mapsto \operatorname{mult}_x(\gamma)\) records both the coefficients of the cycle and the local
singular multiplicities of its irreducible components. 

For a pure-dimensional analytic set $Z$ with irreducible components $Z_i$, we write
\begin{equation}\label{eq:reduced-fundamental-cycle}
 [Z]_{\mathrm{red}}:=\sum_i[Z_i]
\end{equation}
for its \emph{reduced fundamental cycle}: every component has multiplicity one.
The sum is locally finite, and finite when $Z$ is compact. If $Z$ is irreducible, we simply denote $[Z]_{\mathrm{red}}=[Z]$.

Throughout the paper, by a \emph{complex space} we mean a reduced complex
space in the sense of \cite[Definition 2.4.3]{cycles1}. Thus, a complex space
is a Hausdorff topological space, countable at infinity, which is locally
isomorphic, as a complex space, to a complex analytic subset of an open set in
some \(\C^N\). Complex analytic subsets are always endowed with their reduced
induced structure. Reduced complex spaces generalize complex manifolds by allowing singularities.
A point \(y\) of a reduced complex space \(Y\) is called \emph{regular}, or
\emph{smooth}, if it admits a neighborhood which is isomorphic, as a complex
space, to a complex submanifold of an open set in some \(\C^N\). The set of
smooth points is denoted by \(Y_{\mathrm{reg}}\), and its complement is the
singular locus \(Y_{\mathrm{sing}}\).

By an \emph{analytic subset} of a complex space $X$ we mean a
closed subset locally defined by the vanishing of finitely many
holomorphic functions on $X$. We always endow it with its reduced
induced complex-space structure. 

In our setting, the cycle space \(\mathcal{C}_q(M)\) carries a
structure of reduced complex space \cite[Corollary 5.8.2.6]{cycles2}.
If \(\alpha\) is a continuous \((q,q)\)-form on \(M\), its integral over \(\gamma\) is defined by
\begin{equation}\label{eq:int_of_form}
\int_\gamma \alpha
:=
\sum_{j=1}^N m_j\int_{(Z_j)_{\mathrm{reg}}}\alpha,
\end{equation}
where \((Z_j)_{\mathrm{reg}}\) denotes the smooth locus of \(Z_j\); the integral is finite by \cite[Theorem 3.2.26]{cycles1}.  In
other words, one integrates over the smooth part of each irreducible component and weights the
result by its coefficient in the cycle. In particular, the K\"ahler form \(\omega\) defines the
\(\omega\)-volume of a cycle \(\gamma\) by
\[
\operatorname{vol}_\omega(\gamma):=\int_\gamma \frac{\omega^q}{q!}.
\]
Note that, by Wirtinger's theorem \cite[p.~31]{GH}, this is the total \(q\)-dimensional K\"ahler volume
of \(\gamma\), counted with multiplicities.

Since \(M\) is compact, the topology of \(\mathcal{C}_q(M)\) (as introduced in \cite[\S4.2.2]{cycles1}) admits the following convenient
characterization: it is the coarsest topology such that, for every closed subset \(F\subseteq M\),
the set
\[
\{\gamma\in\mathcal{C}_q(M): |\gamma|\cap F=\varnothing\}
\]
is open, and for every continuous \((q,q)\)-form \(\alpha\) on \(M\), the function
\[
\gamma\longmapsto \int_\gamma \alpha
\]
is continuous \cite[Theorem~4.2.26]{cycles1}.

Throughout the paper, we fix a connected component
\[
B\subset \mathcal{C}_q(M).
\]
We assume that $B$ does not consist of the zero cycle, so its cycles have positive volume.
Note that $B$ inherits a reduced complex-space
structure and is an analytic subset of $\mathcal C_q(M)$.

Whenever needed, we denote by \(C\subset B\) an irreducible component of \(B\).
Since \(M\) is a compact K\"ahler manifold, the connected components of \(\mathcal{C}_q(M)\) are
compact \cite[Corollary~2.15]{scv7}. Moreover, the function
\begin{equation}
\gamma\longmapsto \operatorname{vol}_\omega(\gamma)=\int_\gamma \frac{\omega^q}{q!}\label{eq:VB}
\end{equation}
is constant on connected components of \(\mathcal{C}_q(M)\)
\cite[proof of Theorem~12.2.3.5]{cycles2}, and we denote by
\begin{equation}\label{eq:V}
V_B:=\operatorname{vol}_\omega(\gamma), \qquad \gamma\in B,
\end{equation}
this common positive value.

\begin{remark}\label{remark:projective}
If \(M\) is projective, the Barlet cycle space agrees with the classical Chow
variety. More precisely, let \(M\subset \CP^N\) be a complex projective
variety. By Chow's theorem, every compact analytic subvariety of \(M\) is algebraic.
Thus the points of \(\cC_q(M)\) are precisely the effective algebraic \(q\)-cycles
on \(M\). Moreover, the locus
of cycles of degree \(d\) is naturally identified, as a reduced complex space, with
the analytification of the usual Chow variety
\(\operatorname{Chow}_{q,d}(M)\); see \cite[Chapter~9]{cycles2}. In this
sense, in the projective case the Barlet space is the analytic incarnation of the
Chow variety:
\[
   \cC_q(M)
   \simeq
   \bigsqcup_{d\geq 0}
   \operatorname{Chow}_{q,d}(M) .
\]
For example, when \(M=\CP^n\) and \(q=n-1\), the degree-\(d\) part of the
cycle space is the projective space
\[
   \cC_{n-1,d}(\CP^n)
   \simeq
   \PP(H^0(\CP^n,\mathcal O_{\CP^n}(d)))
   \simeq
   \CP^{\binom{n+d}{d}-1}.
\]
Indeed, its points are effective hypersurface cycles of degree \(d\), equivalently
divisors of homogeneous polynomials of degree \(d\), modulo multiplication by a
nonzero scalar. See \cref{ex:hypersurfaces-cpn} below for a more detailed discussion of this example in the context of this paper.
\end{remark}
\subsubsection{The tautological family}\label{sec:tautological}
The cycle space $\cC_q(M)$ carries a \emph{tautological} family of
compact $q$-cycles, defined as follows:
\begin{equation}\label{eq:tautologicalX}
\mathcal X:=\{(\gamma,x)\in \cC_q(M)\times M:\ x\in |\gamma|\}.
\end{equation}
By \cite[Theorem 4.6.1]{cycles1}, the set $\mathcal{X}$ is a closed analytic subset of $\cC_q(M)\times M$, and the projections on the two factors induce holomorphic maps
\begin{equation}\label{eq:tautological-diagram}
\begin{tikzcd}[row sep=2.2em, column sep=2.8em]
\mathcal X \arrow[r,"j"] \arrow[d,"\pi"'] & M \\
\cC_q(M)&
\end{tikzcd}
\end{equation}
with $\pi$ proper (by compactness of $M$).

We recall what analyticity means for a family $(\gamma_s)_{s\in S}$ of $q$-cycles,
where $S$ is a reduced complex space. Locally, choose product coordinates
$U\times V\subset\C^q\times\C^{n-q}$ in $M$ in which the projection of a cycle
onto $U$ is a finite branched covering. The family is \emph{analytic} if, in every
such chart adapted to $\gamma_{s_0}$, the same projection remains finite of a fixed
degree $d$ for $s$ near $s_0$, and its $d$ points over $u\in U$, counted with their
cycle multiplicities, define a holomorphic map
\[
 S_0\times U\longrightarrow\operatorname{Sym}^d(V).
\]
Here $S_0$ is a neighborhood of $s_0$, and $\operatorname{Sym}^d(V)$ denotes the
space of unordered $d$-tuples in $V$; see \cite[Chapter~4]{cycles1}.
The family is \emph{proper} if its incidence set
$\mathcal X_S:=\{(s,x):x\in|\gamma_s|\}$ projects properly to $S$, that is, the
inverse image of every compact subset of $S$ is compact. For an analytic family
in our compact ambient manifold $M$, this follows from the closedness of
$\mathcal X_S\subset S\times M$.

Moreover, the analytic structure of the Barlet space is built in such a way that  the tautological family is universal in the following sense. If \(S\) is a
reduced complex space and \((\gamma_s)_{s\in S}\) is a proper analytic family  of
compact \(q\)-cycles in \(M\), then the associated classifying map
\[
   f:S\longrightarrow \cC_q(M),\qquad s\longmapsto \gamma_s,
\]
is holomorphic. For this reason the tautological family is sometimes also called \emph{universal}.

\subsubsection{Graph-cycles and multiplicities}
\label{subsubsec:graph-cycle-resolution}

An analytic family has both an incidence set, which records its supports, and a
graph-cycle, which also records multiplicities. For a holomorphic arc $\alpha:\Delta\to\mathcal C_q(M)$, where $\Delta\subset\C$ is an open disk, its graph-cycle is a locally finite analytic cycle
\[
 \mathscr G_\alpha=\sum_a m_a[G_a]
 \quad\text{on }\Delta\times M,
 \qquad
 |\mathscr G_\alpha|=\{(t,x):x\in|\alpha(t)|\},
\]
of pure dimension $q+1$; see \cite[Theorem~4.3.11]{cycles1}.
The word ``locally'' allows its support to be noncompact in the parameter direction.
Its intersection with the parameter fiber is the original cycle, with its multiplicities:
\[
 \mathscr G_\alpha\cap[\{t\}\times M]=\{t\}\times\alpha(t)
 \qquad\text{\cite[Proposition~7.1.5.1]{cycles2}}.
\]
The geometric properties of this graph needed for the length estimate are established in
\Cref{ss:arc-graph}.

We also recall the action of a holomorphic map on cycles. Let $h:X\to X'$ be
holomorphic and proper on the support of a $p$-cycle $\gamma=\sum_i m_i[Z_i]$.
Its \emph{push-forward} is the $p$-cycle
\[
 h_*\gamma:=\sum_{\dim h(Z_i)=p}m_i\,d_i[h(Z_i)],
 \qquad d_i:=\deg(h|_{Z_i}:Z_i\to h(Z_i)),
\]
where $d_i$ is the number of points in a generic fiber, counted with multiplicity.
Components whose images have smaller dimension contribute zero. In particular,
an isomorphism preserves the multiplicities of the components.
This is how we will identify cycles in $\{t\}\times M$ with cycles in $M$.
We say that a map \emph{dominates} its target if its image is dense in that target;
for a proper holomorphic map, whose image is closed, this means that it is surjective.

\subsection{Regular points of the cycle space}
A point of the cycle space carries two types of regularity properties: one concerns the
parameter space itself, the other the cycle represented by that point.

\begin{definition}\label{def:two-regularities}
Let \(B\subset \mathcal{C}_q(M)\) be a connected component, and let \(\gamma\in B\).

\begin{enumerate}
\item
We say that \(\gamma\) is \emph{regular} if \(B\) is smooth at \(\gamma\). We denote the
corresponding locus by
\[
B_{\mathrm{reg}}:=\{\gamma\in B : B \text{ is smooth at } \gamma\}.
\]

\item
We say that \(\gamma\) is \emph{smooth and reduced} if the cycle represented by \(\gamma\) is a
smooth reduced compact complex submanifold of \(M\). We denote the corresponding locus by
\[
B^{\mathrm{sm,red}}
:=
\{\gamma\in B : \gamma \text{ represents a smooth reduced } q\text{-cycle in }M\}.
\]
\end{enumerate}

In the sequel we will often need to consider points of \(B\) satisfying both conditions at
once. We therefore introduce also the following notation.

\begin{enumerate}
\setcounter{enumi}{2}
\item
We denote by
\[
\widetilde{B}:=B_{\mathrm{reg}}\cap B^{\mathrm{sm,red}}.
\]
We call \(\widetilde{B}\) the \emph{good locus} of \(B\) and its elements the \emph{good points}.
\end{enumerate}
\end{definition}
Similarly, if \(C\subset B\) is an irreducible component, we define
\[
C_{\mathrm{reg}}:=\{\gamma\in C : C \text{ is smooth at } \gamma\},\qquad
C^{\mathrm{sm,red}}:=C\cap B^{\mathrm{sm,red}},\qquad
\widetilde C:=C_{\mathrm{reg}}\cap C^{\mathrm{sm,red}}.
\]
Thus \(\widetilde C\) is the locus of points of \(C\) which are regular as points of \(C\) and
which represent smooth reduced cycles. Notice that \(\widetilde C\) need not coincide with
\(\widetilde B\cap C\). Indeed, \(\widetilde B\cap C\) records the points of \(C\) which are
regular in \(B\) and represent smooth reduced cycles, whereas \(\widetilde C\) records the
points which are regular in \(C\) and represent smooth reduced cycles. In particular one always has $
\widetilde B\cap C \subset \widetilde C,$
and equality holds exactly when every point of \(\widetilde C\) is regular as a point of \(B\).

\begin{remark}
The two notions introduced in \Cref{def:two-regularities} are independent in general. For example,
let \(B\) be the component of \(\cC_{n-1}(\CP^n)\) parametrizing
degree-\(d\) hypersurfaces, and let \(\Sigma_{n,d}\subset B\) be the discriminant, namely the
subset of points representing singular hypersurfaces. Since \(B\simeq \CP^N\) is smooth, one has
$B_{\mathrm{reg}}=B.$
On the other hand, a point of \(\Sigma_{n,d}\) represents a singular cycle, hence not a smooth
reduced one. Therefore, every point of \(\Sigma_{n,d}\) belongs to \(B_{\mathrm{reg}}\) but not to
\(B^{\mathrm{sm,red}}\), showing that the two conditions are genuinely different. 
Conversely, the other implication also fails. Indeed, by Vakil's Murphy's Law for
Chow varieties, the Chow variety of nonsingular curves in projective space may be
singular \cite[Theorem~1.1, M4]{vakil}. Since in the projective case the
Barlet space agrees with the Chow variety, this gives points $
   \gamma\in B^{\mathrm{sm,red}}\setminus B_{\mathrm{reg}}.$
Thus the two conditions are independent.
\end{remark}
An irreducible component of the Barlet space can consist entirely of nonsmooth or nonreduced cycles (see \cref{rem:scope-bad-locus-theorem}). When this is not the case, the complement of the good locus is a proper analytic subset, in particular of real codimension at least two.
\begin{lemma}\label{propo:analyticbad}
In the above setting, let
\(S\subset \mathcal C_q(M)\) be an analytic subset and for $s\in S$ let \(\gamma_s\) denote the cycle
parametrized by \(s\). Then
\[
S^{\mathrm{sing}}
:=
\{s\in S:\ \gamma_s \text{ is not smooth and reduced}\}
\]
is a closed analytic subset of \(S\).
\end{lemma}

\begin{proof}
Let
\[
\mathcal{X}_S=\{(s,x)\in S\times M:\ x\in |\gamma_s|\}
\]
be the restriction of the tautological family to $S$, and let \(\pi:\mathcal{X}_S\to S\) be the projection.
Since \(M\) is compact, \(\pi\) is proper. For \(m\ge 1\), define the set
\[
F_m
=
\{(s,x)\in \mathcal{X}_S:\operatorname{mult}_x(\gamma_s)\ge m\},\]
where the pointwise multiplicity is defined in \eqref{eq:multiplicity}. By \cite[Proposition 4.3.10]{cycles1}, \(F_m\)
is a closed analytic subset of \(\mathcal{X}_S\). In particular \(F_2\) is analytic. Since \(\operatorname{mult}_x(\gamma)=1\) if and only if, near \(x\), the cycle $\gamma$ is smooth and reduced, we have that $S^{\mathrm{sing}}=\pi(F_2).$ By Remmert's proper mapping theorem, \(S^\mathrm{sing}\) is therefore analytic in \(S\). 
\end{proof}
\begin{corollary}\label{prop:bad-locus-analytic}
In the above setting, let \(C\subset B\) be an irreducible component, and assume that \(C\)
contains a smooth reduced cycle. Then $C\setminus\widetilde C$ is a proper closed analytic subset of $C$, and
\(\widetilde C\) is a dense, connected, complex
submanifold of \(C\).\end{corollary}

\begin{proof}Observe first that $\widetilde{C}$ is the complement of the union of the two analytic sets consisting of singular points of $C$ and, respectively, points of $C$ parametrizing cycles that are not smooth and reduced. Irreducible components of the cycle space are themselves analytic. Therefore the set of singular points of $C$ is proper analytic by \cite[p.~117]{Grauert1984}; the set of cycles in $C$ that are not smooth and reduced is analytic by \cref{propo:analyticbad} and proper by assumption.  Since proper analytic sets have real codimension at least two, the conclusion follows.
\end{proof}

\begin{example}\label{rem:scope-bad-locus-theorem}

The bad locus may contain whole irreducible components of the cycle space. Indeed, let \(M\) be
the second Hirzebruch surface, that is, the projectivization of the vector bundle
\(E=\mathcal{O}\oplus \mathcal{O}(-2)\to \CP^1\)
(see \cite[Chapter V, Section 2]{hartshorne}). We denote by
\(\pi:M\to \CP^1\) the induced projection. Its fibers are projective lines, and we denote by
\(S\subset M\) the negative section. Thus \(S\) is an analytic curve, \([S]\) is its reduced
fundamental cycle, and \(S^2=-2\)
(see \cite[Chapter V, Section 2, Proposition 2.9]{hartshorne}). Let \(F_t\subset M\) be the
fiber over \(t\in\CP^1\). Consider the family of cycles
\[
\Gamma:=\{[S]+[F_t]:t\in\CP^1\}\subset \mathcal C_1(M).
\]
It is an analytic family of compact cycles parametrized by \(\CP^1\), and hence defines a compact
irreducible analytic subset of the Barlet space. We claim that \(\Gamma\) is exactly the locus of
effective cycles with homology class
\(\operatorname{cl}(S)+\operatorname{cl}(F)\), where \(F\) denotes any fiber. Recall that
\[ 
F^2=0,\qquad S\cdot F=1,
\]
see \cite[Chapter V, Section 2, Proposition 2.3]{hartshorne}. Since
\((\operatorname{cl}(S)+\operatorname{cl}(F))\cdot\operatorname{cl}(S)=-1<0\), a standard
intersection-number argument implies that any effective cycle in this class must contain
\([S]\) as an irreducible component. The residual cycle has homology class
\(\operatorname{cl}(F)\), and therefore is the reduced cycle of a fiber. Thus
\[
\Gamma=
\left\{
\gamma\in\mathcal C_1(M):
\operatorname{cl}(\gamma)=\operatorname{cl}(S)+\operatorname{cl}(F)
\right\}.
\]
Since the homology class of a cycle is locally constant in analytic families
(see \cite[Remark (ii), p. 420]{cycles1}), any irreducible component of
\(\mathcal C_1(M)\) meeting \(\Gamma\) must be contained in this fixed homology class, and hence
in \(\Gamma\). Therefore \(\Gamma\) is both a connected component and an irreducible component of
\(\mathcal C_1(M)\). Finally, every cycle \([S]+[F_t]\) is reduced, but it is singular because
\(S\) and \(F_t\) meet transversely at one point. Thus \(\Gamma\) is an irreducible curve in the
Barlet space all of whose members are reduced and singular.

\end{example}

\subsection{The tangent space at good points}\label{Ss:tangentgood}

We describe tangent vectors to an irreducible component $C\subset B$ at a good
point $\gamma=[Z]\in\widetilde C$, where $Z\subset M$ is a smooth compact complex
submanifold. The basic fact is the existence of a canonical injective complex-linear map
\begin{equation}\label{eq:tangent-normal}
 \kappa_\gamma:T_\gamma C\hookrightarrow H^0(Z,N_{Z/M}),
 \qquad v\longmapsto\nu_v,
\end{equation}
where $N_{Z/M}:=(TM|_Z)/TZ$ is the holomorphic normal bundle. Thus an infinitesimal
motion of the cycle determines a holomorphic normal field along $Z$, and this field
determines the tangent vector. Surjectivity onto all holomorphic normal fields is
not required. We will shortly describe $\nu_v$ explicitly using isotopies.

\begin{remark}\label{rem:douady}For completeness, we explain the deformation-theoretic origin of this map.
The auxiliary \emph{Douady space} $D_q(M)$ parametrizes compact complex subspaces
of $M$ near $Z$ and is allowed to be nonreduced. Its role here is to justify
\eqref{eq:tangent-normal} and the smooth local family described below; these two
facts can be taken as standard inputs for the rest of the paper.
Near $[Z]$, the Barlet space is canonically isomorphic to $D_q(M)_{\mathrm{red}}$
\cite[Corollary~10.3.4.4]{cycles2}. The inclusions of $C$ into $B$ and of the
reduction into the Douady space give
\[
 \begin{aligned}
 T_\gamma C&\hookrightarrow T_\gamma B
       \simeq T_{[Z]}D_q(M)_{\mathrm{red}},\\
 T_{[Z]}D_q(M)_{\mathrm{red}}&\hookrightarrow T_{[Z]}D_q(M)
       \simeq H^0(Z,N_{Z/M}),
 \end{aligned}
\]
where the last isomorphism is the standard description of the Douady tangent space
\cite[Proposition~1.7, p.~323]{scv7}. Their composition is $\kappa_\gamma$.
The first inclusion is an isomorphism when $C$ and $B$ coincide near $\gamma$,
in particular when $\gamma$ is regular in $B$. Our definition of $\widetilde C$
requires regularity only in $C$, so that inclusion, as well as the inclusion
induced by reduction, must in general be kept as an injection. On a sufficiently small neighborhood $U\subset\widetilde C$, the pullback of the
universal Douady family coincides with the tautological family. Smoothness of its
fibers and the local product description in \cite[Corollary~10.3.4.3]{cycles2}
show that $\mathcal X_U$ is a complex manifold and $\pi:\mathcal X_U\to U$ is a
proper holomorphic submersion.

A standard sufficient condition for $\kappa_\gamma$ to be an isomorphism is
the vanishing of $H^1(Z,N_{Z/M})$, the first sheaf cohomology group of the
holomorphic normal bundle. This group contains the obstructions to extending
infinitesimal embedded deformations; Kodaira's completeness theorem gives
smoothness of the local deformation space when it vanishes
\cite{kodaira-completeness}. More generally, smoothness of the full Douady
space at $[Z]$ suffices: its reduction then coincides locally with it, and
there is a unique local irreducible component. Smoothness of the reduction
alone does not imply this conclusion.
\end{remark}

We give now an explicit description of the injection \eqref{eq:tangent-normal} using the
language of differential topology. Forgetting the complex structure, we regard \(T_{[Z]}C\) as
a real vector space. Given \(v\in T_{[Z]}C\), let \(t\mapsto\gamma_t\) be a smooth curve in the
good locus with \(\gamma_0=[Z]\) and \(\dot\gamma_0=v\), and write
\(Z_t:=|\gamma_t|\). After possibly shrinking the parameter interval
\(I=(-\epsilon,\epsilon)\), the union
\(\{(t,x)\in I\times M:x\in Z_t\}\subset I\times M\) is a smooth submanifold and its projection
to \(I\) is a proper submersion. By Ehresmann's fibration theorem, we may choose an isotopy,
that is, a smooth family of embeddings \(F_t:Z\to M\), such that
\[
F_{0}=\mathrm{id}_{Z},
\qquad
F_t(Z)=Z_t.
\]
We define the velocity field \(\xi_t\in \Gamma\bigl(Z_t,TM|_{Z_t}\bigr)\) along \(Z_t\) by
\[
\xi_t\bigl(F_t(z)\bigr):=\frac{d}{dt}F_t(z)\in T_{F_t(z)}M,
\qquad z\in Z.
\]
Note that, if \(F_t\) is replaced by
\(F_t\circ\phi_t\), where \(\phi_t\) is a smooth family of diffeomorphisms of \(Z\) with
\(\phi_0=\mathrm{id}_{Z}\), then
\[
\frac{d}{dt}(F_t\circ\phi_t)(z)=\xi_t(F_t(\phi_t(z)))+dF_t|_{\phi_t(z)}\bigl(\dot\phi_t(z)\bigr),
\]
and the second term is tangent to \(Z_t\). Therefore, the class of $\xi_0$ modulo $TZ$ is independent of the chosen isotopy and gives a well defined  section of the normal bundle:
\[
z\longmapsto[\xi_0(z)]\in N_{Z/M}.
\]
\begin{lemma}\label{lemma:normalnu}
The normal section \(\nu_v\) is represented by the normal component of the infinitesimal
generator of any isotopy \(F_t:Z\to M\) satisfying \(F_t(Z)=Z_t\). Equivalently,
\(\nu_v=[\xi_0]\) in \(N_{Z/M}\).
\end{lemma}
\begin{proof}
Throughout the proof we use the standard identification of real tangent vectors
with their \((1,0)\)-parts; in particular, the class of a real vector field along
\(Z\) in the holomorphic normal bundle \(N_{Z/M}\) means the class of its
\((1,0)\)-part.

Following the construction from \cite[Section 2]{axelsson-schumacher}, we first recall explicitly how the map \(\kappa_\gamma\) is defined. Let
\(U\subset \widetilde{C}\) be a sufficiently small neighborhood of \([Z]\), and
let
\[
\mathcal X_U\subset U\times M
\]
be the restriction of the tautological family. We denote by
\[
\pi:\mathcal X_U\to U,
\qquad
j:\mathcal X_U\to M
\]
the two projections. Since \(U\subset\widetilde{C}\), the map \(\pi\) is a smooth
proper holomorphic submersion near the fiber \(Z=\pi^{-1}([Z])\), and
\(j|_Z:Z\hookrightarrow M\) is the inclusion.

Along \(Z\), the tangent sequence of the family and the normal sequence of the
embedding fit into the commutative diagram 
\[
\begin{tikzcd}
0 \arrow[r] & TZ \arrow[r] \arrow[d, equal] & T\mathcal X_U|_Z \arrow[r, "d\pi"] \arrow[d, "dj"] & T_{[Z]}C\otimes_{\mathbb C}\mathcal O_Z \arrow[r] \arrow[d] & 0 \\
0 \arrow[r] & TZ \arrow[r] & TM|_Z \arrow[r] & N_{Z/M} \arrow[r] & 0
\end{tikzcd}.
\]
The right vertical arrow is the sheaf morphism whose value on a constant section
\(v\in T_{[Z]}C\) is \(\nu_v\). Equivalently, \(\nu_v\) is obtained by lifting
\(v\) to a section of \(T\mathcal X_U|_Z\), applying \(dj\), and then projecting
modulo \(TZ\). This is precisely the embedded Kodaira--Spencer description of
the tangent--normal map.

Let now \(t\mapsto\gamma_t\) be a smooth curve such that
\(\gamma_0=[Z]\), \(\dot\gamma_0=v\), and \(F_t:Z\to M\) is as in the statement. Consider the map
\[
\widehat F:I\times Z\longrightarrow \mathcal X_U,
\qquad
\widehat F(t,z):=(\gamma_t,F_t(z)).
\]
This map is well-defined because \(F_t(z)\in Z_t=|\gamma_t|\). Define
\[
\eta(z):=\frac{d}{dt}\widehat F(t,z)\Big|_{t=0}.
\]
Then \(\eta\) is a tangent field to \(\mathcal X_U\) along \(Z\), and, by
construction,
\[
d\pi(\eta)=v.
\]
Thus \(\eta\) is a lift of \(v\) through the tangent sequence of the family.
Moreover, since \(j\circ\widehat F(t,z)=F_t(z)\), we have
\[
dj(\eta)(z)=\frac{d}{dt}F_t(z)\Big|_{t=0}=\xi_0(z).
\]
Therefore, by the definition of \(\kappa_\gamma\),
\[
\nu_v=[dj(\eta)]=[\xi_0]\in N_{Z/M}.\qedhere
\]
\end{proof} 
\subsection{The Weil--Petersson K\"ahler structure of the cycle space}\label{sec:WPgood}
We now recall the construction of a natural K\"ahler structure defined on the good points of each component of the cycle space of a K\"ahler manifold. Actually, this K\"ahler structure exists, in the more general sense of complex spaces, on the whole cycle space (see \cite[Chapter~12]{cycles2} and \cref{remark:current} below), but we will not need it in this generality. 

On the good locus, this structure admits a simple description in terms of the normal deformation fields introduced above. 

\begin{remark}From now on, we will use the Hermitian structure given by the K\"ahler form $\omega$ to identify the smooth bundle underlying $N_{Z/M}$ with the orthogonal complement $TZ^\perp$; in particular, for every $v\in TC$, we can identify the normal field $\nu_v$ with a section of  $TZ^\perp$.
\end{remark}
\begin{definition}\label{def:WP}
Following \cite[Definition~2.1]{axelsson-schumacher}, let
\(C\subset B\) be an irreducible component, let
\(\gamma=[Z]\in\widetilde{C}\), and let \(v,w\in T_\gamma C\). The
\emph{Weil--Petersson K\"ahler form} is defined by
\begin{equation}\label{eq:kahler-form-normal-fields}
\owp_\gamma(v,w):=\frac{1}{V_B}\int_Z \omega(\nu_v,\nu_w)\,\frac{\omega^q}{q!},
\end{equation}
where $V_B$ is the normalization constant from \eqref{eq:V}.
Similarly, we will call \emph{Weil--Petersson} also the associated Hermitian and Riemannian metrics on $\widetilde{C}$.

\end{definition}

At the same point $\gamma=[Z]$, the associated Riemannian metric is
\[
 \Theta_{\gamma}(v,w)=\frac1{V_B}\int_Z g(\nu_v,\nu_w)\,\frac{\omega^q}{q!}.
\]

It is positive definite because $\kappa_\gamma$ is injective, and it is Hermitian because
$\nu_{Jv}=J\nu_v$.

The name Weil--Petersson is motivated by the following fact. If $\pi:\mathcal{X}_C\to C$ and $j:\mathcal{X}_C\to M$ denote the two projections, then 
\begin{equation}\label{eq:barlet-fiber-integral-form}
\owp=\frac{1}{V_B}\pi_*\!\left(j^*\frac{\omega^{q+1}}{(q+1)!}\right)\, \textrm{on $\widetilde{C}$},
\end{equation}
where $\pi_{*}$ denotes integration along the fibers\footnote{Since
\(\pi:\mathcal X_U\to U\) is a smooth proper holomorphic submersion with fibers of complex dimension
\(q\), one can integrate differential forms on \(\mathcal X_U\) along the fibers of \(\pi\). In
particular, the \((q+1,q+1)\)-form \(j^*\omega^{q+1}\) gives rise to a \((1,1)\)-form
\(\pi_*(j^*\omega^{q+1})\) on \(U\). Concretely, for \(\gamma\in U\) and \(v,w\in T_\gamma U\), one
chooses lifts \(\widetilde v,\widetilde w\) along the fiber \(\pi^{-1}(\gamma)\) and sets
\[
\bigl(\pi_*(j^*\omega^{q+1})\bigr)_\gamma(v,w)
=
\int_{\pi^{-1}(\gamma)}
\iota_{\widetilde w}\iota_{\widetilde v}\,j^*\omega^{q+1},
\]
which is independent of the choice of lifts.}, see \cite[Lemma 2.3]{axelsson-schumacher} (notice, in particular, that this shows that $\owp$ is indeed a K\"ahler form on $\widetilde{C}$).

\begin{remark}\label{remark:current}
Formula \eqref{eq:barlet-fiber-integral-form} is first a statement on the good locus, where
the Weil--Petersson form is a K\"ahler form. Beyond this locus, one should distinguish
two facts. First, the cycle space is K\"ahler in the sense of complex spaces. This is due to Fujiki, Varouchas, and
Barlet--Varouchas, see \cite[Chapter~12]{cycles2} and \cite{axelsson-schumacher}. Indeed locally,
if \(i\partial\bar\partial\varphi=\omega^{q+1}\), then
\[
  \gamma\longmapsto \int_\gamma \varphi
\]
gives continuous strictly plurisubharmonic potentials on the cycle space; see
Theorem~11.5.3.2 and the proof of Theorem~12.3.4.1 in \cite{cycles2}. Strict
plurisubharmonicity is understood in the sense of reduced complex spaces
\cite[Definition~11.1.3.1 and Theorem~11.1.3.4]{cycles2}. Second, if \(C\subset B\) is irreducible and \(\widetilde C\neq\varnothing\), the smooth
Weil--Petersson form on \(\widetilde C\) extends across \(C\) as a positive closed
\((1,1)\)-current with continuous local \(\partial\bar\partial\)-potentials, strictly positive
at reduced cycles \cite[Theorem~2.4]{axelsson-schumacher}. 

\end{remark}

\subsection{Subanalytic geometry}\label{subsect:subanalytic}

We regard every reduced complex space as a reduced real-analytic space in the standard way.
Concretely, if a local model \(A\subset\C^N\) is defined by holomorphic equations 
$f_1=\cdots=f_r=0,$
then, under the identification \(\C^N\simeq\R^{2N}\), its underlying set is defined by the
real-analytic equations $
\operatorname{Re}f_j=\operatorname{Im}f_j=0$ for $j=1,\ldots,r. $
This is the real-analytic structure used below when discussing definability and real-analytic
curves.

We shall use a few standard facts from subanalytic and o-minimal geometry. We recall them
here in the form needed below. A subset of some Euclidean space is called \emph{globally
subanalytic} if it is definable in the o-minimal structure \(\R_{\mathrm{an}}\) of restricted\footnote{A restricted real analytic function is the extension
to \(\R^n\) by zero of a function which is real analytic on a neighborhood of the cube
\([-1,1]^n\). Equivalently, \(\R_{\mathrm{an}}\) is the expansion of the real ordered field by all
such restricted analytic functions.}
real analytic functions. Concretely, this means that it can be described, in Euclidean
coordinates, by a first-order formula whose terms are built from polynomial functions and
restricted real analytic functions. We refer the reader to \cite{vandendries} for more details.
In the rest of the paper we shall use the terms \emph{globally subanalytic} and
\emph{definable} interchangeably, always referring to definability in \(\R_{\mathrm{an}}\). A
map between definable sets is called definable if its graph is definable.
We shall use the following standard properties of globally subanalytic sets. The properties are
numbered for later reference.
\begin{enumerate}
\item\label{item:one} \emph{Stability under definable constructions.}
Definable sets are stable under Boolean operations, Cartesian products, projections, and
fibers. 

\item\label{item:five} \emph{Definable path connectedness.}
Any two points of a connected definable set can be joined by a continuous definable arc; see
\cite[Proposition~2.2.4]{valette}.

\item\label{item:six} \emph{One-variable preparation and Puiseux expansions.}
A definable real-valued function is real analytic on the cells of a suitable finite cell
decomposition.  Moreover, if \(h:(0,\eta)\to\R\) is definable, then for some
\(0<\varepsilon<\eta\), \(p\in\N_{\geq1}\), \(m\in\mathbb Z\), and real coefficients \(a_i\),
\[
h(t)=\sum_{i=m}^{\infty}a_i t^{i/p}\qquad(0<t<\varepsilon)
\]
with a convergent Puiseux series.  If \(h\) is bounded at the origin, no negative power occurs,
and \(t\mapsto h(t^p)\) extends real analytically across \(0\).  See
\cite[Theorem~1.8.2, Remark~1.8.3 and Proposition~1.8.4]{valette}.

\item\label{item:seven} \emph{Gradient {\L}ojasiewicz inequality.}
If \(h\) is a \(C^1\) definable function on a definable \(C^1\)-manifold $N\subset\R^P$ and extends
continuously at $x_0\in\overline N$, then there exist \(C>0\) and
\(\theta\in(0,1)\cap\mathbb Q\) such that
\[
 |h(x)-h(x_0)|^\theta\leq C|\nabla h(x)|
\]
near \(x_0\); see \cite[Corollary~2.2.9]{valette}. Here the gradient uses the induced Euclidean metric, and $x_0$ need not belong to $N$.

\item\label{item:finite-subanalytic-volume} \emph{Finite volume.} A relatively compact subanalytic set of dimension $k$ has finite $k$-dimensional Hausdorff measure; see \cite[Section~4.3]{valette}.
\end{enumerate}
In the sequel we will need the following technical results.

Let \(Y\) be a reduced complex space. A map \(c:I\to Y\), where \(I\subset\R\) is an interval,
is called \emph{real analytic} if, for every \(t_0\in I\), there are a neighborhood \(I_0\) of
\(t_0\), a complex-space chart
\[
 \psi:U\stackrel{\sim}{\longrightarrow} A\subset V\subset\C^N,
 \qquad c(I_0)\subset U,
\]
and a real-analytic map \(h:I_0\to\C^N\simeq\R^{2N}\) such that
\(h=\psi\circ c\). This definition is independent of the chosen chart.

\begin{lemma}\label{lem:complexification-real-analytic-arc}

Let $Y$ be a reduced complex space and let $c:(-\delta,\delta)\to Y$ be real analytic in the
preceding sense. After decreasing $\delta$, there are a disc
$\Delta\subset\C$ centered at the origin and a holomorphic map $\alpha:\Delta\to Y$ such that
\[
 \alpha(t)=c(t)\qquad(-\delta<t<\delta,\ t\in\R).
\]

\end{lemma}

\begin{proof}
Choose a complex-space chart identifying a neighborhood of $c(0)$ with a closed analytic subset
$A$ of an open set $U\subset\C^N$.  The coordinate functions of $c$ have convergent real power
series and hence extend holomorphically to a sufficiently small disc; call the resulting map  $\widehat\alpha:\Delta \to U$.  If $g_1,\ldots,g_r$ cut out $A$ locally, every holomorphic function
$g_j\circ\widehat\alpha$ vanishes on a real interval and therefore, by the identity theorem for holomorphic functions, on the disk.  Thus $\widehat\alpha(\Delta)\subset A$, and transporting the map back through the
chart proves the assertion.
\end{proof}

We isolate the one-dimensional consequence that will later convert a subanalytic curve into a
finite collection of genuinely real-analytic arcs.  

\begin{proposition}
\label{prop:subanalytic-analytic-reparametrization}

Let \(A\subset\R^N\) be a compact real-analytic set (and hence a definable set). If
$c:[a,b]\to A$ is continuous and definable, then there are finitely many compact subintervals
$I_1,\ldots,I_L$ with pairwise disjoint interiors and union $[a,b]$, compact intervals
$K_1,\ldots,K_L$, and monotone homeomorphisms
\[
 \rho_\ell:K_\ell\longrightarrow I_\ell
\]
such that every $c\circ\rho_\ell$ extends to a real-analytic map from a neighborhood of
$K_\ell$ into $A$.

\end{proposition}

\begin{proof}
Set $h:=c=(h_1,\ldots,h_N)$. Apply the analytic cell-decomposition part of
Property~\ref{item:six} simultaneously to
the finitely many coordinate functions.  We obtain a partition
\[
 a=s_0<s_1<\cdots<s_m=b
\]
such that every $h_i$ is real analytic on each open interval $(s_{j-1},s_j)$.

Fix one such interval and consider its left endpoint $s=s_{j-1}$.  By the one-variable
preparation and Puiseux-expansion part of Property~\ref{item:six},
after decreasing $\delta>0$, each coordinate has a convergent Puiseux expansion
\begin{equation}\label{eq:puiseux-coordinate-expansion}
 h_i(s+t)=\sum_{k=0}^{\infty}a_{i,k}t^{k/p_i}
 \qquad(0<t<\delta).
\end{equation}
There are no negative exponents because $h$ is bounded, and the constant terms are the coordinates
of $h(s)$ because $h$ is continuous.  Let $p$ be a common multiple of the finitely many
denominators $p_i$.  Substituting $t=\tau^p$ in
\eqref{eq:puiseux-coordinate-expansion} gives convergent ordinary power series in $\tau$.
Consequently
\[
 \tau\longmapsto h(s+\tau^p),
 \qquad 0\leq\tau\leq\delta^{1/p},
\]
extends real analytically through $\tau=0$.  The same argument at the right endpoint, applied to
$h(s_j-t)$, produces a power substitution there.  After choosing the two endpoint neighborhoods
disjoint, the remaining closed middle interval needs no change of parameter: every coordinate of
$h$ is already real analytic on a neighborhood of it.  Thus each original cell is the union of at
most three compact pieces with the asserted monotone parametrizations.  Repeating this for the
finitely many cells produces the intervals $I_\ell,K_\ell$ and maps $\rho_\ell$.

It remains to check that the analytic ambient extensions still take values in \(A\). Work
near one endpoint and let $\widetilde h$ denote the extension just constructed.  Locally the
real-analytic set \(A\) is the common zero set of finitely many real-analytic functions
$g_1,\ldots,g_r$.  For every $k$, the real-analytic function $g_k\circ\widetilde h$ vanishes on a
nonempty one-sided interval, because there it agrees with $h$.  The one-variable identity theorem
therefore makes it vanish on a full neighborhood of the endpoint. Hence \(\widetilde h\) takes
values in \(A\), so it is the required extension of \(c\circ\rho_\ell\).
\end{proof}

We will also use the following property of the tautological family of cycles.

\begin{proposition}\label{prop:definable-tautological-diagram}
Let \((M,\omega)\) be a compact K\"ahler manifold, and let \(B\) be a connected
component of the Barlet space \(\mathcal C_q(M)\). Let
\[
\mathcal X_B:=\{(\gamma,x)\in B\times M:\ x\in |\gamma|\}
\]
be the tautological family over \(B\). Then, after suitable real-analytic embeddings into
Euclidean spaces, \(B\), \(M\), and \(\mathcal X_B\) become compact definable sets, and the
projections
\(
\mathcal X_B\longrightarrow B
\)
and
\(
\mathcal X_B\longrightarrow M
\)
become definable maps.
\end{proposition}

\begin{proof}

The manifold \(M\) and the Barlet component \(B\) are compact reduced complex spaces. By the
realification convention recalled above, their underlying spaces are compact reduced
real-analytic spaces: locally, equations for a complex embedding are replaced by their real and
imaginary parts in \(\R^{2N}\).
Since the compact reduced real-analytic space \(B\), and each of the finitely many connected
components of \(M\), have uniformly bounded Zariski tangent dimension, the embedding theorem of
Acquistapace--Broglia--Tognoli
\cite[Theorem~1]{acquistapace-broglia-tognoli} gives real-analytic embeddings
\[
\iota_B:B\hookrightarrow \R^N,
\qquad
\iota_M:M\hookrightarrow \R^P,
\]
whose images are real-analytic subsets of Euclidean space. Since \(B\) and \(M\) are compact,
these images are compact real-analytic sets. In particular, they are globally subanalytic, hence
definable. The tautological family \(\mathcal X_B\) is a closed analytic subset of \(B\times M\).
Therefore, under the product embedding
\[
\iota_B\times \iota_M:B\times M\hookrightarrow \R^N\times\R^P,
\]
its image is a compact real-analytic subset of \(\R^N\times\R^P\). It is therefore globally
subanalytic, hence definable. Finally, with respect to these embeddings, the tautological
projections are restrictions to \(\mathcal X_B\) of the two coordinate projections of
\(\R^N\times\R^P\), and hence are definable.
\end{proof}

\begin{corollary}\label{cor:piecewise-analytic-connectivity}
Any two points of a connected analytic subset $A\subset B$ can be joined in $A$ by
a finite concatenation of real-analytic arcs, each extending real analytically through its
endpoints. In particular, this holds for $A=B$.
\end{corollary}
\begin{proof}
By \Cref{prop:definable-tautological-diagram}, the closed analytic set $A$ is a compact
real-analytic subset of a Euclidean space. Property~\ref{item:five} gives a continuous
definable connecting curve, and \Cref{prop:subanalytic-analytic-reparametrization}
gives its analytic pieces.
\end{proof}

\subsection{Integration currents and Stokes' theorem}
\label{ss:currents}

In the proof of finiteness, we will integrate over the regular part of the graph
of a family above a real parameter interval. This is a smooth manifold, but it
need not be closed in the ambient space: singular points have been omitted.
The lemma below will justify Stokes' theorem there by showing that this singular
frontier gives no contribution.

If $Y$ is an oriented smooth $k$-dimensional submanifold of a smooth manifold $U$, with
locally finite volume, its integration current is
\[
 [Y](\beta):=\int_Y\beta,
 \qquad \beta\in C_c^\infty(U;\Lambda^kT^*U).
\]
The boundary of a $k$-current $T$ is defined by $\partial T(\eta):=T(\dd\eta)$.

This notation agrees with the earlier cycle notation through integration:
to an analytic cycle $\gamma=\sum_i m_i[Z_i]$ we associate the current
$\beta\mapsto\sum_i m_i\int_{(Z_i)_{\mathrm{reg}}}\beta$, using the complex
orientation. We use the same symbol for a cycle and its associated current when
the meaning is clear. For a real manifold $Y$, such as the real slices used
below, $[Y]$ denotes its integration current only. 

\begin{lemma}\label{lem:gmt-removable-frontier}
Let $U$ be a smooth Riemannian manifold, let $\Sigma\subset U$ be closed, and let
$Y\subset U\setminus\Sigma$ be an oriented smooth $k$-dimensional submanifold with $k\geq1$, without
manifold boundary. Assume that $Y$ is closed in $U\setminus\Sigma$, has locally finite
$\mathcal H^k$-measure, and that $\mathcal H^{k-1}(\Sigma)=0$.
Then $[Y]$ is a locally integral current on $U$ and $\partial[Y]=0$.
\end{lemma}

\begin{proof}[Proof sketch]
We use the measure support theorem for locally flat currents
\cite[Theorem~2.4.2]{king1971currents}: a locally flat $r$-current
supported on a set of zero $\mathcal H^r$-measure vanishes.

The current $[Y]$ is integer rectifiable and has locally finite mass.
Exhausting $Y$ by compact smooth domains approximates $[Y]$ locally in mass,
hence in flat norm, by their integration currents. Thus $[Y]$ is locally
flat. The point of this regularity is that $\partial[Y]$ is locally flat
as well, without requiring a prior bound on its mass.

Ordinary Stokes' theorem away from $\Sigma$ gives
$\supp(\partial[Y])\subset\Sigma$: compactly supported test forms there
meet a closed smooth submanifold without boundary. The support theorem,
applied in local coordinates with $r=k-1$, now gives $\partial[Y]=0$.
Since $[Y]$ is integer rectifiable with locally finite mass and zero
boundary, it is locally integral.
\end{proof}

The geometric principle is that a singular frontier of real codimension
at least two cannot contribute a boundary to the integration current,
provided the regular part has locally finite volume. More precisely,
under the hypotheses above, a subanalytic set $\Sigma$ of dimension at
most $k-2$ has $\mathcal H^{k-1}(\Sigma)=0$, so the lemma applies.
This is why complex singularities fit the argument: a loss of one complex
dimension is a loss of two real dimensions. Allowing one real parameter
adds one dimension to both the regular part and its singular frontier,
and preserves this gap. We will verify these conditions for the real
slices of the graph in \Cref{lem:closed-real-slice}.

\subsection{Length spaces}\label{sec:length-spaces}

We recall in this section some basic notions from metric geometry that will be used later. Let $(X, \mathsf{d})$ be a metric space. If $A\subset X$ is a subspace, we denote by $\mathsf d_A:=\mathsf d|_{A\times A}$ the restricted metric. Following \cite{burago}, for a continuous curve $\eta:[a, b]\to (X, \mathsf{d})$, we define its length by:
\[L_{\mathsf{d}}(\eta):=\sup\left\{\sum_{i=0}^{N-1}\mathsf{d}(\eta(t_i), \eta(t_{i+1}))\,\bigg|\, a=t_0<\cdots <t_N=b\right\}.\]
A curve is called \emph{rectifiable} if its length is finite. 
\begin{definition}\label{def:innerdistance}Let $(X, \mathsf{d})$ be a metric space. We define the \emph{inner distance} as the function $\mathsf{d}^{\mathrm{in}}:X\times X\to [0, \infty]$, defined by
\[
\mathsf{d}^{\mathrm{in}}(x,y) : = 
\inf \left\{L_{\mathsf{d}}(\eta) \, 
\bigg|\, \eta : [a,b] \stackrel{C^0}{\longrightarrow} X, \, \eta(a) = x, \, \eta(b) = y \right\}. 
\]
\end{definition}
Note that  
$\mathsf{d}^{\mathrm{in}}$ is a possibly infinite distance over $X$
(connected components of the topology induced by the inner distance are distinguished by the finiteness of the inner distance).

In the sequel it will be convenient to use absolutely continuous curves. Recall that a curve \(\eta:I\to X\) belongs to \(\mathsf{AC}^p(I;X)\),
\(1 \leq p<\infty\), and is called \emph{absolutely continuous} if there exists
\(g\in L^p(I)\) such that
\begin{equation}\label{eq:AC}
\mathsf d(\eta(s),\eta(t))\le \int_s^t g(r)\,dr
\qquad\text{for all } s\le t \text{ in } I.
\end{equation}
We call any such function $g$ an \emph{admissible control} for $\eta$.
For every absolutely continuous curve, the metric derivative
\begin{equation}\label{eq:metricspeed}
|\dot\eta_t|_{\mathsf d}:=\lim_{h\to 0}\frac{\mathsf d(\eta(t+h),\eta(t))}{|h|}
\end{equation}
exists for almost every \(t\in I\) \cite[Theorem~1.1.2]{ags} and the function
\(t\mapsto |\dot\eta_t|_{\mathsf d}\) is itself admissible in the above inequality, that is:
\[
\mathsf d(\eta(s),\eta(t))\le \int_s^t |\dot\eta_r|_{\mathsf d}\,dr
\qquad\text{for all } s\le t.
\]
If $\eta \in \mathsf{AC}^p(I;X)$ 
then its length has the following expression:
$
L_{\mathsf{d}}(\eta) = \int_a^b |\dot\eta_t|_{\mathsf d} \,dt. 
$

On a compact parameter interval, every absolutely continuous curve is rectifiable. Conversely,
every rectifiable curve has a $1$-Lipschitz arc-length parametrization with the same endpoints,
image, and length (the constant case being immediate); see
\cite[Proposition~2.5.9]{burago}. Since Lipschitz curves are absolutely continuous, the infimum
in \Cref{def:innerdistance} is unchanged if one restricts it to rectifiable, Lipschitz, or
absolutely continuous curves.

We shall use without further mention four elementary consequences of the partition definition of
length and of \eqref{eq:AC}: length is invariant under a monotone homeomorphic change of
parameter; length is additive under a finite subdivision; finitely many absolutely continuous
pieces with matching endpoints glue to an absolutely continuous curve; and an estimate
\[
 \mathsf{d}(\sigma(s),\sigma(t))\leq\int_s^t g(r)\,\dd r,
 \qquad 0<s<t\leq\varepsilon,
\]
with $g\in L^1(0,\varepsilon)$, remains valid at $s=0$ whenever
$\sigma(s)\to\sigma(0)$. The last assertion follows simply by letting $s\to 0$.

The following result is folklore.

\begin{proposition}\label{propo:folk}
Let $(X,\mathsf d)$ be a complete metric space and assume that the inner
distance $\mathsf d^{\mathrm{in}}$ is finite. Then
$(X,\mathsf d^{\mathrm{in}})$ is complete.
If, in addition, $(X,\mathsf d)$ is compact, then
$(X,\mathsf d^{\mathrm{in}})$ is geodesic.
\end{proposition}

\begin{proof}
Let $(x_n)$ be a $\mathsf d^{\mathrm{in}}$-Cauchy sequence. Since
$\mathsf d\leq\mathsf d^{\mathrm{in}}$, completeness of $(X,\mathsf d)$
gives a limit $x$ in the original metric. Choose a subsequence with
$\mathsf d^{\mathrm{in}}(x_{n_k},x_{n_{k+1}})<2^{-k}$ and rectifiable
curves $\gamma_k$ joining consecutive terms, with
$L_{\mathsf d}(\gamma_k)<2^{1-k}$. Concatenate them on consecutive
intervals accumulating at $1$. The lengths of the tails tend to zero,
so the concatenation extends continuously to $1$, with endpoint $x$.
Each tail, including this endpoint, is rectifiable and gives
\[
 \mathsf d^{\mathrm{in}}(x_{n_k},x)
 \leq\sum_{j=k}^{\infty}L_{\mathsf d}(\gamma_j)
 \leq 2^{2-k}\longrightarrow0.
\]
The Cauchy property then gives convergence of the whole sequence in
$\mathsf d^{\mathrm{in}}$, proving completeness.

If $(X,\mathsf d)$ is compact, finiteness of the inner distance ensures
that every pair of points can be joined by a rectifiable curve.
By \cite[Proposition~2.5.19]{burago}, a shortest such curve exists.
Every subpath is also shortest; hence its constant-speed parametrization
is a geodesic for $\mathsf d^{\mathrm{in}}$.
\end{proof}

\subsection{The Wasserstein space}

We now return to the fixed compact K\"ahler manifold \(M\). Its K\"ahler form \(\omega\)
determines the smooth Riemannian metric
\[
 g(u,v):=\omega(u,Jv)
\]
and the associated Riemannian distance \(\dM\). Thus \((M,\dM)\) is, in
particular, a compact metric space. Let \(\cP(M)\) be the space of Borel probability measures on
\(M\). We use the \emph{narrow topology}: a sequence
\((\mu_k)_{k\in\N}\subset \cP(M)\) converges narrowly to \(\mu\in \cP(M)\) if and only if
\[
\int_M \varphi\,\dd\mu_k \longrightarrow \int_M \varphi\,\dd\mu
\qquad\text{for every }\varphi\in C^0(M).
\]
Thus a curve $t\mapsto\mu_t$ is \emph{narrowly continuous} precisely when
\(t\mapsto\int_M\varphi\,\dd\mu_t\) is continuous for every \(\varphi\in C^0(M)\).

\begin{definition}\label{def:wasserstein-distance}
For \(\mu,\nu\in \cP(M)\), let \(\Pi(\mu,\nu)\) denote the set of probability measures on
\(M\times M\) whose first and second marginals are \(\mu\) and \(\nu\). The \emph{$2$-Wasserstein
distance} between \(\mu\) and \(\nu\) is defined by\footnote{This is usually denoted by $\W_2$, but since we will only be concerned with the $2$-Wasserstein distance, we omit the subscript from the notation.}
\[
\W(\mu,\nu):=
\left(\inf_{\pi\in \Pi(\mu,\nu)}
\int_{M\times M} \dM(x,y)^2\,\dd\pi(x,y)\right)^{\frac{1}{2}}.
\]
\end{definition}

Since \(M\) is compact, the cost function \(\dM^2\) is bounded. In particular,
\(\W(\mu,\nu)\) is finite for every \(\mu,\nu\in \cP(M)\). Moreover, \((\cP(M),\W)\) is a
compact metric space, and the topology induced by \(\W\) agrees with the
narrow topology \cite[Chapter~7, Section~7.1]{ags}.

\subsection{Tangent vectors and metric speed in Wasserstein space}\label{sec:wasserstein-tangents}
We now recall the part of the infinitesimal structure of \((\cP(M),\W)\) that will be used later;
see \cite[Chapter~8]{ags}. Let \(t\mapsto \mu_t\) be an absolutely continuous curve in the metric space \((\cP(M),\W)\). Then there exists a
Borel time-dependent vector field
\[
v_t\in L^2(\mu_t;TM)
\]
such that
\begin{equation}\label{eq:continuity}
\partial_t\mu_t+\operatorname{div}(v_t\mu_t)=0
\end{equation}
in the sense of distributions \cite[Chapter~8]{ags}. In general such a field is not unique.
However, among all velocity fields generating the same curve, there is a unique one of minimal
\(L^2(\mu_t)\)-norm. Moreover, for almost every \(t\in I\), this distinguished field recovers
exactly the metric speed:
\[
|\dot\mu_t|_{\W}=\|v_t\|_{L^2(\mu_t)}
\]
\cite[Chapter~8 and Proposition~8.4.5]{ags}.

This leads to the infinitesimal description of the Wasserstein space. For \(\mu\in \cP(M)\), the
Wasserstein tangent space at \(\mu\) is defined as the \(L^2(\mu)\)-closure of gradients of smooth functions on
\(M\), namely
\[
T_\mu\cP(M):=\overline{\{\nabla \phi:\phi\in C^\infty(M)\}}^{\,L^2(\mu)},
\]
where $\nabla \phi$ denotes the Riemannian gradient of the function $\phi$.
If \(t\mapsto \mu_t\) is an absolutely continuous curve in \((\cP(M),\W)\), then the distinguished
velocity field of minimal \(L^2(\mu_t)\)-norm belongs to \(T_{\mu_t}\cP(M)\) for almost every
\(t\), and it is therefore natural to regard it as the tangent vector to the curve \((\mu_t)\) at
time \(t\). This reflects the infinitesimal behavior of optimal transport for the quadratic cost:
among all velocity fields producing the same first-order variation of \(\mu\), the gradients are the
distinguished ones \cite[Section~8.4]{ags}.

When the measure is supported on a smooth closed submanifold, the tangent space admits a more concrete
description. \begin{proposition}[{\cite[Proposition~2.8]{ACL}}]Let \(S\subset M\) be a smooth closed submanifold, let \(TS^\perp\subset TM|_S\) be its Riemannian normal
bundle, and let \(\mu\in \cP(M)\) satisfy \(\supp(\mu)\subset S\). Then
\begin{equation}\label{eq:tangentsplitting}
T_\mu\cP(M)
=
\overline{\{\nabla_S f:f\in C^\infty(S)\}}^{\,L^2(\mu)}
\oplus
L^2(TS^\perp,\mu),
\end{equation}
where \(\nabla_S\) denotes the gradient on \(S\) with respect to the induced Riemannian metric. In
other words, every Wasserstein tangent vector at such a measure splits uniquely into a tangential
gradient part and a normal part. 
\end{proposition}

We will also use the converse direction of the continuity-equation
description: a velocity field that transports the measures provides an
upper bound for their Wasserstein speed. The following lemma is the
converse part of Erbar's criterion \cite[Proposition~2.5]{Erbar}, together
with its immediate consequence for length. We distinguish the local
square-integrability required by the criterion from the integrability
needed for finite total length.

\begin{lemma}[Continuity equation and speed]\label{lem:continuity-speed-criterion}
Let $M$ be a compact Riemannian manifold, let $I\subset\R$ be a bounded
open interval, and let $t\mapsto\mu_t\in\cP(M)$ be narrowly continuous.
Suppose that $v:I\times M\to TM$ is Borel, with $v_t(x)\in T_xM$, and put
\[
 A(t):=\|v_t\|_{L^2(\mu_t)}.
\]
Assume $A\in L^2_{\mathrm{loc}}(I)$ and
\begin{equation}\label{eq:field-continuity-weak}
 \int_I\int_M
 \bigl(\partial_t\varphi+\langle\nabla_M\varphi,v_t\rangle_g\bigr)
 \,\dd\mu_t\,\dd t=0
 \qquad\bigl(\varphi\in C_c^\infty(I\times M)\bigr).
\end{equation}
Then $\mu\in\mathsf{AC}^2_{\mathrm{loc}}(I;(\cP(M),\W))$ and
\[
 |\dot\mu_t|_\W\leq A(t)\quad\text{for almost every }t\in I,
 \qquad
 \W(\mu_s,\mu_t)\leq\int_s^t A(r)\,\dd r\quad(s<t\text{ in }I).
\]
If also $A\in L^1(I)$, then $L_\W(\mu;I)\leq\int_I A(t)\,\dd t<\infty$.
\end{lemma}

\begin{proof}
The converse part of \cite[Proposition~2.5]{Erbar} gives local absolute
continuity and the speed bound on the Riemannian manifold $M$.
Integration of the metric speed gives the distance estimate. Summing
that estimate over partitions gives the length bound.
\end{proof}

\section{Wasserstein geometry of the cycle space}\label{sec:wasserstein-cycle-geometry}

\subsection{The cycle-to-measure map}
We now pass from cycles to probability measures. 
This is the central construction of the paper: once cycles are embedded into Wasserstein
space through their normalized volume measures, one can transfer optimal transport geometry
back to the cycle space itself and use it to define intrinsic distances on Barlet components. The embedding procedure is produced by the following map.

\begin{definition}\label{def:mu-gamma}

Let $(M, \omega)$ be a compact K\"ahler manifold and let $B\subset \mathcal{C}_q(M)$ be a connected component of the Barlet space. We define the \emph{cycle-to-measure map}
\[
\mug\colon B\longrightarrow \cP(M)
\]
as follows. If
\[
\gamma=\sum_{j=1}^N m_j[Z_j]\in B,
\]
where the \(Z_j\) are the irreducible components of \(|\gamma|\) and \(m_j\in \mathbb N\), then
\(\mug(\gamma)\) is the  measure characterized by
\[
\int_M \varphi\, d\mug(\gamma)
=
\frac{1}{V_B}
\sum_{j=1}^N m_j \int_{(Z_j)_{\mathrm{reg}}} \varphi\,\frac{\omega^q}{q!}
\qquad \forall\,\varphi\in C^0(M).
\]
Observe that this is a probability measure thanks to \eqref{eq:int_of_form} and \eqref{eq:V}.
\end{definition}

The next point is that this construction is faithful: it depends
continuously on the cycle, and different cycles give different measures.

\begin{proposition}\label{prop:mu-gamma}
In the above setting, the map \(\mug\colon B\to \cP(M)\) is continuous and injective. In particular, it is a
homeomorphism onto its image.
\end{proposition}

\begin{proof}

We first prove continuity. For every continuous \((q,q)\)-form \(\alpha\) on \(M\), the map
\[
\gamma\longmapsto \int_\gamma \alpha
\]
is continuous for the topology of the cycle space
\cite[Theorem~4.2.26]{cycles1}.
 Applying this to \(\alpha=\varphi\,\frac{\omega^q}{q!}\), with
\(\varphi\in C^0(M)\), we obtain that
\[
\gamma\longmapsto \int_M \varphi\, d\mug(\gamma)
=
\frac{1}{V_B}\int_\gamma \varphi\,\frac{\omega^q}{q!}
\]
is continuous for every \(\varphi\in C^0(M)\). Hence \(\mug\) is continuous for the
narrow topology on \(\cP(M)\).

We now prove injectivity. Let \(\gamma,\eta\in B\), and assume that
\(\mug(\gamma)=\mug(\eta)\). Since the normalizing constant \(V_B\) is the same on
all of \(B\), this means that
\[
\int_\gamma \varphi\,\frac{\omega^q}{q!}=\int_{\eta} \varphi\,\frac{\omega^q}{q!}
\qquad \forall\,\varphi\in C^0(M).
\]
In other words, \(\gamma\) and \(\eta\) define the same positive Radon measure on \(M\). Write
\[
\gamma=\sum_{i=1}^N m_i[Y_i],
\qquad
\eta=\sum_{j=1}^P n_j[W_j].
\]
Since \(\omega^q\) restricts to a strictly positive volume form on the smooth locus of every
irreducible \(q\)-dimensional analytic set, the support of the measure associated with a cycle is
exactly the support of the cycle. Therefore
\[
|\gamma|=|\eta|.
\]

Let \(Y\) be an irreducible component of this common support, and choose a smooth point of \(Y\)
which does not belong to any other irreducible component. In a sufficiently small neighborhood of
such a point, both cycles reduce to multiples of \([Y]\), say \(m[Y]\) and \(n[Y]\). Equality of the
corresponding measures on that neighborhood then gives
\[
m\int_Y \varphi\,\frac{\omega^q}{q!}
=
n\int_Y \varphi\,\frac{\omega^q}{q!}
\]
for every nonnegative \(\varphi\in C^0(M)\) supported in that neighborhood. Since the integral is
strictly positive for some such \(\varphi\), it follows that \(m=n\). Repeating the argument for
every irreducible component, we conclude that \(\gamma=\eta\).

Finally, \(B\) is compact by the compact K\"ahler hypothesis, while \(\cP(M)\) is Hausdorff for the
narrow topology. Therefore
$\mug\colon B\to\cP(M)$, which we have just shown to be continuous and injective,
is a homeomorphism onto its image.
\end{proof}

\subsection{The inner Wasserstein distance}
We now use the cycle-to-measure map to measure deformations of cycles by their
Wasserstein length. Thus a subset of the Barlet space inherits an intrinsic
distance from its image in \((\cP(M),\W)\).

\begin{definition}
\label{def:inner-wasserstein-distance}
Let \(A\subset B\) be a connected subset of the fixed nonzero Barlet component. We define the
 \emph{outer Wasserstein distance} on \(A\) by
\[
\W_A(\gamma_0,\gamma_1)
:=
\W\bigl(\mug(\gamma_0),\mug(\gamma_1)\bigr).
\]
The \emph{inner Wasserstein distance} on \(A\) is the extended inner distance
\[\W_A^{\mathrm{in}}\colon A\times A\to[0,\infty]\] associated to $\W_A$  in the sense of \cref{def:innerdistance}.
\end{definition}

\begin{remark}Clearly, the inner distance need not coincide with the outer distance. For example, on a smooth submanifold \(A\subset\mathbb R^n\) it gives the induced Riemannian distance, whereas the outer Euclidean distance is the chordal one. Similarly, here
\[
\W_A(\gamma_0,\gamma_1)\le
\W_A^{\mathrm{in}}(\gamma_0,\gamma_1),
\]
and the inequality may be strict, since the ambient Wasserstein geodesic need not remain in \(\mug(A)\), hence need not come from a deformation through cycles in \(A\). This phenomenon disappears for components of the space of zero-dimensional cycles; see \Cref{sec:zerocycles}.
\end{remark}

\subsection{Finiteness of the inner Wasserstein distance}

The aim of this subsection is to prove the following.
\begin{theorem}\label{thm:finiteness}
Let $A\subset B$ be a connected analytic subset. Then the inner Wasserstein distance
$\W_A^{\mathrm{in}}$ is finite. In particular this holds for $A=B$.
\end{theorem}

The assertion concerns all of $A$, including its singular and nonreduced
cycles; its good locus may even be empty. By
\Cref{cor:piecewise-analytic-connectivity}, any two points of $A$ can be
joined by finitely many real-analytic arcs. The remaining issue is to
show that these arcs have finite Wasserstein length, especially near
parameters where their cycles become singular or acquire multiplicities.

This is a local question in the parameter. A real-analytic arc extends
locally to a holomorphic family, and we may translate any parameter under
consideration to $t=0$. We thus study a holomorphic arc
$\alpha:\Delta\to B$, where $\Delta$ is a disk centered at zero, allowing
$\alpha(0)$ to be singular or nonreduced. Writing
$\mu_t:=\mug(\alpha(t))$, our goal is to prove that there are $C>0$,
$\varepsilon>0$ and $\theta\in(0,1)$ such that
\[
 L_{\W}(\mu;[0,\varepsilon])
 \leq C\int_0^\varepsilon r^{-\theta}\,\dd r<\infty.
\]
The integrability of $r^{-\theta}$ is what allows the path to reach the
possibly degenerate cycle $\alpha(0)$ with finite length.
Let us explain how we obtain this estimate.

For a smoothly moving submanifold, one can describe its motion by a
velocity field and use that field to estimate the Wasserstein speed of
its volume measure. Here the fibers may be singular, so there is no
globally smooth family on which to carry out this construction.
Instead, we work with the graph
\[
 \{(t,x)\in\Delta\times M:x\in|\alpha(t)|\},
\]
which records the moving supports. We use the regular part of this
\emph{whole reduced graph}, thereby excluding both singular points of
its components and their intersections. After shrinking the disk, the
remaining pieces are disjoint smooth families for $t\ne0$, and the
excluded set has zero volume in every such fiber. The pieces therefore
recover the entire cycle measures when weighted by the graph-cycle
multiplicities (\Cref{ss:arc-graph}).

On each piece we construct a tangent field $H_a=(1,\nu_a)$: its first
component advances time at unit speed, and its second component
$\nu_a(t,x)\in T_xM$ is the normal velocity of the moving fiber.
A subanalytic gradient inequality bounds $|\nu_a(t,x)|_g$ by
$Ct^{-\theta}$, with $\theta<1$. Since the pieces are disjoint, their
normal velocities define a single time-dependent Borel field $v_t(x)$
on a set of full $\mu_t$-measure. Its $L^2(\mu_t)$-norm has the same
integrable bound (\Cref{ss:arc-velocities}).

We must then verify that $v_t$ describes the variation of the
\emph{whole} cycle measure. Although the omitted set has zero fiber
volume, we must also rule out a boundary contribution there. Stokes'
theorem for the real slices of the graph does this: summing the identities
for the components, with their multiplicities, gives
$\partial_t\mu_t+\operatorname{div}(v_t\mu_t)=0$.
The continuity-equation criterion then shows that $\mu$ is locally
absolutely continuous for $t>0$, with
\[
 |\dot\mu_t|_\W\leq\|v_t\|_{L^2(\mu_t)}\leq Ct^{-\theta}
 \qquad\text{for almost every }t>0.
\]
Integrating this speed bound controls the length on every interval
$[s,\varepsilon]$ with $s>0$, uniformly as $s\downarrow0$.
Continuity of the cycle-to-measure map gives
$\mu_s\to\mug(\alpha(0))$, so the bound extends to the closed interval
$[0,\varepsilon]$ and yields the desired finite-length estimate
(\Cref{ss:arc-continuity}).

Finally, applying the same argument on both sides of each parameter
gives finite Wasserstein length locally along every real-analytic arc.
Compactness of its parameter interval then gives finite length for the
whole arc. The finitely many arcs joining the prescribed endpoints
therefore concatenate to a path of finite Wasserstein length, proving
finiteness of the inner distance (\Cref{ss:arc-finiteness}).

\subsubsection{The graph of a holomorphic family}\label{ss:arc-graph}
Let $\alpha:\Delta_0\to B$ be holomorphic, with $\Delta_0$ a disk centered
at zero. Its graph-cycle is the $(q+1)$-dimensional analytic cycle
\[
 \mathscr G_\alpha=\sum_a m_a[\widehat G_a],\qquad m_a\geq1,
\]
on $\Delta_0\times M$. Denote its reduced support by
$\widehat G:=\bigcup_a\widehat G_a$, and denote the two projections by
$\pi:\Delta_0\times M\to\Delta_0$ and $j:\Delta_0\times M\to M$.
Each irreducible component $\widehat G_a$ is proper over $\Delta_0$ and
dominates the disk, meaning that $\pi(\widehat G_a)=\Delta_0$.
Indeed, its image is a closed analytic subset by Remmert's theorem, hence
either the disk or a point. The latter is impossible: $\widehat G_a$ has
dimension $q+1$, whereas a parameter fiber of $\widehat G$ has dimension $q$.
All components meet the compact central fiber, so local finiteness of the
irreducible decomposition makes their number finite. Write $a=1,\ldots,N$.

After choosing a smaller disk $\Delta\Subset\Delta_0$, set
\[
 G:=\widehat G\cap\pi^{-1}(\Delta),\qquad
 G_a:=\widehat G_a\cap\pi^{-1}(\Delta),\qquad S:=\operatorname{Sing}(G).
\]
The index $a$ continues to label an original irreducible component
$\widehat G_a$; its restriction $G_a$ may split into several pieces,
which we keep together under that label. Arguments using irreducibility
below are made before this restriction.

We work on the regular part of the whole reduced support $G$. Thus $S$
includes both the singularities of individual components and all their
pairwise intersections. Put
\[
 \Delta^*:=\Delta\setminus\{0\},\qquad
 G_a^\circ:=(G_a\setminus S)\cap\pi^{-1}(\Delta^*),
 \qquad
 G_{a,t}:=G_a\cap\pi^{-1}(t),\quad
 G_{a,t}^\circ:=G_a^\circ\cap\pi^{-1}(t).
\]
At every point of $G\setminus S$ exactly one component is present.
Consequently the sets $G_a^\circ$ are pairwise disjoint, and each is an
open submanifold of $(G_a)_{\mathrm{reg}}$.

\begin{lemma}\label{lem:arc-incidence-model}
The disk $\Delta$ can be chosen so that, for every $a$ and $t\in\Delta^*$:
\begin{enumerate}[label=(\roman*)]
\item $\pi:G_a^\circ\to\Delta^*$ is a holomorphic submersion with nonempty
$q$-dimensional fibers;
\item the fiber cycle is reconstructed by
\begin{equation}\label{eq:arc-fiber-decomposition}
 \alpha(t)=\sum_{a=1}^N m_a\,(j|_{G_{a,t}})_*[G_{a,t}]_{\mathrm{red}};
\end{equation}
\item $\dim_\C(S\cap G_{a,t})\leq q-1$. In particular,
$G_{a,t}^\circ$ is a dense open submanifold of the pure $q$-dimensional
fiber $G_{a,t}$.
For $q\geq1$, its complement has real dimension at most
$2q-2$, hence zero $\mathcal H^{2q-1}$-measure and, in particular,
zero $2q$-dimensional volume. For $q=0$, the complement is empty.
\end{enumerate}

\end{lemma}
Here $[G_{a,t}]_{\mathrm{red}}$ assigns multiplicity one to every
irreducible component of the fiber, as in
\eqref{eq:reduced-fundamental-cycle}. The map $j|_{G_{a,t}}$ identifies
that fiber with an analytic subset of $M$; its push-forward places the
reduced cycle in $M$, while $m_a$ retains the graph-cycle multiplicity.
\begin{proof}
Every fiber of $\widehat G_a$ is pure $q$-dimensional. Indeed, intersecting
the pure $(q+1)$-dimensional space with the nonconstant parameter
hypersurface gives dimension at least $q$
\cite[Proposition~7.1.1.3]{cycles2}, while containment in
$\{t\}\times|\alpha(t)|$ gives dimension at most $q$.

Consider on the original support $\widehat G$ the closed analytic sets
\[
 \widehat S:=\operatorname{Sing}(\widehat G),\qquad
 \widehat R:=\widehat G\setminus
 \{z\in\widehat G_{\mathrm{reg}}:
       \dd\pi_z|_{T_z\widehat G}\ne0\}.
\]

To see why $\widehat R$ is analytic, write
$T_z^{\mathrm{Zar}}\widehat G$ for the Zariski tangent space, obtained by
linearizing local holomorphic defining equations of $\widehat G$ at $z$.
Since $\widehat G$ is reduced and pure $(q+1)$-dimensional, the Jacobian
criterion says that this tangent space has dimension $q+1$ precisely at
regular points. At a singular point its dimension is therefore at least
$q+2$. The differential $\dd\pi_z$ takes values in the one-dimensional
space $T_{\pi(z)}\Delta_0$, so its kernel on
$T_z^{\mathrm{Zar}}\widehat G$ has dimension at least $q+1$ at every
singular point. At a regular point, the same kernel has dimension at
least $q+1$ exactly when $\dd\pi_z=0$. Consequently
\[
 \widehat R
 =\bigl\{z\in\widehat G:
   \dim_\C\ker\bigl(\dd\pi_z|_{T_z^{\mathrm{Zar}}\widehat G}\bigr)
        \geq q+1\bigr\}.
\]
In a local embedding, this is a rank condition on the matrix obtained
by adjoining the row $\dd\pi$ to the Jacobian of the defining equations.
It is thus defined by the vanishing of holomorphic minors, so
$\widehat R$ is a closed analytic subset; this is the relative
singular-locus criterion of \cite[\S2.17]{Fischer1976}.

Let $W$ be an irreducible component of $\widehat R$ not contained in
$\widehat S$. Since $W\cap\widehat S$ is a proper analytic subset of
$W$, the set $W_{\mathrm{reg}}\setminus\widehat S$ is open and dense in
$W_{\mathrm{reg}}$. At its points, $\widehat G$ is regular and the
definition of $\widehat R$ gives $\dd\pi_z|_{T_z\widehat G}=0$.
Because $T_zW\subset T_z\widehat G$, the differential
$\dd(\pi|_{W_{\mathrm{reg}}})$ vanishes there as well, and hence vanishes
everywhere on $W_{\mathrm{reg}}$ by continuity. The regular part of an
irreducible analytic set is connected, so $\pi$ is constant on
$W_{\mathrm{reg}}$, and then on $W$ by continuity. Thus $W$ is
\emph{vertical}, meaning that it is contained in a single parameter
fiber.

Every irreducible component of $\widehat S$ has dimension at most $q$
and is either vertical or dominates the disk, by properness and
Remmert's theorem.

Choose a disk $\Delta_1$ centered at zero with
$\overline{\Delta_1}\subset\Delta_0$. The set
$K:=\widehat G\cap\pi^{-1}(\overline{\Delta_1})$ is compact, since it is
closed in the compact product $\overline{\Delta_1}\times M$.
The irreducible decompositions of the analytic sets $\widehat R$ and
$\widehat S$ are locally finite: every point has a neighborhood meeting
only finitely many components. A finite cover of $K$ by such
neighborhoods shows that only finitely many components of either set
meet $K$.

Each vertical component just identified lies over a single parameter
value. Let $E\subset\overline{\Delta_1}$ be the finite set of these
values for the vertical components of $\widehat R$ and $\widehat S$
that meet $K$. Choose $\delta>0$ so that the disk
$\Delta:=\{t\in\C:|t|<\delta\}$ lies in $\Delta_1$ and
$\delta<|e|$ for every $e\in E\setminus\{0\}$. If this last set is empty,
there is no additional restriction on $\delta$. No nonzero fiber over
$\Delta$ can then meet a vertical component of $\widehat R$ or
$\widehat S$; such components may remain only over the central
parameter $t=0$.

Then $\pi$ has no critical points on $G\setminus S$ over $\Delta^*$.
Every component of $\widehat S$ meeting a nonzero fiber dominates the
disk; intersecting it with that parameter hypersurface lowers its
dimension by one. Thus
$\dim_\C(S\cap\pi^{-1}(t))\leq q-1$ for $t\ne0$.
This proves (iii), including the estimate for intersections of distinct
graph components. Purity of the fibers now gives their nonempty dense
regular parts and proves (i).

Finally, the graph-cycle identity recalled in the preliminaries gives
\[
 \{t\}\times\alpha(t)
 =\sum_a m_a\bigl([G_a]\cap[\{t\}\times M]\bigr).
\]
For $t\ne0$, each irreducible component of $G_{a,t}$ has a dense open
subset in $G_a^\circ$. There the intersection with $\{t\}\times M$ is
transverse, hence has multiplicity one
\cite[Proposition~7.1.3.2]{cycles2}. Therefore
$[G_a]\cap[\{t\}\times M]=[G_{a,t}]_{\mathrm{red}}$.
Identifying the parameter fiber with $M$ proves (ii).
\end{proof}

When integrating in $x$, we henceforth identify $G_{a,t}^\circ$ with its
image in $M$ under $j$. These images are pairwise disjoint for each
$t\ne0$. Writing $\mu_t:=\mug(\alpha(t))$, the lemma gives
\begin{equation}\label{eq:arc-measure}
 \int_M f\,\dd\mu_t
 =\frac1{V_B}\sum_a m_a\int_{G_{a,t}^\circ}f(x)\,\frac{\omega^q}{q!}
 \qquad\text{for every bounded Borel }f:M\to\R.
\end{equation}
In particular, their union has full $\mu_t$-measure, and setting $f=1$
gives $\sum_a m_a\int_{G_{a,t}^\circ}\omega^q/q!=V_B$.

\subsubsection{Controlled normal velocities}\label{ss:arc-velocities}
We first construct a tangent field on the graph that advances the
parameter at unit speed. Its component in $M$ will give the normal
velocity of the moving fibers. Write the disk coordinate as
$\tau+i\sigma$, choose $\varepsilon>0$ with $[0,\varepsilon]\subset\Delta$,
and set
\begin{equation}\label{eq:arc-real-slice}
 U:=(0,\varepsilon)\times M,\qquad
 Y_a:=G_a^\circ\cap U.
\end{equation}
Here $U$ is viewed as the real slice $\sigma=0$ in $\Delta\times M$.
Each $Y_a$ is a smooth real-analytic manifold of dimension $2q+1$.
Use the product of the Euclidean metric on the disk and the K\"ahler
metric $g$ on $M$, and the induced metrics on the regular graph pieces.

Define the tangent field
\begin{equation}\label{eq:arc-horizontal-field}
 H_a(z):=\frac{\nabla^{G_a}\tau(z)}{|\nabla^{G_a}\tau(z)|^2},
 \qquad z=(t,x)\in Y_a.
\end{equation}
The denominator is nonzero by the submersion property. If $J$ denotes
the complex structure on $\Delta\times M$, the Cauchy--Riemann equations
give $\nabla^{G_a}\sigma=J\nabla^{G_a}\tau$. Thus
$\nabla^{G_a}\tau$ is tangent to $\sigma=0$, and
\begin{equation}\label{eq:arc-horizontal-identities}
 H_a(z)\in T_zY_a,\qquad
 \dd\tau(H_a)=1,\quad \dd\sigma(H_a)=0,\qquad
 |H_a|=|\nabla^{G_a}\tau|^{-1}.
\end{equation}
The lower-case $j$ still denotes the projection to $M$. Define
\begin{equation}\label{eq:arc-spatial-velocity}
 \nu_a(t,x):=\dd j_{(t,x)}H_a(t,x)\in T_xM,
 \qquad H_a(t,x)=(1,\nu_a(t,x))\in\R\times T_xM.
\end{equation}
Thus $H_a$ is a vector field on $Y_a$, whereas $\nu_a$ is a vector field
along $j|_{Y_a}$, or equivalently a section of $(j|_{Y_a})^*TM$.
At a fixed time it is a field along the embedded fiber $G_{a,t}^\circ$.
A local integral curve $(t(s),x(s))$ of $H_a$ satisfies
$t'(s)=1$ and $x'(s)=\nu_a(t(s),x(s))$, which explains its interpretation
as the spatial velocity. The symbol $z$ in \eqref{eq:arc-horizontal-field}
denotes the point $(t,x)$ of the graph.

Let $V_{a,z}:=\dd j_z(\ker\dd\pi_z\cap T_zG_a)$ be the complex tangent
plane of the fiber at $x$.
\begin{lemma}\label{lem:arc-normality}
The field $\nu_a(z)$ is orthogonal to $V_{a,z}$ and satisfies
$|\nu_a(z)|_g\leq|H_a(z)|$.
\end{lemma}
\begin{proof}
For $u\in V_{a,z}$, the vector $(0,u)$ is vertical and tangent to $G_a$.
Consequently
\[
 \langle\nu_a,u\rangle_g
 =\langle H_a,(0,u)\rangle
 =\frac{\dd\tau(0,u)}{|\nabla^{G_a}\tau|^2}=0.
\]
The product decomposition gives $|H_a|^2=1+|\nu_a|_g^2$, hence the bound.
\end{proof}

\begin{lemma}\label{lem:uniform-subanalytic-gradient}
After decreasing $\varepsilon$, there are $c>0$ and $\theta\in(0,1)$,
independent of $a$, such that
\begin{equation}\label{eq:arc-gradient-estimate}
 |\nabla^{G_a}\tau(z)|\geq c\,\tau(z)^\theta\qquad(z\in Y_a).
\end{equation}
Consequently, with $C=c^{-1}$,
\begin{equation}\label{eq:arc-velocity-estimate}
 |\nu_a(t,x)|_g\leq |H_a(t,x)|\leq Ct^{-\theta}.
\end{equation}
\end{lemma}
\begin{proof}
Use the real-analytic embedding of $M$ from
\Cref{prop:definable-tautological-diagram}, together with the disk
coordinate, to embed a relatively compact disk times $M$ into a Euclidean
space. Since the graph and its singular locus are analytic on a larger
disk, every $Y_a$ is a relatively compact definable real-analytic
submanifold. Denote its induced Euclidean metric by $h$.

The function $\tau|_{Y_a}$ is definable and extends continuously, with
value zero, at every point of the compact set
\[
 K_a:=\overline{Y_a}\cap(\{0\}\times M).
\]
Apply the boundary version of the gradient inequality in
Property~\ref{item:seven} at each point $z_0\in K_a$. It gives local
bounds $|\nabla_h^{Y_a}\tau|_h\geq c_{z_0}\tau^{\theta_{z_0}}$,
with $\theta_{z_0}<1$. A finite cover of $K_a$ makes these bounds uniform
near $K_a$: take $\tau<1$, choose $\theta\in(0,1)$ larger than the finitely
many exponents, and take the smallest constant.

All points of $Y_a$ with sufficiently small positive $\tau$ lie in this
neighborhood. Otherwise a sequence outside it with $\tau\to0$ would
have a limit in $K_a$, by compactness of $M$. There are finitely many
indices $a$, so the constants and exponent can be chosen simultaneously.

On the fixed compact ambient set the Euclidean and product metrics are
uniformly equivalent. Their dual norms on every tangent subspace are
therefore uniformly equivalent as well, giving the same gradient bound
for the product metric on $Y_a$. This comparison does not require $g$
to be real analytic. Finally, $\nabla^{G_a}\tau$ is tangent to $Y_a$,
so it equals $\nabla^{Y_a}\tau$ for that metric. This proves
\eqref{eq:arc-gradient-estimate}; the velocity bound follows from
\eqref{eq:arc-horizontal-identities} and \Cref{lem:arc-normality}.
\end{proof}

We can now assemble the spatial velocities into an ordinary field on
$M$ depending on time. Define
\begin{equation}\label{eq:arc-global-field}
 v_t(x):=
 \begin{cases}
  \nu_a(t,x),& (t,x)\in Y_a,\\
  0\in T_xM,& (t,x)\notin\bigcup_aY_a.
 \end{cases}
\end{equation}
There is no choice of component in this definition because the $Y_a$
are disjoint. Each is Borel and each $\nu_a$ is smooth there, so
$(t,x)\mapsto v_t(x)$ is Borel. For every $t\in(0,\varepsilon)$ it
agrees with the normal velocities on a set of full $\mu_t$-measure.
By \eqref{eq:arc-measure} and \eqref{eq:arc-velocity-estimate},
\begin{equation}\label{eq:arc-field-energy}
 \|v_t\|_{L^2(\mu_t)}^2
 =\frac1{V_B}\sum_a m_a\int_{G_{a,t}^\circ}
       |\nu_a(t,x)|_g^2\,\frac{\omega^q}{q!}
 \leq C^2t^{-2\theta}.
\end{equation}
The coefficients $m_a$ weight the measures; they do not create additional
velocity values. It remains to prove that this field satisfies the
continuity equation for the whole family of measures.

\subsubsection{The continuity equation and the arc estimate}\label{ss:arc-continuity}
The field $v_t$ is now defined and its $L^2(\mu_t)$-norm is at most
$Ct^{-\theta}$. To use \Cref{lem:continuity-speed-criterion}, we must
show that it describes the variation of $\mu_t$ across the omitted
singularities and component intersections. We first prove Stokes' theorem
on each real slice, then sum the resulting identities with the weights $m_a$.

With $U$ and $Y_a$ as in \eqref{eq:arc-real-slice}, set
\[
 \Sigma_a:=S\cap G_a\cap U.
\]
This includes the intersections with other graph components, as well as
the singularities of $G_a$ itself.
\begin{lemma}\label{lem:closed-real-slice}
Orient $Y_a$ by $\dd\tau\wedge j^*(\omega^q/q!)>0$. Then $T_a:=[Y_a]$
is a locally integral $(2q+1)$-current on $U$ and $\partial T_a=0$.
\end{lemma}
\begin{proof}
The submersion property makes $Y_a$ smooth without manifold boundary.
Its vertical tangent spaces are complex $q$-planes, and
$\dd\tau(H_a)=1$, so the displayed form is a positive volume form.
The set $\Sigma_a$ is closed in $U$, and
$Y_a=(G_a\cap U)\setminus\Sigma_a$ is closed in $U\setminus\Sigma_a$.
Relative compactness and subanalyticity give finite $(2q+1)$-dimensional
volume, also for the given metric by uniform metric comparison.

For $q\geq1$, \Cref{lem:arc-incidence-model}(iii) and the subanalytic
fiber-dimension inequality give
\[
 \dim_\R\Sigma_a\leq 2(q-1)+1=2q-1,
 \qquad \mathcal H^{2q}(\Sigma_a)=0.
\]
For $q=0$, the set $\Sigma_a$ is empty. In both cases
\Cref{lem:gmt-removable-frontier}, with $k=2q+1$, proves the assertion.
\end{proof}

Before deriving the variation of the fiber measures, we record an elementary
K\"ahler identity. It shows that only the normal component of a velocity contributes
when differentiating the integral of a test function over a complex fiber.
Here $\zeta$ is an arbitrary vector at $x$; in the next proof
we apply the identity with $\zeta=\nu_a(t,x)$.

\begin{lemma}
\label{lem:kahler-normal-contraction}
Let $x\in M$, let $P\subset T_xM$ be a complex $q$-plane, and write
$\zeta=\zeta^T+\zeta^\perp$ according to the orthogonal decomposition
$T_xM=P\oplus P^\perp$. Then, for every $\varphi\in C^\infty(M)$,
\begin{equation}\label{eq:kahler-normal-contraction}
 \left.\iota_\zeta(\dd\varphi\wedge\frac{\omega^q}{q!})\right|_P
 =
 \langle\nabla\varphi,\zeta^\perp\rangle_g\left.\frac{\omega^q}{q!}\right|_P.
\end{equation}
\end{lemma}

\begin{proof}
For $q=0$, both sides equal $\dd\varphi(\zeta)$. Assume $q\geq1$.

The tangential part gives no contribution, because
$\iota_{\zeta^T}(\dd\varphi\wedge\frac{\omega^q}{q!})|_P$ evaluates a $(2q+1)$-form on vectors contained in
the real $2q$-plane $P$. Since $P$ is $J$-invariant, so is $P^\perp$, and therefore
$\omega(\zeta^\perp,v)=0$ for every $v\in P$. It follows that
$\iota_{\zeta^\perp}\frac{\omega^q}{q!}|_P=0$. Restricting
\[
 \iota_{\zeta^\perp}(\dd\varphi\wedge\frac{\omega^q}{q!})
 =\dd\varphi(\zeta^\perp)\frac{\omega^q}{q!}
  -\dd\varphi\wedge\iota_{\zeta^\perp}\frac{\omega^q}{q!}
\]
to $P$ proves \eqref{eq:kahler-normal-contraction}.
\end{proof}

\begin{lemma}\label{lem:branchwise-variation}
For every $\varphi\in C_c^\infty((0,\varepsilon)\times M)$,
\begin{equation}\label{eq:arc-branchwise-variation}
 \int_0^\varepsilon\int_{G_{a,t}^\circ}
 \bigl(\partial_t\varphi(t,x)+\dd_M\varphi(t,x)[\nu_a(t,x)]\bigr)
 \,\frac{\omega^q}{q!}\,\dd t=0.
\end{equation}
\end{lemma}

\begin{proof}
Since $\partial T_a=0$ and $\dd(\omega^q/q!)=0$,
\[
 0=T_a\bigl(\dd(\varphi j^*(\omega^q/q!))\bigr)
  =\int_{Y_a}\dd\varphi\wedge j^*(\omega^q/q!).
\]
At $(t,x)$ evaluate the integrand on $H_a=(1,\nu_a)$ followed by an
oriented basis of the vertical tangent plane. Its temporal contribution
is $\partial_t\varphi\,\omega^q/q!$. Its spatial contribution is
$\iota_{\nu_a}(\dd_M\varphi\wedge\omega^q/q!)$ restricted to that plane.
By normality and \Cref{lem:kahler-normal-contraction}, their sum is
$(\partial_t\varphi+\dd_M\varphi[\nu_a])\,\omega^q/q!$.
Since $\dd\tau(H_a)=1$, integration along the submersion
$\tau:Y_a\to(0,\varepsilon)$ gives
\eqref{eq:arc-branchwise-variation}.
These integrals are absolutely convergent: the velocity bound is uniform
on the compact time support of $\varphi$, and
\eqref{eq:arc-measure} bounds the fiber volumes.
\end{proof}

Multiply \eqref{eq:arc-branchwise-variation} by $m_a/V_B$ and sum over $a$.
The measure identity \eqref{eq:arc-measure} and the definition of $v_t$
give the continuity equation for the entire cycle measures:
\begin{equation}\label{eq:arc-continuity-equation}
 \int_0^\varepsilon\int_M
 \bigl(\partial_t\varphi(t,x)+\langle\nabla_M\varphi(t,x),v_t(x)\rangle_g\bigr)
 \,\dd\mu_t(x)\,\dd t=0.
\end{equation}
Thus the geometric construction supplies a single Borel field satisfying
both the continuity equation and the required energy bound.

\begin{proposition}\label{prop:transport-estimate-holomorphic-arc}
Let $\alpha:\Delta_0\to B$ be holomorphic on a disk centered at zero and
put $\mu_t:=\mug(\alpha(t))$. There are $C>0$, $\varepsilon>0$ and
$\theta\in(0,1)$, with $[0,\varepsilon]\subset\Delta_0$, such that
\begin{equation}\label{eq:arc-wasserstein-estimate}
 \W(\mu_s,\mu_t)\leq C\int_s^t r^{-\theta}\,\dd r
 \qquad(0\leq s<t\leq\varepsilon).
\end{equation}
Consequently $\mu\in\mathsf{AC}^1([0,\varepsilon];(\cP(M),\W))$,
$|\dot\mu_t|_\W\leq Ct^{-\theta}$ for almost every $t\in(0,\varepsilon)$, and
\begin{equation}\label{eq:arc-length-estimate}
 L_\W(\mu;[0,\varepsilon])
 \leq C\int_0^\varepsilon r^{-\theta}\,\dd r
 =\frac{C}{1-\theta}\varepsilon^{1-\theta}<\infty.
\end{equation}
\end{proposition}

\begin{proof}
Apply the preceding construction after shrinking the disk.
The curve $t\mapsto\mu_t$ is narrowly continuous by \Cref{prop:mu-gamma}.
The field \eqref{eq:arc-global-field} is Borel, satisfies
\eqref{eq:arc-continuity-equation}, and has
$\|v_t\|_{L^2(\mu_t)}\leq Ct^{-\theta}$ by \eqref{eq:arc-field-energy}.
This bound is square-integrable on every compact subinterval of
$(0,\varepsilon)$ and integrable down to zero. Hence
\Cref{lem:continuity-speed-criterion} gives the metric speed bound and
\eqref{eq:arc-wasserstein-estimate} for $0<s<t<\varepsilon$.

As $s\downarrow0$, continuity of the cycle-to-measure map and compactness
of $M$ give $\mu_s\to\mu_0$ in $\W$. Passing to this limit extends the
distance estimate to $s=0$; continuity also includes $t=\varepsilon$.
The resulting integrable upper bound on the closed interval proves
$\mathsf{AC}^1$ regularity. Finally, for every partition of
$[0,\varepsilon]$, summing \eqref{eq:arc-wasserstein-estimate} over its
successive subintervals gives the same upper bound
$C\int_0^\varepsilon r^{-\theta}\,\dd r$. Taking the supremum over
partitions proves \eqref{eq:arc-length-estimate}.
\end{proof}

\subsubsection{From analytic arcs to finiteness}\label{ss:arc-finiteness}

We now pass from the holomorphic estimate to the curves that occur in the definition of the inner
distance. Complexification treats real-analytic curves; analytic power reparametrization then
treats continuous subanalytic curves.

\begin{theorem}\label{thm:finite-length-curves}
Let \(J=[a,b]\subset\R\). Suppose that \(c:J\to B\) is the restriction of a real-analytic map defined
on an open neighborhood of \(J\). Then
\[
 \mug\circ c\in\mathsf{AC}^1(J;(\cP(M),\W)).
\]
More precisely, for every \(t_0\in J\) and every
\(\sigma\in\{-1,1\}\) for which \(t_0+\sigma(0,\varepsilon)\subset J\)
for some \(\varepsilon>0\), there exist
\(C>0\), \(\varepsilon>0\), and \(\theta\in(0,1)\) such that
\[
 \W\bigl(
   \mug(c(t_0+\sigma s)),
   \mug(c(t_0+\sigma t))
 \bigr)
 \leq
 C\int_s^t r^{-\theta}\,\dd r
\]
for all \(0\leq s<t\leq\varepsilon\).

Consequently, if \(c:J\to B\) is continuous and subanalytic, then
\[
 L_{\W}(\mug\circ c;J)<\infty.
\]
\end{theorem}
\begin{proof}
We first consider the real-analytic case. Fix \(t_0\in J\). By
\Cref{lem:complexification-real-analytic-arc}, after restricting to a
neighborhood of \(t_0\), there is a holomorphic map
\(\alpha:\Delta\to B\) such that
\[
 \alpha(r)=c(t_0+r)
\]
for every sufficiently small real \(r\).

Applying \Cref{prop:transport-estimate-holomorphic-arc} to \(\alpha\), and
to the reflected map \(z\mapsto\alpha(-z)\), shows that on each side of
\(t_0\) there are \(C>0\), \(\varepsilon>0\), and
\(\theta\in(0,1)\) such that
\[
 \W\bigl(
   \mug(c(t_0+\sigma s)),
   \mug(c(t_0+\sigma t))
 \bigr)
 \leq
 C\int_s^t r^{-\theta}\,\dd r,
 \qquad 0\leq s<t\leq\varepsilon,
\]
where \(\sigma=1\) on the right of \(t_0\) and \(\sigma=-1\) on the left.

Since \(r^{-\theta}\in L^1(0,\varepsilon)\), these inequalities provide an
\(L^1\) upper gradient on each side of \(t_0\). If two parameters lie on
opposite sides, we insert the intermediate measure \(\mug(c(t_0))\) and
use the triangle inequality. Thus \(\mug\circ c\) is
\(\mathsf{AC}^1\) on a neighborhood of every \(t_0\in J\).

By compactness, finitely many such neighborhoods cover \(J\). After taking
a subdivision subordinate to this cover and using the finite-gluing
property of absolutely continuous curves, we conclude that
\[
 \mug\circ c\in\mathsf{AC}^1(J;(\cP(M),\W)).
\]

Now let \(c:J\to B\) be continuous and subanalytic. Apply
\Cref{prop:subanalytic-analytic-reparametrization} with respect to the fixed
embedding \(\iota_B:B\hookrightarrow\R^N\). This gives a finite subdivision
of \(J\) into compact intervals \(I_\ell\), compact intervals \(K_\ell\),
and monotone homeomorphisms
\[
 \rho_\ell:K_\ell\longrightarrow I_\ell
\]
such that \(c\circ\rho_\ell\) extends real analytically to a neighborhood of
\(K_\ell\). By the first part,
\[
 L_\W(\mug\circ c\circ\rho_\ell;K_\ell)<\infty.
\]
Length is invariant under monotone homeomorphic reparametrizations, and
therefore
\[
 L_\W(\mug\circ c;I_\ell)
 =
 L_\W(\mug\circ c\circ\rho_\ell;K_\ell)<\infty.
\]
Summing over the finite subdivision gives
\[
 L_\W(\mug\circ c;J)
 =
 \sum_\ell L_\W(\mug\circ c;I_\ell)<\infty,
\]
as claimed.
\end{proof}

\begin{proof}[Proof of \Cref{thm:finiteness}]
By \Cref{cor:piecewise-analytic-connectivity}, any two points of $A$ can be joined by
finitely many real-analytic arcs contained in $A$. Each has finite Wasserstein length by
\Cref{thm:finite-length-curves}. Their concatenation is therefore an admissible curve of
finite length in the definition of $\W_A^{\mathrm{in}}$.
\end{proof}

The finiteness theorem also turns the general compactness result for inner metrics into a concrete
existence statement for minimizing deformations of cycles.

\begin{theorem}\label{thm:inner-wasserstein-complete-geodesic}
Let $A\subset B$ be a connected analytic subset. Then
$(A,\W_A^{\mathrm{in}})$ is a complete geodesic metric space. In particular, this holds for $A=B$.
\end{theorem}

\begin{proof}
Let $\W_A$ denote the pullback of $\W$ by $\mug$. By \Cref{prop:mu-gamma}, it is a metric whose
topology is the Barlet topology. The analytic subset $A$ is compact, and
\Cref{thm:finiteness} shows that the inner distance associated with $\W_A$ is finite. The result
therefore follows from \Cref{propo:folk}.
\end{proof}

\subsection{Comparison between the Weil--Petersson and the Wasserstein geometry}
\label{subsubsec:differential-mug}

The preceding proof used a Borel normal velocity field to bound the Wasserstein
speed along holomorphic arcs, including near singular cycles.
Fix now an irreducible component $C\subset B$ with $\widetilde C\neq\varnothing$.
On $\widetilde C$, a smooth curve has a normal deformation field that is a Wasserstein
tangent vector and solves the continuity equation with minimal norm. This identifies the
Wasserstein speed exactly with the Weil--Petersson norm. We then compare the resulting
length distances and their completions.

We denote by $\lVert\cdot\rVert_{\mathrm{WP}}$ the norm induced by the normalized
Weil--Petersson metric $\Theta$, so $\lVert v\rVert_{\mathrm{WP}}^2=\Theta(v,v)$.
\begin{theorem}
\label{thm:isometry}
Let $I$ be a compact interval and let $t\mapsto\gamma_t$ be a smooth curve in
$\widetilde{C}$. Write
\[
 \mu_t:=\mug(\gamma_t),
 \qquad
 \nu_t:=\nu_{\dot\gamma_t}.
\]
Then $t\mapsto\mu_t$ belongs to
$\mathsf{AC}^2(I;(\cP(M),\W))$. Its metric derivative exists at every $t\in I$, with the
one-sided interpretation at the endpoints, and
\begin{equation}\label{eq:metriccompare}
 |\dot\mu_t|_\W^2
 =\int_M|\nu_t|_g^2\,\dd\mu_t
 =\|\dot\gamma_t\|_{\mathrm{WP}}^2.
\end{equation}
Consequently, every piecewise smooth curve in $\widetilde{C}$ has the same Wasserstein and
Weil--Petersson lengths.
\end{theorem}

\begin{proof}
Put $Z_t:=|\gamma_t|$.

By the smooth-family description in \Cref{Ss:tangentgood}, locally on $I$ the pullback of the
tautological family is a smooth proper family of embedded complex submanifolds.

After choosing a smooth trivialization, we may write it as
\[
 F_t:Z\longrightarrow Z_t\subset M.
\]
Let $W_t$ be the velocity field along $Z_t$ determined by
\[
 W_t\circ F_t=\partial_tF_t.
\]
By \Cref{lemma:normalnu}, its orthogonal normal component is $\nu_t$.

\smallskip
\noindent\emph{Step 1: the continuity equation.}
For every $\varphi\in C^\infty(M)$, the definition of the cycle-to-measure map gives
\[
 \int_M\varphi\,\dd\mu_t
 =\frac1{V_B}\int_{Z_t}\varphi\frac{\omega^q}{q!}.
\]
Cartan's formula, $\dd(\omega^q/q!)=0$, and Stokes' theorem imply
\begin{equation}\label{eq:good-family-first-variation}
 \frac{\dd}{\dd t}\int_M\varphi\,\dd\mu_t
 =\frac1{V_B}\int_{Z_t}\iota_{W_t}(\dd\varphi\wedge\frac{\omega^q}{q!}).
\end{equation}
At $x\in Z_t$, apply \Cref{lem:kahler-normal-contraction} with
$P=T_xZ_t$ and $\zeta=W_t(x)$. Since the normal component of $W_t$ is $\nu_t$, we obtain
\[
 \left.\iota_{W_t}(\dd\varphi\wedge\frac{\omega^q}{q!})\right|_{T_xZ_t}
 =\langle\nabla\varphi,\nu_t\rangle_g\left.\frac{\omega^q}{q!}\right|_{T_xZ_t}.
\]
Substitution into \eqref{eq:good-family-first-variation} gives
\begin{equation}\label{eq:good-family-continuity-equation}
 \frac{\dd}{\dd t}\int_M\varphi\,\dd\mu_t
 =\int_M\langle\nabla\varphi,\nu_t\rangle_g\,\dd\mu_t.
\end{equation}
Thus $(\mu_t,\nu_t)$ solves the continuity equation. The smoothness of the family implies that
$\|\nu_t\|_{L^2(\mu_t)}$ is bounded on compact parameter intervals. The
continuity-equation characterization of Wasserstein curves therefore gives
$\mu\in\mathsf{AC}^2$ and
\begin{equation}\label{eq:good-family-speed-upper-bound}
 |\dot\mu_t|_\W\leq\|\nu_t\|_{L^2(\mu_t)}
\end{equation}
for almost every $t$.

\smallskip
\noindent\emph{Step 2: minimality of the normal velocity.}

Because $\nu_t$ is normal to $Z_t$, the tangent-space decomposition
\eqref{eq:tangentsplitting} gives
\[
 \nu_t\in T_{\mu_t}\cP(M)
\]
for every $t$. Since $\nu_t$ satisfies the continuity equation
\eqref{eq:good-family-continuity-equation}, the characterization of tangent velocities in
\cite[Proposition~2.5]{Erbar} shows that it is the distinguished velocity of minimal
$L^2(\mu_t)$-norm. Therefore
\begin{equation}\label{eq:good-family-speed-ae}
 |\dot\mu_t|_\W=\|\nu_t\|_{L^2(\mu_t)}
\end{equation}
for almost every $t$.

\smallskip

\noindent\emph{Step 3: the identity at every time.}
We give a direct argument that also treats the endpoints of $I$.
Set $a(t):=\|\nu_t\|_{L^2(\mu_t)}$. Smooth dependence of the family
shows that $a$ is continuous. Integrating
\eqref{eq:good-family-speed-upper-bound}, we obtain
\[
 \W(\mu_s,\mu_t)\leq\int_s^t a(r)\,\dd r
 \qquad(s<t).
\]
In particular, at every fixed $t\in I$,
\[
 \limsup_{h\to0}\frac{\W(\mu_{t+h},\mu_t)}{|h|}\leq a(t),
\]
where $t+h\in I$.

For the reverse inequality, choose $f\in C^\infty(M)$ with
$\nabla f|_{Z_t}=\nu_t$. Such a function exists because $\nu_t$ is
normal: in a tubular neighborhood, set
$f(\exp_x n)=\langle\nu_t(x),n\rangle_g$ and multiply by a cutoff
equal to one near the zero section. Taylor's formula along a minimizing
geodesic gives
\[
 |f(y)-f(x)|
 \leq |\nabla f(x)|_g\,\dM(x,y)+K_f \dM(x,y)^2,
\]
where $K_f$ depends only on a bound for the Hessian of $f$.
Integrating against an optimal coupling of $\mu_t$ and $\mu_{t+h}$
and using Cauchy--Schwarz yields
\[
 \left|\int_M f\,\dd(\mu_{t+h}-\mu_t)\right|
 \leq a(t)\W(\mu_{t+h},\mu_t)
       +K_f\W(\mu_{t+h},\mu_t)^2.
\]
By \eqref{eq:good-family-continuity-equation}, the expression inside
the absolute value is $h\,a(t)^2+o(|h|)$. The upper bound already proved
gives $\W(\mu_{t+h},\mu_t)=O(|h|)$. If $a(t)>0$, division by $|h|$
therefore gives
\[
 a(t)\leq
 \liminf_{h\to0}\frac{\W(\mu_{t+h},\mu_t)}{|h|}.
\]
If $a(t)=0$, the upper bound alone gives the desired conclusion.
Thus the metric derivative exists at every time and equals $a(t)$.

\smallskip
\noindent\emph{Step 4: identification with the Weil--Petersson norm.}
By the definition of $\mu_t$, the Hermitian identity
$|\nu_t|_g^2=\omega(\nu_t,J\nu_t)$, and using \cref{def:WP}, we get
\begin{align*}
 \int_M|\nu_t|_g^2\,\dd\mu_t
 &=\frac1{V_B}\int_{Z_t}|\nu_t|_g^2\frac{\omega^q}{q!}\\
 &=\frac1{V_B}\int_{Z_t}\omega(\nu_t,J\nu_t)\frac{\omega^q}{q!}\\
 &=\owp_{\gamma_t}(\dot\gamma_t,J\dot\gamma_t)
 =\|\dot\gamma_t\|_{\mathrm{WP}}^2.
\end{align*}
Combining this identity with Step~3 proves \eqref{eq:metriccompare}. Integrating the speed
identity on each smooth piece proves the equality of lengths.
\end{proof}
\begin{remark}

\label{rem:infinitesimal-wasserstein-kahler}
\Cref{thm:isometry} shows that, on the good locus, the cycle-to-measure
map carries more than the topology of the cycle space into Wasserstein space.
It also transports its infinitesimal K\"ahler geometry. Indeed, if
$\gamma\in\widetilde{C}$, $Z:=|\gamma|$, and $v\in T_\gamma\widetilde{C}$, then
the normal deformation field $\nu_v$ may be regarded, through the splitting
\eqref{eq:tangentsplitting}, as a Wasserstein tangent vector at
$\mug(\gamma)$. This gives a natural linear differential
\[
D_\gamma\mug:T_\gamma\widetilde{C}\longrightarrow T_{\mug(\gamma)}\cP(M),
\qquad
D_\gamma\mug(v):=\nu_v .
\]
The tangent-space decomposition \eqref{eq:tangentsplitting}, as used in the proof of
\Cref{thm:isometry}, shows directly that
$D_\gamma\mug(v)$ belongs to the Wasserstein tangent space.  Its norm is
\[
 \bigl\|D_\gamma\mug(v)\bigr\|_{L^2(\mug(\gamma))}
 =\|v\|_{\mathrm{WP}}.
\]
Along a smooth curve this vector is the minimal continuity-equation velocity. By
\Cref{thm:isometry}, its norm gives the metric speed at every time.

The symplectic form has a direct infinitesimal description as well. Every smooth normal
field $\nu_v$ along the complex submanifold $Z$ is the restriction of a Hamiltonian
vector field on $M$. Indeed, the tubular-neighborhood construction gives a smooth
function $f_v$, vanishing on $Z$, with $\nabla f_v|_Z=-J\nu_v$.
With the convention $\iota_{X_f}\omega=-\dd f$, one has $X_f=J\nabla f$, so
$X_{f_v}|_Z=\nu_v$. Consequently, on the image of $D_\gamma\mug$,
\[
 \int_M\omega\bigl(D_\gamma\mug(v),D_\gamma\mug(w)\bigr)\,\dd\mug(\gamma)
 =\frac1{V_B}\int_Z\omega(\nu_v,\nu_w)\,\frac{\omega^q}{q!}
 =\owp_\gamma(v,w).
\]
Thus both the metric and the symplectic form of the normalized Weil--Petersson structure
are recovered by integrating the corresponding ambient tensors on these normal
deformation fields.

\end{remark}

\subsubsection{Bi-Lipschitz properties}
The proof of the next theorem is a Wasserstein analogue of the Barlet--Magnusson
principle that cycles can be locally separated by finitely many integration
functionals. More precisely, in \cite[Theorem 11.3.2.1]{cycles2} the authors construct local embeddings of the
cycle space by taking finitely many functions defined by integration of local holomorphic forms (see also their interpretation in terms of the
Andreotti--Norguet map \cite[Section 11.4]{cycles2}). 
Here, instead of using
cohomological integration coordinates, we choose smooth test functions
\(f_1,\ldots,f_m\in C^\infty(M)\) whose normal gradients are the images, under \eqref{eq:tangent-normal}, of a real basis of $T_\gamma C$. As we will see, the  map
\[
F(\eta)=
\left(
\int_M f_1\,d\mug(\eta),
\ldots,
\int_M f_m\,d\mug(\eta)
\right)
\]
is a local coordinate system near the good cycle. Kantorovich duality
then turns these coordinates into a lower bound for the ambient Wasserstein
distance.

\begin{theorem}\label{thm:local-bilipschitz}
Let $C\subset B$ be an irreducible component with $\widetilde C\neq\varnothing$. For every \(\gamma\in \widetilde{C}\)  there exist an open
neighborhood \(U\subset \widetilde{C}\) of \(\gamma\) and a constant
\(C_U\ge 1\) such that, for every \(\eta,\theta\in U\),
\begin{equation}\label{eq:bilipschitz}
C_U^{-1} \mathsf{d}_{\mathrm{WP}}(\eta,\theta)
\le
\W\bigl(\mug(\eta),\mug(\theta)\bigr)
\le
\mathsf{d}_{\mathrm{WP}}(\eta,\theta).
\end{equation}

\end{theorem}

\begin{proof}
If $\dim_\C C=0$, the assertion is immediate. Assume $\dim_\C C>0$.
Let \(Z:=|\gamma|\), and put
\[
m:=\dim_{\mathbb R}T_\gamma\widetilde{C}=2\dim_\C T_\gamma\widetilde{C} .
\]

Pick a real basis $\{e_1,\ldots,e_m\}$ of $T_\gamma C$ and put
$\sigma_i:=\nu_{e_i}\in\Gamma(Z,TZ^\perp)$. These normal fields are linearly independent
by the injectivity of \eqref{eq:tangent-normal}.

Using the tubular neighborhood theorem (as in the proof of \cite[Proposition~2.8]{ACL}), for every $i=1, \ldots, m$ there exists a smooth function $f_i:M\to \R$ such that
\[
(\nabla f_i)|_Z=\sigma_i.
\]
We will use the functions $f_1, \ldots, f_m$ to build a smooth embedding of a neighborhood of $\gamma$ into $\R^m$, as follows.
Define
\[
F:\widetilde{C}\longrightarrow \mathbb R^m,
\qquad
F(\eta)
:=
\left(
\int_M f_1\,d\mug(\eta),
\ldots,
\int_M f_m\,d\mug(\eta)
\right).
\]
Equivalently, if \(\eta\) is represented by the smooth reduced cycle
\(Z_\eta\), then
\[
F_i(\eta)
=
\frac{1}{V_B}
\int_{Z_\eta} f_i\,\frac{\omega^q}{q!}.
\]

We show that \(F\) is a local coordinate system at \(\gamma\). The tautological family is a
smooth proper submersion over a neighborhood of \(\gamma\), so differentiation under fiber
integration shows directly that $F$ is smooth there.

We show now that the differential of $F$ at $\gamma$ is an isomorphism. Let
\(v\in T_\gamma\widetilde{C}\), and denote by \(\nu_v\) the corresponding normal
deformation field along \(Z\). Applying the first-variation identity
\eqref{eq:good-family-continuity-equation} to a curve tangent to $v$ gives
\[
(D_\gamma F_i)(v)
=
\frac{1}{V_B}
\int_Z
\left\langle \nabla f_i,\nu_v\right\rangle
\frac{\omega^q}{q!}=\frac{1}{V_B}
\int_Z
\left\langle \sigma_i,\nu_v\right\rangle
\frac{\omega^q}{q!}.
\]
In particular,\[
(D_\gamma F_i)(e_j)
=
\frac{1}{V_B}
\int_Z
\langle \sigma_i,\sigma_j\rangle
\frac{\omega^q}{q!}.
\]
Thus the matrix of \(D_\gamma F\) is the \(L^2\)-Gram matrix of the linearly independent normal fields
\(\sigma_1,\ldots,\sigma_m\). This matrix is positive definite, hence \(D_\gamma F\) is an
isomorphism. Therefore, by the inverse function theorem, after shrinking around \(\gamma\), there is an
open neighborhood \(U\subset \widetilde{C}\) such that \(F|_U\) is a
diffeomorphism onto its image. 

Choose a Euclidean ball \(Q\), centered at \(F(\gamma)\), whose closure is contained in the
coordinate image, and replace \(U\) by \(F^{-1}(Q)\). The differential of \(F^{-1}\) is bounded
on \(\overline Q\). Lifting the straight segment from \(F(\eta)\) to \(F(\theta)\) therefore gives
a constant \(c_U>0\) such that, for all \(\eta,\theta\in U\),
\begin{equation}\label{eq:WPF}
c_{U}\, \mathsf{d}_{\mathrm{WP}}(\eta,\theta)
\le
\|F(\eta)-F(\theta)\|_{\R^m},
\end{equation}
where $\|\cdot\|_{\R^m}$ denotes the standard Euclidean norm on $\R^m$. 

We now use \eqref{eq:WPF} to compare the Weil--Petersson distance with the Wasserstein distance on $U$. Let
\[
L_F
:=
\left(
\sum_{i=1}^m \operatorname{Lip}(f_i)^2
\right)^{1/2}.
\]
By Kantorovich duality for  the Wasserstein $1$-distance \(W_1\), for each \(i=1, \ldots, m\) we have
\[
\left|
F_i(\eta)-F_i(\theta)
\right|
=
\left|
\int_M f_i\,d\mug(\eta)
-
\int_M f_i\,d\mug(\theta)
\right|
\le
\operatorname{Lip}(f_i)\,
W_1\bigl(\mug(\eta),\mug(\theta)\bigr).
\]
Consequently,
\[
\|F(\eta)-F(\theta)\|_{\R^m}
\le
L_F\,W_1\bigl(\mug(\eta),\mug(\theta)\bigr).
\]
Since \(W_1\le \W\) on \(\cP(M)\), we get
\[
\|F(\eta)-F(\theta)\|_{\R^m}
\le
L_F\,\W\bigl(\mug(\eta),\mug(\theta)\bigr).
\]
Combining this with the coordinate estimate \eqref{eq:WPF} gives
\[
\W\bigl(\mug(\eta),\mug(\theta)\bigr)
\ge
\frac{c_{U}}{L_F}\,
\mathsf{d}_{\mathrm{WP}}(\eta,\theta)
\]
for all $\eta,\theta\in U$.
Setting $
C_U
:=
\max\left\{
1,\frac{L_F}{c_{U}}
\right\}, $ this gives the left inequality in \eqref{eq:bilipschitz}.

It remains to prove the upper bound for $\W$ on $U$ in terms of $\mathsf{d}_{\mathrm{WP}}$. Let
\(\alpha:[0,1]\to \widetilde{C}\) be a piecewise smooth curve joining
\(\eta\) to \(\theta\). By \Cref{thm:isometry}, the curve
\(\mug\circ\alpha\) is absolutely continuous in
\((\cP(M),\W)\). Writing \(\mu_t:=\mug(\alpha(t))\), its metric speed satisfies
\[
|\dot\mu_t|_{\W}
=
\|\dot\alpha_t\|_{\mathrm{WP}}
\]
for almost every \(t\in [0,1]\). Therefore
\[
\W\bigl(\mug(\eta),\mug(\theta)\bigr)
\le
L_{\W}(\mug\circ\alpha)
=
L_{\mathrm{WP}}(\alpha).
\]
Since the Riemannian distance can be obtained by taking infimum on smooth curves only, taking the infimum over all such curves \(\alpha\) in the above inequality gives
\[
\W\bigl(\mug(\eta),\mug(\theta)\bigr)
\le
\mathsf{d}_{\mathrm{WP}}(\eta,\theta).\qedhere
\]
\end{proof}

\begin{corollary}
\label{cor:mu-preserves-ac2}
The map
\(\mug|_{\widetilde{C}}:\widetilde{C}\longrightarrow \cP(M)\)
and its inverse
\((\mug|_{\widetilde{C}})^{-1}:\mug(\widetilde{C})\longrightarrow \widetilde{C}\)
send \(\mathsf{AC}^2\)-curves to \(\mathsf{AC}^2\)-curves.
\end{corollary}
\begin{proof}
We prove only the implication from \(\mug(\widetilde{C})\) back to
\(\widetilde{C}\); the other one is identical, using the opposite inequality in
\cref{thm:local-bilipschitz}.

Let
\(\beta\in \mathsf{AC}^2([0,1];(\mug(\widetilde{C}),\W))\), and set
\(\alpha:=(\mug|_{\widetilde{C}})^{-1}\circ\beta\). Since
\(\mug|_{\widetilde{C}}\) is a homeomorphism onto its image by
\cref{prop:mu-gamma}, the curve \(\alpha\) is continuous, and therefore
\(K:=\alpha([0,1])\) is compact. By \cref{thm:local-bilipschitz}, we can cover
\(K\) by finitely many open sets \(U_1,\ldots,U_N\subset\widetilde{C}\) such that
on each \(U_i\) one has
\(\mathsf{d}_{\mathrm{WP}}\le C_i \W(\mug(\cdot),\mug(\cdot))\). Let
\(C_K:=\max_i C_i\).

Let \(g\in L^2(0,1)\) be an admissible control for \(\beta\) in the sense of
\eqref{eq:AC}. Fix \(s<t\). Since \(\alpha\) is continuous, we may choose a
partition \(s=t_0<\cdots<t_L=t\) such that each subcurve
\(\alpha([t_{j-1},t_j])\) is contained in one of the \(U_i\)'s. Hence
\[
\mathsf{d}_{\mathrm{WP}}(\alpha(s),\alpha(t))
\le
\sum_j \mathsf{d}_{\mathrm{WP}}(\alpha(t_{j-1}),\alpha(t_j))
\le
C_K\sum_j \W(\beta(t_{j-1}),\beta(t_j))
\le
C_K\int_s^t g(r)\,dr .
\]
Since \(C_Kg\in L^2(0,1)\), the definition \eqref{eq:AC} gives
\(\alpha\in \mathsf{AC}^2([0,1];(\widetilde{C},\mathsf{d}_{\mathrm{WP}}))\).
\end{proof}

\subsubsection{Inner distances and metric completion}
The domain in which paths are allowed to move matters in this comparison.
We first identify the inner Wasserstein distance computed within
$\widetilde C$ with the Weil--Petersson distance. We then allow paths to
pass through the singular or nonreduced cycles of $C$ and compare the
resulting complete metric with the completion of the good locus.

\begin{proposition}\label{prop:local-wp-inner-wasserstein}
Let \(C\subset B\) be an irreducible component with
\(\widetilde C\neq\varnothing\). Then
\[
 (\widetilde C,\W_{\widetilde C}^{\mathrm{in}})
 =
 (\widetilde C,\mathsf d_{\mathrm{WP}})
\]
as metric spaces. Moreover,
\[
 \W_C^{\mathrm{in}}\leq\mathsf d_{\mathrm{WP}}
 \quad\text{on }\widetilde C\times\widetilde C,
\]
and the two distances agree on a neighborhood of the diagonal in
\(\widetilde C\times\widetilde C\).
\end{proposition}

\begin{proof}
If \(t\mapsto\gamma_t\) is a smooth curve in \(\widetilde C\), then
\Cref{thm:isometry} shows that \(t\mapsto\mug(\gamma_t)\) belongs to
\(\mathsf{AC}^2\) and
\[
 L_{\mathrm{WP}}(\gamma)
 =
 L_{\W}(\mug\circ\gamma).
\]
Consequently,
\[
 \mathsf d_{\mathrm{WP}}
 =
 \inf_{\gamma\in C^\infty}L_{\mathrm{WP}}(\gamma)
 =
 \inf_{\gamma\in C^\infty}L_{\W}(\mug\circ\gamma)
 \geq
 \inf_{\gamma\in\mathsf{AC}^2}L_{\W}(\mug\circ\gamma)
 =
 \W_{\widetilde C}^{\mathrm{in}},
\]
where all infima are taken over curves with the prescribed endpoints.

To prove the reverse inequality, we extend the speed identity of
\Cref{thm:isometry} to \(\mathsf{AC}^2\)-curves. By
\Cref{cor:mu-preserves-ac2}, the maps
\[
 \mug|_{\widetilde C}
 \quad\text{and}\quad
 (\mug|_{\widetilde C})^{-1}
\]
preserve the class of \(\mathsf{AC}^2\)-curves. Let
\(\gamma\in\mathsf{AC}^2(J;\widetilde C)\), set
\(\mu_t:=\mug(\gamma_t)\), and fix a time \(t\) at which \(\gamma\) is
differentiable. Write
\[
 x:=\gamma_t,
 \qquad
 v:=\dot\gamma_t.
\]
Differentiability in the Riemannian manifold
\((\widetilde C,\mathsf d_{\mathrm{WP}})\) gives
\[
 \mathsf d_{\mathrm{WP}}
 \bigl(\gamma_{t+h},\exp_x(hv)\bigr)
 =o(|h|).
\]
Since \(\mug\) is locally Lipschitz by
\Cref{thm:local-bilipschitz}, it follows that
\[
 \W\bigl(\mu_{t+h},\mug(\exp_x(hv))\bigr)
 =o(|h|).
\]
The reverse triangle inequality therefore yields
\[
 \left|
 \W(\mu_{t+h},\mu_t)
 -
 \W\bigl(\mug(\exp_x(hv)),\mu_t\bigr)
 \right|
 \leq
 \W\bigl(\mu_{t+h},\mug(\exp_x(hv))\bigr).
\]
Dividing by \(|h|\), letting \(h\to0\), and applying
\Cref{thm:isometry} to the smooth curve
\(h\mapsto\exp_x(hv)\), we obtain
\[
 |\dot\mu_t|_{\W}
 =
 \lim_{h\to0}
 \frac{\W(\mu_{t+h},\mu_t)}{|h|}
 =
 \|v\|_{\mathrm{WP}}
\]
for almost every \(t\). Thus every
\(\mathsf{AC}^2\)-curve in \(\widetilde C\) has the same Wasserstein and
Weil--Petersson lengths. Hence
\[
 \mathsf d_{\mathrm{WP}}
 =
 \inf_{\gamma\in C^\infty}L_{\mathrm{WP}}(\gamma)
 =
 \inf_{\gamma\in\mathsf{AC}^2}L_{\mathrm{WP}}(\gamma)
 =
 \inf_{\gamma\in\mathsf{AC}^2}L_{\W}(\mug\circ\gamma)
 =
 \W_{\widetilde C}^{\mathrm{in}},
\]
where the second equality follows from
\cite[Corollary~3.13]{burtscher}.

Allowing curves to travel in all of \(C\) enlarges the class of
competitors. Therefore
\[
 \W_C^{\mathrm{in}}
 \leq
 \W_{\widetilde C}^{\mathrm{in}}
 =
 \mathsf d_{\mathrm{WP}}
 \quad\text{on }\widetilde C\times\widetilde C.
\]

It remains to prove the local equality. Fix
\(\gamma\in\widetilde C\), and choose open neighborhoods \(W\) and \(V\)
of \(\gamma\) such that
\[
 \overline W\subset V,
 \qquad
 \overline V\subset\widetilde C.
\]

Since $\mug$ is a homeomorphism onto its image and $C$ is compact,
there is $\delta>0$ such that
\[
 \W(\mug(\eta),\mug(\theta))\geq\delta
 \qquad(\eta\in\overline W,\ \theta\in C\setminus V).
\]
If $C\setminus V$ is empty, any $\delta>0$ has this property.

Choose an open neighborhood \(U_\gamma\) of \(\gamma\), with
\(\overline{U_\gamma}\subset W\), whose
\(\mathsf d_{\mathrm{WP}}\)-diameter is smaller than \(\delta\).

Let \(\eta,\theta\in U_\gamma\). Since
\[
 \W_C^{\mathrm{in}}(\eta,\theta)
 \leq
 \mathsf d_{\mathrm{WP}}(\eta,\theta)
 <\delta,
\]
competitors whose lengths converge to
\(\W_C^{\mathrm{in}}(\eta,\theta)\) may be chosen with length strictly
smaller than \(\delta\). Such a curve cannot meet
\(\mug(C\setminus V)\), because it starts in \(\mug(\overline W)\).
It therefore remains in \(\mug(V)\). After arclength reparametrization on a compact interval, it is Lipschitz in Wasserstein distance. By \Cref{cor:mu-preserves-ac2}, its pullback by
\((\mug|_{\widetilde C})^{-1}\) is an
\(\mathsf{AC}^2\)-curve in \(\widetilde C\), and the speed identity gives
\[
 L_{\W}(\beta)
 =
 L_{\mathrm{WP}}
 \bigl((\mug|_{\widetilde C})^{-1}\circ\beta\bigr)
 \geq
 \mathsf d_{\mathrm{WP}}(\eta,\theta).
\]
Taking the infimum over these competitors proves
\[
 \W_C^{\mathrm{in}}(\eta,\theta)
 \geq
 \mathsf d_{\mathrm{WP}}(\eta,\theta).
\]
Together with the opposite inequality already established, this gives
equality on \(U_\gamma\times U_\gamma\), and hence on a neighborhood of
the diagonal.
\end{proof}

\begin{corollary}\label{thm:wp-completion-comparison}
Let \(C\subset B\) be an irreducible component with
\(\widetilde C\neq\varnothing\), and let
\(\widehat C_{\mathrm{WP}}\) denote the metric completion of
\((\widetilde C,\mathsf d_{\mathrm{WP}})\). Then the identity on
\(\widetilde C\) extends uniquely to a canonical \(1\)-Lipschitz
surjection
\[
 \overline\iota:
 \widehat C_{\mathrm{WP}}
 \longrightarrow
 (C,\W_C^{\mathrm{in}}).
\]
\end{corollary}

\begin{proof}
By \Cref{prop:local-wp-inner-wasserstein},
\(\W_C^{\mathrm{in}}\leq\mathsf d_{\mathrm{WP}}\) on
\(\widetilde C\times\widetilde C\). Hence the inclusion
\[
 (\widetilde C,\mathsf d_{\mathrm{WP}})
 \longrightarrow
 (C,\W_C^{\mathrm{in}})
\]
is \(1\)-Lipschitz. Since \(C\) is irreducible, it is connected, and
\Cref{thm:inner-wasserstein-complete-geodesic} shows that
\((C,\W_C^{\mathrm{in}})\) is complete. The inclusion therefore extends
uniquely to a \(1\)-Lipschitz map
\(\overline\iota:\widehat C_{\mathrm{WP}}\to C\).

It remains to prove that this extension is surjective. Let
\(\gamma\in C\); there is nothing to prove if
\(\gamma\in\widetilde C\). Otherwise,
\Cref{prop:bad-locus-analytic} shows that \(\widetilde C\) is a dense
definable open subset of \(C\). By the definable curve-selection lemma
\cite[Curve Selection Lemma]{valette}, there are \(\varepsilon>0\) and a
continuous definable curve \(c:[0,\varepsilon]\to C\) such that
\(c(0)=\gamma\) and \(c(t)\in\widetilde C\) for \(t>0\).

By \Cref{thm:finite-length-curves}, the curve \(\mug\circ c\) has finite
Wasserstein length. On every interval \([s,t]\subset(0,\varepsilon]\),
\Cref{prop:local-wp-inner-wasserstein} and the invariance of the length of
a rectifiable curve under passage to the induced inner metric give
\[
 \mathsf d_{\mathrm{WP}}\bigl(c(s),c(t)\bigr)
 \leq
 L_{\W}(\mug\circ c;[s,t]).
\]
Since \(\mug\circ c\) has finite length, the right-hand side tends to zero
as \(s,t\downarrow0\). Thus \(c(t)\) is Cauchy for
\(\mathsf d_{\mathrm{WP}}\) as \(t\downarrow0\), and hence determines a
point \(p\in\widehat C_{\mathrm{WP}}\).

Finally, the restriction of \(c\) to \([0,t]\) is an admissible curve in
\(C\), so
\[
 \W_C^{\mathrm{in}}\bigl(c(t),\gamma\bigr)
 \leq L_{\W}(\mug\circ c;[0,t])
 \longrightarrow0
 \qquad (t\downarrow0).
\]
By the definition of the extension, this means that
\(\overline\iota(p)=\gamma\). Therefore \(\overline\iota\) is
surjective.
\end{proof}

\begin{remark}\label{rem:completion-domains}
The completion in \Cref{thm:wp-completion-comparison} is also the
completion of $(\widetilde C,\W_{\widetilde C}^{\mathrm{in}})$,
because the two metrics on $\widetilde C$ are identical. Thus completing
the inner Wasserstein metric on the good locus does not produce a
different space from completing the Weil--Petersson metric there.
The comparison is with $(C,\W_C^{\mathrm{in}})$, where paths may use
all cycles from the outset. The corollary provides a canonical
surjection onto this space; it does not assert that the surjection is
injective or isometric.

For a connected component $B$ containing several irreducible components,
there is a further distinction. Its good locus may be disconnected,
whereas every two points of $B$ have finite inner Wasserstein distance.
Completing each
connected component of the good locus separately does not identify
boundary points belonging to different components. The following
example makes this distinction explicit.
\end{remark}

\begin{example}[Conics and pairs of lines in $\CP^3$]
\label{ex:conics-p3-completion}
Equip $M=\CP^3$ with the Fubini--Study metric and let
\[
 B=\operatorname{Chow}_{1,2}(\CP^3)
\]
be the space of effective one-dimensional cycles of degree two.
It has two irreducible components: $C_{\mathrm{con}}$, the closure of
the locus of smooth plane conics, and $C_{\mathrm{lines}}$, the locus
of unordered pairs of lines, with repetition allowed
\cite[Section~2]{quadratic-space-curves}. Their intersection consists
of pairs of intersecting lines and double lines. In particular, $B$
is connected. The good loci are
\begin{align*}
 U_{\mathrm{con}}=\widetilde C_{\mathrm{con}}
   &=\{\text{smooth plane conics}\},\\
 U_{\mathrm{lines}}=\widetilde C_{\mathrm{lines}}
   &=\{\text{pairs of skew lines}\}.
\end{align*}
Each is a single orbit of the connected group
$\operatorname{PGL}_4(\C)$, hence a connected smooth parameter space, and
$\widetilde B=U_{\mathrm{con}}\sqcup U_{\mathrm{lines}}$.
A pair of skew lines is a smooth reduced cycle with disconnected
support, as allowed in our definition of the good locus.

Use homogeneous coordinates $[x_0:x_1:x_2:x_3]$ and set
\[
 L_1=\{x_0=x_3=0\},\qquad
 L_2=\{x_1=x_3=0\},\qquad
 \gamma_0=[L_1]+[L_2].
\]
Consider the holomorphic arcs
\begin{align*}
 \alpha(t)&=\bigl[\{x_3=0,\ x_0x_1=t x_2^2\}\bigr],\\
 \beta(t)&=[L_1]+\bigl[\{x_1=0,\ x_3=t x_2\}\bigr].
\end{align*}
For $t\neq0$, the first cycle is a smooth conic and the two lines
in the second cycle are skew. At $t=0$, both arcs reach $\gamma_0$.
By \Cref{thm:finite-length-curves}, both curves of measures have
finite Wasserstein length on $[0,\varepsilon]$. Hence
\begin{align*}
 \W_B^{\mathrm{in}}(\alpha(t),\beta(t))
 &\leq L_\W(\mug\circ\alpha;[0,t])
       +L_\W(\mug\circ\beta;[0,t])\\
 &\longrightarrow0\qquad(t\downarrow0).
\end{align*}
The admissible connecting curve follows $\alpha$ backwards to
$\gamma_0$ and then follows $\beta$ forwards.

On every compact subinterval of $(0,\varepsilon]$, the two arcs lie
in their respective good loci. By \Cref{thm:isometry}, their
Weil--Petersson lengths equal their Wasserstein lengths. The vanishing
of the lengths of the tails therefore shows that $\alpha(t)$ and
$\beta(t)$ define points
\[
 p_{\mathrm{con}}\in\widehat C_{\mathrm{con},\mathrm{WP}},
 \qquad
 p_{\mathrm{lines}}\in\widehat C_{\mathrm{lines},\mathrm{WP}}.
\]
The two completion maps, followed by the inclusions of the irreducible
components into $B$, send both points to $\gamma_0$. Consequently,
the canonical map
\[
 \widehat C_{\mathrm{con},\mathrm{WP}}
 \sqcup\widehat C_{\mathrm{lines},\mathrm{WP}}
 \longrightarrow (B,\W_B^{\mathrm{in}})
\]
is not injective. There is no path in $\widetilde B$ joining a smooth
conic to a pair of skew lines; their Weil--Petersson length distance,
if regarded as an extended distance on $\widetilde B$, is infinite.
The inner Wasserstein geometry on $B$ instead joins the two components
through the common singular cycle.

This example concerns the connected, reducible space $B$. It does not
assert a failure of injectivity of the completion map for either
individual irreducible component.
\end{example}

\section{Examples}
\label{subsec:examples}

\subsection{Optimal transport between fibers}\label{ex:transport-between-fibers}
Let $f:M\to Y$ be a proper surjective open holomorphic map onto a connected complex
manifold, and put $q=\dim_\C M-\dim_\C Y$. This gives a direct application of the
cycle-space construction: the ambient manifold is the total space of the family.

The openness assumption is equivalent here to flatness and to every fiber being pure
$q$-dimensional \cite[p.~158, Corollary]{Fischer1976}. With their natural fiber
multiplicities, the fibers therefore form an analytic family of compact cycles
\cite[Proposition~4.3.3]{cycles1}, classified by
\[
 \gamma_f:Y\longrightarrow B\subset\mathcal C_q(M),\qquad
 y\longmapsto[M_y].
\]
Different fibers have disjoint nonempty supports, so $\gamma_f$ is injective. Since $M$
is compact and $f$ is surjective, $Y$ is compact. Hence $\gamma_f$ is a homeomorphism
onto the connected analytic subset $\gamma_f(Y)$ of $B$.

The normalized fiber measure $\mu_y:=\mug(\gamma_f(y))$ is characterized by
\[
 \int_M\varphi\,\dd\mu_y=\frac1{V_B}\int_{[M_y]}\varphi\,\frac{\omega^q}{q!}.
\]
This includes singular and nonreduced fibers with their cycle multiplicities. The induced
transport distance on the base is
\[
 \mathsf d_f(y_0,y_1):=\inf_c L_\W(\mu\circ c)
 =\W_{\gamma_f(Y)}^{\mathrm{in}}(\gamma_f(y_0),\gamma_f(y_1)),
\]
where $c$ ranges over continuous curves in $Y$ joining the prescribed endpoints.
By \Cref{thm:inner-wasserstein-complete-geodesic}, $(Y,\mathsf d_f)$ is a complete
geodesic metric space. If a smooth curve $c$ has classifying cycles in the good locus
of an irreducible component $C\subset B$, then
\[
 L_\W(\mu\circ c)=L_{\mathrm{WP}}(\gamma_f\circ c)
\]
by \Cref{thm:isometry}. This relates the induced metric to the classical transfer of
K\"ahler structures through the classifying map of a family
\cite[Proposition~12.3.5.1 and Corollary~12.3.5.2]{cycles2}.
\subsection{An irreducible component with empty good locus}
The following component has no good points, yet its inner Wasserstein
metric is induced by an explicit smooth Riemannian metric. Thus the finiteness
theorem applies even when there is no good-locus Weil--Petersson metric to complete.

\begin{example}
\label{ex:hirzebruch-inner-geometry}
A concrete model for $M=\mathbb F_2$ is obtained by adding one point at infinity to
each fiber of the holomorphic tangent bundle $T\CP^1$. Since $T\CP^1$ is the line
bundle $\mathcal O(2)$, its fiberwise projective compactification is the second
Hirzebruch surface. The added points form the section $S$. Its normal bundle is
$(T\CP^1)^*\simeq\mathcal O(-2)$: reciprocating the fiber coordinate near infinity
gives the dual transition functions. Thus $S^2=-2$, and $S$ is the negative section.
Let $\pi:M\to\CP^1$ be the projection, whose fibers are projective lines, and put
$F_t=\pi^{-1}(t)$.  In \Cref{rem:scope-bad-locus-theorem} we proved that
\[
 \Gamma:=\{[S]+[F_t]:t\in\CP^1\}
\]
is an irreducible connected component of $\mathcal C_1(M)$ and that every one of its cycles is
reduced and singular.  The classifying map
\[
 \Phi:\CP^1\longrightarrow\Gamma,
 \qquad t\longmapsto [S]+[F_t],
\]
is an analytic isomorphism: the linear system has the fixed component $S$, and its moving part is
the pencil $|F|\simeq\CP^1$.

Fix a K\"ahler form $\omega$ on $M$ and put
\[
 a:=\int_S\omega,
 \qquad b:=\int_{F_t}\omega,
 \qquad V_\Gamma=a+b.
\]
The number $b$ is independent of $t$, and the cycle-to-measure map is
\begin{equation}\label{eq:hirzebruch-cycle-measure}
 \mug([S]+[F_t])
 =\frac1{V_\Gamma}
 \left(\omega\lfloor_S+\omega\lfloor_{F_t}\right).
\end{equation}
The first summand is fixed.  To describe the moving summand, let
$\mathcal H:=(\ker d\pi)^\perp$ be the $g$-orthogonal horizontal distribution.  Since $\pi$ is
holomorphic, both $\ker d\pi$ and $\mathcal H$ are complex line bundles.  For
$v\in T_t\CP^1$, denote by $v^H(x)\in\mathcal H_x$, $x\in F_t$, its horizontal lift.

Define a Riemannian metric on the base by
\begin{equation}\label{eq:hirzebruch-base-metric}
 g_\Gamma(v,w)
 :=\frac1{V_\Gamma}\int_{F_t}g(v^H,w^H)\,\omega.
\end{equation}
This is smooth and positive definite.  If $t=t(s)$ is a smooth base curve, the section $S$ has
zero normal velocity and the normal velocity of $F_{t(s)}$ is $\dot t(s)^H$.  The two components
meet at one point, which has zero area for the cycle measure. We now justify the speed formula
directly, since the support itself is not a smooth submanifold.

Define a vector field $u_s$ along $S\cup F_{t(s)}$, up to a set of measure zero, by
\[
 u_s=0\quad\text{on }S,
 \qquad
 u_s=\dot t(s)^H\quad\text{on }F_{t(s)}.
\]
The fixed component contributes no derivative. Applying Cartan's formula and the K\"ahler
contraction calculation from the proof of \Cref{thm:isometry} to the moving fiber gives, for
every $\varphi\in C^\infty(M)$,
\begin{equation}\label{eq:hirzebruch-continuity-equation}
 \frac{\dd}{\dd s}\int_M\varphi\,\dd\mug([S]+[F_{t(s)}])
 =\int_M\langle\nabla\varphi,u_s\rangle\,
   \dd\mug([S]+[F_{t(s)}]).
\end{equation}

The field $u_s$ belongs to the Wasserstein tangent space of this singularly supported measure.
Indeed, let $p=S\cap F_{t(s)}$. Choose cutoffs on $F_{t(s)}$ which vanish near $p$ and converge
pointwise to $1$ away from $p$. On the part of the fiber where such a cutoff is supported, the
tubular-potential construction produces a smooth function on $M$ whose gradient along the fiber
is the cutoff times $\dot t(s)^H$ and whose support is disjoint from $S$. The gradients therefore
converge to $u_s$ in $L^2(\mug([S]+[F_{t(s)}]))$, because the omitted fiber neighborhoods shrink to
the zero-area set $\{p\}$. Thus $u_s$ lies in the $L^2$-closure of smooth gradients. As in the
minimal-velocity argument in the proof of \Cref{thm:isometry}, the continuity equation
\eqref{eq:hirzebruch-continuity-equation} then identifies $u_s$ with the minimal velocity.

The speed identity also holds at every parameter of a smooth base curve. To see this
without imposing smoothness on the union $S\cup F_{t(s)}$, fix $s$ and
$f\in C^\infty(M)$, and write $\mu_r:=\mug([S]+[F_{t(r)}])$.
Taylor's formula along minimizing ambient geodesics, integrated
against an optimal coupling, gives
\[
 \left|\int_M f\,\dd\mu_{s+h}-\int_M f\,\dd\mu_s\right|
 \leq \|\nabla f\|_{L^2(\mu_s)}\W(\mu_{s+h},\mu_s)
       +C_f\W(\mu_{s+h},\mu_s)^2.
\]
The continuity equation gives
$\W(\mu_{s+h},\mu_s)\leq\int_s^{s+h}\|u_r\|_{L^2(\mu_r)}\,\dd r$
for $h>0$, and the analogous bound for $h<0$. The integrand is continuous.
Dividing the preceding Taylor estimate by $|h|$, using the first-variation formula,
and approximating $u_s$ by smooth gradients proves the matching lower bound
$\liminf_{h\to0}\W(\mu_{s+h},\mu_s)/|h|\geq\|u_s\|_{L^2(\mu_s)}$.

Consequently
\begin{equation}\label{eq:hirzebruch-speed}
 \left|\frac{\dd}{\dd s}\mug([S]+[F_{t(s)}])\right|_{\W}^2
 =\frac1{V_\Gamma}\int_{F_{t(s)}}|\dot t(s)^H|_g^2\,\omega
 =g_\Gamma(\dot t(s),\dot t(s))
\end{equation}
for every $s$.

To pass from smooth curves to the full inner distance, we verify that the
cycle-to-measure map is locally bi-Lipschitz for the metric $g_\Gamma$.

Fix $t_0\in\CP^1$ and a real basis $e_1,e_2$ of $T_{t_0}\CP^1$.  On the two components of
$[S]+[F_{t_0}]$, define the piecewise normal fields
\[
 u_i=0\quad\hbox{on }S,
 \qquad
 u_i=e_i^H\quad\hbox{on }F_{t_0}.
\]
Their values at $S\cap F_{t_0}$ are irrelevant for the cycle measure.  The Gram matrix
\[
 G_{ij}:=\frac1{V_\Gamma}\int_{F_{t_0}}g(e_i^H,e_j^H)\,\omega
       =g_\Gamma(e_i,e_j)
\]
is positive definite. The cutoff and tubular-potential argument just given provides smooth
functions $f_1,f_2\in C^\infty(M)$ whose
gradients, restricted to $S\cup F_{t_0}$, are arbitrarily close to $u_1,u_2$ in
$L^2(\mug([S]+[F_{t_0}]))$.  Choose them close enough that the matrix
\[
 M_{ij}:=\frac1{V_\Gamma}
 \int_{F_{t_0}}g(\nabla f_i,e_j^H)\,\omega
\]
is invertible, and set
\[
 \Psi(t):=\left(
  \int_M f_1\,\dd\mug([S]+[F_t]),
  \int_M f_2\,\dd\mug([S]+[F_t])
 \right).
\]
The first-variation formula \eqref{eq:hirzebruch-continuity-equation} gives
$D\Psi_{t_0}(e_j)=(M_{1j},M_{2j})$.  Thus $\Psi$ is a local coordinate map at $t_0$.
Kantorovich duality for $W_1$, followed by $W_1\leq\W$, now gives on a sufficiently small
neighborhood $U$ of $t_0$ a constant $c_U>0$ such that
\[
 c_U d_{g_\Gamma}(s,t)
 \leq \W\bigl(\mug([S]+[F_s]),\mug([S]+[F_t])\bigr)
 \qquad(s,t\in U).
\]
The reverse local estimate follows by joining $s$ to $t$ by a short $g_\Gamma$-geodesic in
$U$ and using \eqref{eq:hirzebruch-speed}.  Hence $\mug\circ\Phi$ is locally bi-Lipschitz from
$(\CP^1,d_{g_\Gamma})$ to its Wasserstein image.

It follows that a curve in the image is Wasserstein-rectifiable exactly when its pullback to
$\CP^1$ is Riemannian-rectifiable.  
After arclength parametrization the pullback is absolutely continuous. At almost every
parameter it is differentiable in a local coordinate chart. Compare it there with a
smooth curve having the same value and derivative: their coordinate difference is
$o(|h|)$, and the local Lipschitz bound makes the Wasserstein difference $o(|h|)$ as
well. Formula~\eqref{eq:hirzebruch-speed} therefore identifies the two metric speeds
almost everywhere, and consequently the two lengths. Taking infima gives
\begin{equation}\label{eq:hirzebruch-inner-distance}
 \W_\Gamma^{\mathrm{in}}([S]+[F_s],[S]+[F_t])
 =d_{g_\Gamma}(s,t)
\end{equation}
for all $s,t\in\CP^1$.

The metric has a useful intrinsic expression.  Fiber integration gives the positive
$(1,1)$-form
\begin{equation}\label{eq:hirzebruch-base-form}
 \omega_\Gamma
 :=\frac1{V_\Gamma}\,\pi_*\!\left(\frac{\omega^2}{2}\right),
\end{equation}
and the horizontal--vertical splitting shows that
$g_\Gamma(v,w)=\omega_\Gamma(v,Jw)$.  Thus \eqref{eq:hirzebruch-base-form} is the K\"ahler form
of the computed inner geometry.

Its total area is especially simple.  Using
$S^2=-2$, $S\cdot F=1$, and $F^2=0$, the cohomology class of $\omega$ is
\[
 [\omega]=b\,\operatorname{PD}(S)+(a+2b)\,\operatorname{PD}(F).
\]
Consequently
\[
 \int_M\frac{\omega^2}{2}=b(a+b)=bV_\Gamma,
 \qquad
 \int_{\CP^1}\omega_\Gamma=b.
\]
For a general K\"ahler form, \eqref{eq:hirzebruch-base-form} is the most explicit possible
answer: the inner metric need not be round.  If $\omega$ is invariant under the standard lifted
$SU(2)$-action and $\omega_{\mathrm{FS}}$ is normalized by
$\int_{\CP^1}\omega_{\mathrm{FS}}=1$, then
\[
 \omega_\Gamma=b\,\omega_{\mathrm{FS}},
 \qquad
 \W_\Gamma^{\mathrm{in}}([S]+[F_s],[S]+[F_t])
 =\sqrt b\,d_{\mathrm{FS}}(s,t).
\]
This provides an explicit inner geometry on a Barlet component which is entirely contained in the
bad locus.
\end{example}

\subsection{Projective examples}

\begin{example}\label{ex:hypersurfaces-cpn}
Let $M=\CP^n$ carry the Fubini--Study metric, and let $B\subset\mathcal C_{n-1}(M)$
parametrize hypersurfaces of degree $d\geq1$. This is the setting of \cite{ACL}.
Writing $H_{n,d}=H^0(\CP^n,\mathcal O(d))$, the correspondence between nonzero
sections up to scale and effective divisors gives
\[
 B\simeq\mathbb P(H_{n,d})\simeq\CP^N,
 \qquad N=\binom{n+d}{d}-1.
\]
This identification is holomorphic, since the universal zero divisor is an analytic
family; see also \Cref{remark:projective}. In particular, $B$ is smooth and
irreducible. Its good locus is the complement of the discriminant
$\Sigma_{n,d}$, consisting of polynomials with nonsmooth or nonreduced zero divisor.
For $d\geq2$ this is a proper algebraic hypersurface; for $d=1$ it is empty.

All degree-$d$ cycles have volume
$V_B=d\,\operatorname{vol}(\CP^{n-1})$, and
\[
 \int_{\CP^n}f\,\dd\mug([Z(p)])
 =\frac1{V_B}\int_{[Z(p)]}f\,\frac{\omega_{\mathrm{FS}}^{n-1}}{(n-1)!}.
\]
Thus the measure used in \cite{ACL} is precisely the measure of
\Cref{def:mu-gamma}, including at singular and nonreduced hypersurfaces.

The normal deformation field can be expressed without choosing a local defining
function. Fix a nonzero representative $p$ of a good point $[p]$. A tangent vector
$v\in T_{[p]}\mathbb P(H_{n,d})$
is a linear map $\mathbb Cp\to H_{n,d}/\mathbb Cp$; choosing $q\in H_{n,d}$
representing $v(p)$ represents $v$ by $[p+tq]$ to first order.
Along $Z(p)$, differentiation of the section $p$ defines a surjective
complex-linear bundle map
\[
 D_p:T\CP^n|_{Z(p)}\longrightarrow\mathcal O(d)|_{Z(p)},
 \qquad \ker D_p=TZ(p).
\]
This derivative is intrinsic along the zero set: changing the connection changes
the covariant derivative by a term proportional to $p$, which vanishes there.
Give $\mathcal O(d)$ its standard Hermitian metric and write $D_p^*$ for the
Hermitian adjoint. Differentiating $p_t=0$ gives $D_p\nu_q=-q$, whose orthogonal
normal solution is
\[
 \nu_q=-\frac{D_p^*q}{\|D_p\|^2},
 \qquad
 |\nu_q|_g^2=\frac{|q|^2}{\|D_p\|^2}.
\]
Here $\|D_p\|$ is the operator norm; since the target is a line,
$D_pD_p^*=\|D_p\|^2\operatorname{id}$. Replacing $q$ by $q+ap$ does not change
its restriction to $Z(p)$. Hence the formula depends only on $v$, as required.
In particular, the normalized Weil--Petersson norm is
\[
 \|v\|_{\mathrm{WP}}^2
 =\frac1{V_B}\int_{Z(p)}
       \frac{|q|^2}{\|D_p\|^2}\,
       \frac{\omega_{\mathrm{FS}}^{n-1}}{(n-1)!}.
\]
This recovers the energy formula of \cite[Theorem~1.3]{ACL}. The Hermitian form
written in \cite[equation~(1.10)]{ACL} is unnormalized: its real part is
$V_B\Theta$ in the present notation.

If $\pi_1:\mathcal V\to B$ and $\pi_2:\mathcal V\to\CP^n$ are the projections
of the universal hypersurface, then, on the good locus,
\[
 \owp=\frac1{V_B}(\pi_1)_*
       \left(\pi_2^*\frac{\omega_{\mathrm{FS}}^n}{n!}\right).
\]
This is the normalized fiber-integral identity of \cite[Theorem~1.6]{ACL}.
Our general results give a finite, complete geodesic inner distance on all of $B$,
locally equal to the normalized Weil--Petersson distance near every good point.
\end{example}

\begin{example}\label{ex:grassmannian-linear-spaces}
Let \(M=\CP^n\) carry the Fubini--Study form \(\omega_{\mathrm{FS}}\),
normalized by $\int_{\CP^1}\omega_{\mathrm{FS}}=\pi$, and fix
\(q\in\{0,\dots,n-1\}\). Let
\[
B\subset \mathcal C_q(\CP^n)
\]
be the connected component parametrizing projective linear subspaces of dimension \(q\),
equivalently the component of \(\mathcal C_q(\CP^n)\) consisting of cycles of degree \(1\). Then
\[
B \simeq \operatorname{Gr}(q+1,n+1),
\]
the Grassmannian of \((q+1)\)-dimensional complex vector subspaces of \(\C^{n+1}\). More
precisely, a point \(W\in \operatorname{Gr}(q+1,n+1)\) corresponds to the linear cycle
\[
[L_W],\qquad L_W:=\mathbb P(W)\subset \CP^n.
\]

Here there is no discriminant: every
point of \(B\) represents a smooth reduced cycle, and \(B\) itself is smooth. Hence
$\widetilde B=B.$
In particular, the intrinsic transport geometry on this component is globally Riemannian, with no
singular or nonreduced boundary to add. Since all \(q\)-planes are projectively equivalent, the constant \(V_B\) is simply the
volume of \(\CP^q\):
\[
V_B=\vol_{\omega_{\mathrm{FS}}}(\CP^q).
\]
The cycle-to-measure map sends \([L_W]\) to its normalized volume measure
\[
\mug([L_W])=\frac{1}{\vol(\CP^q)}\,\frac{\omega_{\mathrm{FS}}^q}{q!}\bigg|_{L_W}.
\]

Fix the standard Hermitian structure on \(\C^{n+1}\). For
\(W\in \operatorname{Gr}(q+1,n+1)\), write \(\C^{n+1}=W\oplus W^\perp\). Then
\[
T_W\operatorname{Gr}(q+1,n+1)\simeq \operatorname{Hom}(W,W^\perp).
\]
If \(A\in \operatorname{Hom}(W,W^\perp)\), the corresponding first-order deformation of
\(L_W=\mathbb P(W)\) has normal field
\[
\nu_A([z])=\left.\frac{\dd}{\dd t}[z+tAz]\right|_{t=0}
 \qquad ([z]\in\mathbb P(W)).
\]
Since $Az\perp W$, this velocity is normal to $L_W$; the projective
expression is independent of the representative $z$. Its pointwise norm is
\[
|\nu_A([z])|^2=\frac{\|Az\|^2}{\|z\|^2}.
\]
Therefore the normalized Weil--Petersson metric is
\[
\|A\|_{\mathrm{WP}}^2
=
\frac{1}{\vol(\CP^q)}
\int_{\mathbb P(W)}
\frac{\|Az\|^2}{\|z\|^2}\,\frac{\omega_{\mathrm{FS}}^q}{q!}.
\]
Averaging over projective space gives
\[
\|A\|_{\mathrm{WP}}^2=\frac{1}{q+1}\operatorname{tr}(A^*A).
\]
Indeed, the average of the rank-one projector $zz^*/\|z\|^2$ over
$\mathbb P(W)$ commutes with every unitary endomorphism of $W$, so it is a scalar
multiple of the identity. Its trace is $1$, hence this average is
$\operatorname{id}_W/(q+1)$. Pairing with $A^*A$ gives the formula.
Thus the normalized Weil--Petersson metric is exactly the standard
\(U(n+1)\)-invariant Hermitian metric on the Grassmannian, scaled by \(1/(q+1)\).

The pointwise formula shows that the maximal normal speed is the operator norm of $A$,
whereas its transport energy is the sum of the squared singular values divided by $q+1$.
When both $W$ and $W^\perp$ have dimension at least two, maps with the same operator norm
can have different values of $\operatorname{tr}(A^*A)$. Thus the energy distinguishes
motion in one direction from simultaneous motion in several independent directions,
even when their maximal pointwise speeds agree.
For lines in $\CP^3$, singular values $(1,0)$ and $(1,1)$ give maximal speed $1$
in both cases, but squared Weil--Petersson norms $1/2$ and $1$, respectively.

Since \(\widetilde B=B\), \Cref{prop:local-wp-inner-wasserstein} identifies
the inner Wasserstein distance globally with the Riemannian distance of this metric:
\[
\W_B^{\mathrm{in}}=\mathsf d_{\operatorname{Gr}}.
\]

Equivalently, if \(W_0,W_1\subset \C^{n+1}\) are two \((q+1)\)-dimensional subspaces and
\(\theta_1,\dots,\theta_{q+1}\in[0,\pi/2]\) are their principal angles, then
\[
\W_B^{\mathrm{in}}([L_{W_0}],[L_{W_1}])^2
=
\frac{1}{q+1}\sum_{i=1}^{q+1}\theta_i^2.
\]
So, in this homogeneous component, the transport cost is governed exactly by the principal-angle
geometry of the corresponding linear subspaces.
\end{example}

\section{The case of zero-dimensional cycles}\label{sec:zerocycles}

\begin{theorem}\label{thm:zero-cycles-completion}
Let $B\subset\mathcal C_0(M)$ be a connected component of effective zero-cycles of degree $d$,
and assume that $\widetilde B\neq\varnothing$. Then the
metric completion of the good locus $(\widetilde B,\mathsf d_{\mathrm{WP}})$ is canonically
isometric to
\[
 (B,\W_B^{\mathrm{in}})=(B,\W_B).
\]
More precisely,
\[
 \mathsf d_{\mathrm{WP}}(\gamma_0,\gamma_1)
 =\W_B(\gamma_0,\gamma_1)
 =\W_B^{\mathrm{in}}(\gamma_0,\gamma_1)
 \qquad(\gamma_0,\gamma_1\in\widetilde B).
\]
\end{theorem}

The full degree-$d$ locus is the symmetric product
\[
 \operatorname{Sym}^d(M)=M^d/\mathfrak S_d.
\]
Since $M$ is connected, $B=\operatorname{Sym}^d(M)$ is irreducible. A point of $B$ is represented by an unordered tuple
$\gamma=[x_1]+\cdots+[x_d]$, repetitions included, and
\begin{equation}\label{eq:zero-cycle-empirical-measure}
 \mug(\gamma)=\frac1d\sum_{i=1}^d\delta_{x_i}.
\end{equation}
The good locus consists precisely of the configurations with pairwise distinct points:
\[
 \widetilde B
 =
 \{[x_1]+\cdots+[x_d]:x_i\neq x_j\text{ whenever }i\neq j\}.
\]

We first record the global metric computation, where collisions are allowed.

\begin{proposition}
\label{thm:zero-cycles-inner-equals-outer}
For all $\gamma_0,\gamma_1\in B$,
\begin{equation}\label{eq:zero-cycles-matching}
 \W_B(\gamma_0,\gamma_1)^2
 =
 \frac1d\min_{\sigma\in\mathfrak S_d}
 \sum_{i=1}^d \dM(x_i,y_{\sigma(i)})^2,
 \qquad
 \W_B^{\mathrm{in}}(\gamma_0,\gamma_1)=\W_B(\gamma_0,\gamma_1),
\end{equation}
where $\gamma_0=[x_1]+\cdots+[x_d]$ and $\gamma_1=[y_1]+\cdots+[y_d]$.
\end{proposition}

\begin{proof}
A coupling between the two empirical measures is represented by a doubly stochastic transport
matrix. Its quadratic cost is linear in that matrix, so the Birkhoff--von Neumann theorem shows
that a minimizer is a permutation matrix. This proves the matching formula.

Choose an optimal permutation and, for every $i$, a constant-speed minimizing geodesic $c_i$
from $x_i$ to $y_{\sigma(i)}$. The curve
\[
 s\longmapsto\frac1d\sum_{i=1}^d\delta_{c_i(s)}
\]
is the Wasserstein geodesic induced by this optimal coupling. At each time it equals the measure
of the zero-cycle \(\sum_i[c_i(s)]\in B\); collisions merely produce multiplicities and are allowed
in $B$. Thus any two points of $\mug(B)$ can be joined by an ambient Wasserstein
geodesic lying in $\mug(B)$. This geodesic is admissible for the inner distance,
which gives $\W_B^{\mathrm{in}}\leq\W_B$. The reverse
inequality holds for every inner distance.
\end{proof}

\begin{proof}[Proof of \Cref{thm:zero-cycles-completion}]
Let $\gamma_0=[x_1]+\cdots+[x_d]$ and $\gamma_1=[y_1]+\cdots+[y_d]$ lie in $\widetilde B$, choose an
optimal matching, and let $c_i$ be the minimizing geodesics used above. We claim that the
interpolating configuration remains reduced for every $0<s<1$. This is the
Monge--Mather shortening principle; see \cite[Chapter~8]{villani}. We include its short
argument in the present finite setting.

Suppose that $c_i(s)=c_j(s)$ for some $i\neq j$ and $0<s<1$, and write
$L_i=\dM(x_i,y_{\sigma(i)})$ and $L_j=\dM(x_j,y_{\sigma(j)})$. Optimality of the matching,
followed by the triangle inequality through the common point, gives
\begin{align*}
 L_i^2+L_j^2
 &\leq \dM(x_i,y_{\sigma(j)})^2+\dM(x_j,y_{\sigma(i)})^2\\
 &\leq\bigl(sL_i+(1-s)L_j\bigr)^2
      +\bigl(sL_j+(1-s)L_i\bigr)^2\\
 &=L_i^2+L_j^2-2s(1-s)(L_i-L_j)^2.
\end{align*}
Every inequality is therefore an equality and $L_i=L_j$. In particular, the two curves obtained
by exchanging the tails of $c_i$ and $c_j$ at their common point are again minimizing geodesics.
A minimizing broken geodesic is smooth at the joining point, so
$\dot c_i(s)=\dot c_j(s)$. Uniqueness for the geodesic equation then gives $c_i=c_j$, contrary
to $x_i\neq x_j$. Thus no collision occurs.

The optimal Wasserstein geodesic between $\mug(\gamma_0)$ and $\mug(\gamma_1)$ consequently lies
in $\mug(\widetilde B)$. Hence, on $\widetilde B$,
\[
 \W_{\widetilde B}^{\mathrm{in}}=\W_B.
\]
The symmetric product $B$ is irreducible. Since its good locus is nonempty,
\Cref{prop:local-wp-inner-wasserstein} applies and gives
\[
 \mathsf d_{\mathrm{WP}}
 =\W_{\widetilde B}^{\mathrm{in}}
 =\W_B
 =\W_B^{\mathrm{in}}
 \quad\text{on }\widetilde B\times\widetilde B,
\]
where the last equality is \Cref{thm:zero-cycles-inner-equals-outer}.

Finally, $B$ is compact and $\W_B$ induces its Barlet topology by \Cref{prop:mu-gamma}; hence
$(B,\W_B)=(B,\W_B^{\mathrm{in}})$ is complete. If $\dim_\C M>0$, coincident points can be separated by arbitrarily small
perturbations in coordinate neighborhoods, so $\widetilde B$ is dense in $B$.
If $\dim_\C M=0$, connectedness makes $M$ a point, and the assumed nonempty
good locus forces $d=1$; the assertion is then immediate. The preceding
identity of the restricted metric with $\mathsf d_{\mathrm{WP}}$ proves the asserted completion.
\end{proof}

\subsection{A Euclidean realization}
The matching formula also gives a finite-dimensional metric model for the whole component,
including collisions. We record this consequence without repeating the full construction
of the Euclidean embedding for $Q$-valued points.

\begin{theorem}\label{thm:zero-cycles-bilipschitz-embedding}
Let $B=\operatorname{Sym}^d(M)$ be the component of effective zero-cycles of degree $d\geq1$.
There are an integer $N$ and a definable bi-Lipschitz embedding
\[
 (B,\W_B^{\mathrm{in}})\longrightarrow\R^N.
\]
\end{theorem}
\begin{proof}
Choose the real-analytic embedding $\iota_M:M\hookrightarrow\R^r$ of
\Cref{prop:definable-tautological-diagram}. It is bi-Lipschitz for the Riemannian distance
on the compact manifold $M$: boundedness of its differential gives the upper bound,
the inverse function theorem gives a uniform lower bound near the diagonal, and
compactness and injectivity give the lower bound away from the diagonal.

Let $\mathcal A_d(\R^r)$ be the space of unordered $d$-tuples, written as formal sums
$T=\sum_i\llbracket x_i\rrbracket$, with matching metric
\[
 \mathcal G(T,S)^2:=\min_{\sigma\in\mathfrak S_d}\sum_{i=1}^d|x_i-y_{\sigma(i)}|^2.
\]
The embedding theorem of Almgren, in the form of
\cite[Theorem~2.1]{delellis-spadaro}, supplies a bi-Lipschitz map
$\xi:\mathcal A_d(\R^r)\to\R^N$. Its lift to ordered tuples is semialgebraic:
the construction concatenates the ordered lists of projections onto finitely many fixed
directions \cite[Section~2.1.2]{delellis-spadaro}, and sorting is a piecewise linear
semialgebraic map.

Define
\[
 \Phi\left(\sum_{i=1}^d[x_i]\right)
 :=d^{-1/2}\xi\left(\sum_{i=1}^d\llbracket\iota_M(x_i)\rrbracket\right).
\]
The matching formula and the bi-Lipschitz bounds for $\iota_M$ and $\xi$ show that
$\Phi$ is bi-Lipschitz for $\W_B$. By \Cref{thm:zero-cycles-inner-equals-outer},
this is also $\W_B^{\mathrm{in}}$.

For completeness, definability on the given Barlet component can be checked without
choosing an ordering. Let $Q:M^d\to B$ be the quotient map and let
$\widetilde\Phi:M^d\to\R^N$ be the displayed symmetric lift. The map $Q$ is
real analytic between compact real-analytic spaces, hence definable in the realizations
of \Cref{prop:definable-tautological-diagram}; so is $\widetilde\Phi$.
The graph of $\Phi$ is the projection of
\[
 \{(\gamma,z,x)\in B\times\R^N\times M^d:Q(x)=\gamma,\ \widetilde\Phi(x)=z\},
\]
and is therefore definable.
\end{proof}

\begin{remark}
The metric argument uses only the compact Riemannian manifold and its matching distance.
The real-analytic structure makes the embedding definable; the K\"ahler condition is
not needed for this particular result.
\end{remark}
\subsection{Examples}

\begin{example}\label{ex:ruled-surface-component}
Let
\[
M:=\CP^1_u\times \CP^1_v
\]
endowed with the product Fubini--Study form
\[
\omega:=\operatorname{pr}_1^*\omega_{\mathrm{FS}}+\operatorname{pr}_2^*\omega_{\mathrm{FS}}.
\]
For \(u\in \CP^1\), set
\[
F_u:=\{u\}\times \CP^1_v\subset M.
\]
Fix \(d\geq 1\), and let \(B\subset \mathcal C_1(M)\) be the component parametrizing effective
divisors of bidegree \((d,0)\), equivalently divisors in the linear system
\(|\mathcal O_{\CP^1\times\CP^1}(d,0)|\).

Every such divisor is uniquely of the form
\[
\gamma=\sum_{i=1}^\ell m_i[F_{u_i}],
\qquad
u_i\in \CP^1 \text{ distinct},\quad \sum_{i=1}^\ell m_i=d.
\]
Therefore \(B\) is naturally identified with the symmetric product
\[
\operatorname{Sym}^d(\CP^1)\simeq \CP^d.
\]
In particular, \(B\) is smooth and irreducible. Its good locus \(\widetilde B\) consists exactly
of the reduced divisors
\[
[F_{u_1}]+\cdots+[F_{u_d}],
\qquad u_i\neq u_j \text{ for } i\neq j.
\]

This example is worth isolating because the bad locus has a particularly simple meaning:
it comes only from multiplicities. Indeed, distinct fibers of the ruling are disjoint, so every
reduced cycle in this component is already a smooth compact complex submanifold of \(M\),
possibly disconnected. Thus the complement \(B\setminus \widetilde B\) is precisely the collision
locus where two or more fibers merge into a multiple fiber.

Let \(\nu\) denote the normalized Fubini--Study probability measure on \(\CP^1_v\). If
\[
\eta_\gamma:=\frac1d\sum_{i=1}^\ell m_i\,\delta_{u_i}\in \mathcal P(\CP^1_u),
\]
then the cycle-to-measure map is
\[
\mug(\gamma)=\eta_\gamma\otimes \nu \in \cP(M).
\]
Indeed, each fiber \(F_u\) has the same \(\omega\)-volume, and the restriction of \(\omega\) to
\(F_u\) is independent of \(u\).

Projection to the first factor is $1$-Lipschitz, so
$\W(\eta_{\gamma_0},\eta_{\gamma_1})
\leq\W(\mug(\gamma_0),\mug(\gamma_1))$.
Conversely, combine an optimal coupling on the first factor with the coupling
$(\operatorname{id},\operatorname{id})_\#\nu$ on the second. Since the squared
product distance is the sum of the squared factor distances, this coupling has
cost $\W(\eta_{\gamma_0},\eta_{\gamma_1})^2$. Hence
\[
 \W\bigl(\mug(\gamma_0),\mug(\gamma_1)\bigr)
 =\W(\eta_{\gamma_0},\eta_{\gamma_1}).
\]

It follows that the intrinsic Wasserstein geometry of the curve component \(B\) is exactly the
degree-\(d\) zero-cycle geometry of \(\operatorname{Sym}^d(\CP^1)\). In particular,
\[
\W_B^{\mathrm{in}}(\gamma_0,\gamma_1)
=
\W_B(\gamma_0,\gamma_1)
=
\W(\eta_{\gamma_0},\eta_{\gamma_1}),
\]
and \((B,\W_B^{\mathrm{in}})\) is the metric completion of \(\widetilde B\) endowed with
its normalized Weil--Petersson metric.

Here cycles of positive dimension have exactly the same transport geometry as
unordered points on the base. The completion adds the nonreduced fiber cycles.
\end{example}
\begin{example}\label{ex:C-zero-cycles}
Although \(\mathbb C\) is noncompact, the constructions discussed above become completely explicit in
this case, and provide a useful model for the compact situation.

Let \(B\) be the space of effective zero-cycles of degree \(d\) on \(\mathbb C\). Thus a point
of \(B\) is a formal sum \(\gamma=[\lambda_1]+\cdots+[\lambda_d]\), where the \(\lambda_i\in\mathbb
C\) are allowed to coincide and are unordered. Equivalently, \(B\) is the symmetric product
\(\operatorname{Sym}^d(\mathbb C)\).

A convenient set of coordinates on \(B\) is given by the coefficients of the corresponding
monic polynomial. Namely, to \(\gamma=[\lambda_1]+\cdots+[\lambda_d]\) one associates
\[
p_\gamma(z):=\prod_{i=1}^d (z-\lambda_i)
= z^d+c_1z^{d-1}+\cdots+c_d.
\]
The coefficients \((c_1,\dots,c_d)\in\mathbb C^d\) are the elementary symmetric functions of the
roots, so this gives an identification \(B\simeq \mathbb C^d\). Under this identification, the
good locus \(\widetilde{B}\) is the complement of the discriminant hypersurface, that is, the
subset of coefficient vectors for which \(p_\gamma\) has \(d\) distinct roots.

On the ordered cover \(\mathbb C^d\setminus\Delta\), where
\(\Delta=\{(\lambda_1,\dots,\lambda_d):\lambda_i=\lambda_j\text{ for some }i\neq j\}\), the quotient
map is precisely the Vi\`ete map sending the ordered roots to the coefficients
\((c_1,\dots,c_d)\). The cycle-to-measure map is simply
\[
\mug(\gamma)=\frac1d\sum_{i=1}^d \delta_{\lambda_i}.
\]

The tangent directions are equally concrete. At a reduced cycle
\(\gamma=[\lambda_1]+\cdots+[\lambda_d]\in\widetilde{B}\), a tangent vector is represented upstairs
by a tuple \((\dot\lambda_1,\dots,\dot\lambda_d)\in\mathbb C^d\), modulo permutation of the labels.
Since the ambient metric is Euclidean and the cycle is zero-dimensional, the normalized
Weil--Petersson norm is just the averaged Euclidean norm
\[
\|\dot\gamma\|_{\mathrm{WP}}^2=\frac1d\sum_{i=1}^d |\dot\lambda_i|^2.
\]
Thus the Weil--Petersson metric on \(\widetilde{B}\) is the quotient of the flat metric
\(\frac1d\sum_{i=1}^d |d\lambda_i|^2\) on \(\mathbb C^d\setminus\Delta\) by the action of
\(\mathfrak S_d\).

In coefficient coordinates, the metric is therefore the pushforward of this flat metric through the
Vi\`ete map. More explicitly, if \(c_k=(-1)^k e_k(\lambda_1,\dots,\lambda_d)\), then
\(\dot c_k=\sum_{i=1}^d (\partial c_k/\partial\lambda_i)\dot\lambda_i\), so the Weil--Petersson
metric in the \(c\)-coordinates is obtained from the Jacobian of the Vi\`ete map. It is smooth on $\widetilde B$, while its coefficients become singular
at collisions. The metric need not tend to zero there: for example, along
$p_a(z)=z^2-a$, $a\neq0$, one obtains
\[
 \|\dot a\|_{\mathrm{WP}}^2=\frac{|\dot a|^2}{4|a|}.
\]
Although this coefficient diverges as $a\to0$, the radial path to the double root
has finite length.

Finally, for two cycles \(\gamma_0=[\lambda_1]+\cdots+[\lambda_d]\) and
\(\gamma_1=[\mu_1]+\cdots+[\mu_d]\), the ambient Wasserstein distance between the associated empirical
measures is given by the optimal matching formula
\[
\W^2\bigl(\mug(\gamma_0),\mug(\gamma_1)\bigr)
=
\frac1d\min_{\sigma\in\mathfrak S_d}\sum_{i=1}^d |\lambda_i-\mu_{\sigma(i)}|^2.
\]
The corresponding inner distance is realized by optimally matching the roots and
moving each root along a straight segment. The shortening argument in the
proof of \Cref{thm:zero-cycles-completion} applies unchanged to $\C$: when the
endpoint roots are distinct, the interpolating roots do not collide. Thus the
Weil--Petersson completion is the full Euclidean symmetric product, with its
matching distance.
\end{example}

\end{document}